\documentclass[11pt,reqno]{amsart}
\usepackage[a4paper]{anysize}\marginsize{2cm}{2cm}{2cm}{2cm}
\pdfpagewidth=\paperwidth \pdfpageheight=\paperheight
\allowdisplaybreaks
\usepackage[english]{babel}
\usepackage[latin1]{inputenc}
\usepackage[T1]{fontenc} 
\usepackage{latexsym}
\usepackage{amssymb,amsmath,amstext,amsfonts,amsthm}
\usepackage{verbatim}
\usepackage{xcolor}
\usepackage{mathrsfs,caption}
\usepackage{mathtools}
\usepackage{tikz}
\usepackage{relsize}
\usetikzlibrary{arrows}
\usepackage[colorlinks=true, urlcolor=blue, linkcolor=blue, citecolor=blue]{hyperref}

\numberwithin{equation}{section}
\theoremstyle{plain}%
\newtheorem*{thm*}{Main Theorem}
\newtheorem{theorem}{Theorem}
\numberwithin{theorem}{section}
\newtheorem{proposition}[theorem]{Proposition}
\newtheorem{lemma}[theorem]{Lemma}
\newtheorem{corollary}[theorem]{Corollary}

\theoremstyle{definition}
\newtheorem{definition}[theorem]{Definition}
\newtheorem{remark}[theorem]{Remark}
\newtheorem{example}[theorem]{Example}
\newtheorem{problem}[theorem]{Problem}

\newcommand{\C}{\mathbb{C}}

\newcommand{\PP}{\mathbb{P}}
\newcommand{\R}{\mathbb{R}}

\newcommand{\arxiv}[1]{\href{http://arxiv.org/abs/#1}{{\tt arXiv:#1}}}

\date{}
\begin{document}

\author{Luca Sodomaco}
\address{Dipartimento di Matematica e Informatica ``Ulisse Dini'', University of Florence, Italy}
\email{luca.sodomaco@unifi.it}

\subjclass[2010]{14P05, 14M20, 15A72, 15A18, 58K05}

\title{On the product of the singular values of a binary tensor}

\begin{abstract}
A real binary tensor consists of $2^d$ real entries arranged into hypercube format $2^{\times d}$. For $d=2$, a real binary tensor is a $2\times 2$ matrix with two singular values. Their product is the determinant. We generalize this formula for any $d\ge 2$. Given a partition $\mu\vdash d$ and a $\mu$-symmetric real binary tensor $t$, we study the distance function from $t$ to the variety $X_{\mu,\R}$ of $\mu$-symmetric real binary tensors of rank one. The study of the local minima of this function is related to the computation of the singular values of $t$. Denoting with $X_\mu$ the complexification of $X_{\mu,\R}$, the Euclidean Distance polynomial $\mathrm{EDpoly}_{X_\mu^\vee,t}(\epsilon^2)$ of the dual variety of $X_\mu$ at $t$ has among its roots the singular values of $t$. On one hand, the lowest coefficient of $\mathrm{EDpoly}_{X_\mu^\vee,t}(\epsilon^2)$ is the square of the $\mu$-discriminant of $t$ times a product of sum of squares polynomials. On the other hand, we describe the variety of $\mu$-symmetric binary tensors that do not admit the maximum number of singular values, counted with multiplicity. Finally, we compute symbolically all the coefficients of $\mathrm{EDpoly}_{X_\mu^\vee,t}(\epsilon^2)$ for tensors of format $2\times 2\times 2$. 
\end{abstract}

\maketitle

\section{Introduction}
The tensor space $V_\R^{\otimes d}$, with $V_\R\cong\R^2$, contains real binary tensors $t=(t_{i_1\cdots i_d})$ of format $2^{\times d}$. The binary tensors which can be written as $t=x_1\otimes\cdots\otimes x_d$ for some $x_1,\ldots,x_d\in V_\R$ are called {\em rank one} binary tensors. In coordinates, $t_{i_1\cdots i_d}=x_{1,i_1}\cdots x_{d,i_d}$ for all $1\le i_j\le 2$ and $1\le j\le d$. An excellent reference for the algebraic geometry for spaces of tensors is \cite{Lan}. The set $X_{d,\R}\subset V_\R^{\otimes d}$ of real binary tensors of rank one is the cone over the {\em Segre variety}
\[
\mathrm{Seg}(\PP(V_\R)^{\times d})\subset\PP(V_\R^{\otimes d})\cong\PP_\R^{2^d-1}.
\]
By slight abuse of notation, we use the symbol $X_{d,\R}$ also for that Segre variety. We equip $V_\R$ with a scalar product $q\colon V_\R\times V_\R\rightarrow\R$, equivalently with a positive definite quadratic form $q\colon V_\R\rightarrow\R$. In coordinates, for example we use the standard quadratic form $q(x)=x_0^2+x_1^2$ for all $x=(x_0,x_1)\in\R^2$. The quadratic form $q$ induces a quadratic form $\widetilde{q}\colon V_\R^{\otimes d}\rightarrow\R$ by the relation $\widetilde{q}(v_1\otimes\cdots\otimes v_d)=q(v_1)\cdots q(v_d)$, then extended by linearity. The square root of $\widetilde{q}$ is known as the {\em Bombieri norm} of a tensor. In coordinates,
\begin{equation}\label{eq: definition quadratic form on tensors}
\widetilde{q}(t)=\sum_{(i_1,\ldots,i_d)\in\{0,1\}^{\times d}}t_{i_1\cdots i_d}^2\quad\mathrm{for}\ t\in V_\R^{\otimes d}.
\end{equation}
Any tensor $t$ defines a {\em squared distance function} $d_t:X_{d,\R}\longrightarrow \R$ over $X_{d,\R}$ by $d_t(x)\coloneqq\widetilde{q}(t-x)$. The problem of finding the real binary tensor $x$ of rank one that minimizes the squared distance function $d_t$ is called {\em best rank one approximation problem} for $t$. Since $X_{d,\R}$ is a smooth variety, the minimum of $d_t$ is attained among the points $x\in X_{d,\R}$ such that $t-x$ is in the normal space of $X_{d,\R}$ at $x$. Any such $x\in X_{d,\R}$ is called {\em critical} binary tensor of rank one for $t$. Critical binary tensors of rank one may be computed applying a striking result by Lim \cite{Lim} and Qi \cite{Q}, here adapted to the binary case. We stress that, although the best rank one approximation problem arises in a real context, we need to consider the complexified space $V\coloneqq V_\R\otimes\C\cong\C^2$ and extend the quadratic forms $q$ and $\widetilde{q}$ to complex-valued functions, which are not Hermitian forms. In the context of Quantum Information Theory (QIT), binary tensors $t\in V^{\otimes d}$ are called {\em $d$-qubit states}. We say that a tensor $t$ is {\em isotropic} if $\widetilde{q}(t)=0$.
\begin{theorem}[Lim, Qi]\label{thm: Lim Qi}
Given a real binary tensor $t\in V_\R^{\otimes d}$, the non-isotropic critical tensors of rank one for $t$ correspond to tensors $\sigma x=\sigma (x_1\otimes\cdots\otimes x_d)\in V^{\otimes d}$ such that $q(x_j)=1$ for all $1\le j\le d$ and
\begin{equation}\label{eq: singular vector tuple}
\widetilde{q}(t, x_1\otimes\cdots\otimes x_{j-1}\otimes\_\!\otimes x_{j+1}\otimes\cdots\otimes x_d)=\sigma\cdot q(x_j,\_),\quad 1\le j\le d,
\end{equation}
for some $\sigma\in\C$, called {\em singular value} of $t$. The corresponding $d$-ple $(x_1,\ldots,x_d)$ is called {\em singular vector $d$-ple} for $t$. Moreover, we call {\em singular tensor} for $t$ any tensor of rank one written as $\sigma(x_1\otimes\cdots\otimes x_d)$, where $(x_1,\ldots,x_d)$ and $\sigma$ are a singular $d$-ple and a singular value for $t$, respectively.
\end{theorem}

Note that both sides of relation (\ref{eq: singular vector tuple}) correspond to linear operators on $V$, which are represented in coordinates by the column vectors $(y_{j,0},y_{j,1})^\mathsmaller{T}$ and $(x_{j,0},x_{j,1})^\mathsmaller{T}$ respectively, where
\[
y_{j,k} = \sum_{\substack{0\le i_l\le 1 \\ l\neq j}}t_{i_1\cdots k\cdots i_d}x_{1,i_1}\cdots x_{j-1,i_{j-1}}x_{j+1,i_{j+1}}\cdots x_{d,i_d},\quad k\in\{0,1\}.
\]
The above construction generalizes to tensors with partial symmetry. A {\em partition} of the integer $d$ is any tuple $\mu=(\mu_1,\ldots,\mu_s)$ of positive integers such that $\mu_1+\cdots+\mu_s=d$. It is indicated with $\mu\vdash d$. For any positive integer $m$, we denote by $S^mV_\R$ the $m$-th symmetric power of $V_\R$. Besides that, for any $\mu\vdash d$ we use the shorter notation $S^{\mu}V_\R$ for the tensor space $S^{\mu_1}V_\R\otimes\cdots\otimes S^{\mu_s}V_\R$. The vectors in $S^{\mu}V_\R$ are called real {\em $\mu$-symmetric} tensors, or real {\em partially symmetric} tensors if $\mu$ is not specified. The real $\mu$-symmetric binary tensors which can be written as $t=x_1^{\mu_1}\otimes\cdots\otimes x_s^{\mu_s}$ for some $x_1,\ldots,x_s\in V_\R$ are called real $\mu$-symmetric binary tensors of {\em rank one}. The set of all real $\mu$-symmetric binary tensors of rank one is the cone over the {\em Segre-Veronese variety}
\begin{equation}\label{eq: def Segre-Veronese}
X_{\mu,\R}\coloneqq\mathrm{Seg}_\mu(\PP(V_\R)^{\times s})\subset\PP(S^\mu V_\R).
\end{equation}
If $\mu=(1,\ldots,1)\eqqcolon 1^d$ is the trivial partition, we recover the Segre variety $X_{1^d,\R}=X_{d,\R}$. If $\mu=(d)$, we get the {\em Veronese variety} $X_{(d),\R}$ whose points are classes of real {\em symmetric} binary tensors of rank one in $S^dV_\R$. Any real symmetric binary tensor $t\in S^dV_\R$ corresponds to a homogeneous polynomial of degree $d$ in two indeterminates, namely a degree $d$ binary form.

For each $\mu\vdash d$, we make $S^\mu V_\R$ an Euclidean space by considering the restriction of the quadratic form $\widetilde{q}$ to $S^\mu V_\R$. It is important to underline that the distance function on $S^\mu V_\R$ induced by (the restriction of) $\widetilde{q}$ is the only one compatible with the group embedding $\mathrm{SO(V_\R)}^{\times s}\subset \mathrm{SO}(S^\mu V_\R)$.

Analogously to the non-symmetric case, given a $\mu$-symmetric binary tensor $t$ we may consider the distance function $d_t\colon X_{\mu,\R}\rightarrow \R$ and define the critical $\mu$-symmetric binary tensors of rank one for $t$. In particular, the critical $\mu$-symmetric binary tensors for $t$ are characterized by the equations (\ref{eq: singular vector tuple}), restricted from non-symmetric to $\mu$-symmetric tensors. According to Theorem \ref{thm: Lim Qi}, they are called {\em $\mu$-symmetric singular tensors} for $t$. When $\mu=(d)$, equations (\ref{eq: singular vector tuple}) simplify as
\begin{equation}\label{eq: eigenvector eigenvalue}
\widetilde{q}(t,x^{d-1}\cdot\_)=\lambda\cdot q(x,\_)
\end{equation}
and any non-zero solution $x$ of (\ref{eq: eigenvector eigenvalue}) such that $q(x)=1$ is called an {\em E-eigenvector} of $t$, while $\lambda\in\C$ is called an {\em E-eigenvalue of $t$}. Moreover, the tensor $\lambda x^d$ is called an {\em E-eigentensor of $t$} (see also \cite{CS,HHLQ,NQWW,QL}).

In this paper, we are more interested in computing the global distance from a given $\mu$-symmetric binary tensor $t$ to the Segre Veronese variety $X_{\mu,\R}$, rather than focusing on the best rank one approximation of $t$. In other words, we consider all squared distances $\epsilon^2=\widetilde{q}(t-\sigma x)$ between $t$ and a $\mu$-symmetric singular tensor $\sigma x$ for $t$. Such squared distances $\epsilon^2$ satisfy an algebraic relation like
\begin{equation}\label{eq: algebraic relation}
\sum_{j=0}^Na_j(t)\epsilon^{2j}=0,
\end{equation}
where $a_j(t)$ is a homogeneous polynomial in the entries of $t$ and $N=\mathrm{EDdegree}(X_\mu)$ is the {\em Euclidean Distance degree} of $X_\mu\coloneqq X_{\mu,\C}$, introduced by Draisma, Horobe\c{t}, Ottaviani, Sturmfels and Thomas in \cite{DHOST} for any algebraic variety in a finite dimensional Euclidean space. In particular, $\mathrm{EDdegree}(X_\mu)$ corresponds to the number of singular tensors of a general $\mu$-symmetric binary tensor. Friedland and Ottaviani computed in \cite[Theorem 12]{FO} the ED degree of any Segre-Veronese product of projective spaces. In the binary setting, their formula simplifies to
\begin{equation}\label{eq: EDdegree Segre-Veronese}
\mathrm{EDdegree}(X_\mu)=s!\mu_1\cdots\mu_s.
\end{equation}
Observe that $\mathrm{EDdegree}(X_d)=d!$, $\mathrm{EDdegree}(X_{(d)})=d$ and for $d=2$ both formulas agree with the Spectral Theorem. On one hand, up to a scalar factor the univariate polynomial at the left-hand side of (\ref{eq: algebraic relation}) is called {\em Euclidean Distance polynomial} (or {\em ED polynomial}) of $X_\mu$ at $t$ and is denoted by $\mathrm{EDpoly}_{X_\mu,t}(\epsilon^2)$. This definition has been introduced in a recent paper with Ottaviani \cite{OS} in the same setting of \cite{DHOST}. On the other hand, for any fixed $\epsilon\ge 0$, the real part of equation (\ref{eq: algebraic relation}) defines the known $\epsilon$-offset of $X_{\mu,\R}$. Offsets of real algebraic varieties find relevant engineering applications, for example in the study of geometric modeling techniques.

Here we concentrate on ED polynomials of dual varieties of Segre-Veronese varieties. As pointed out in Section \ref{sec: computation}, the roots of the ED polynomial of the dual variety of $X_\mu$ at $t\in S^\mu V$ are the squared singular values of $t$. Let $Q\coloneqq\mathcal{V}(q)=\{[(1,\sqrt{-1})],[(1,-\sqrt{-1})]\}\subset\PP(V)$ be the isotropic quadric. For $\mu=(\mu_1,\ldots,\mu_s)\vdash d$ and $J\subset[s]\coloneqq\{1,\ldots,s\}$, we define
\begin{equation}\label{eq: X mu j}
X_{\mu,J}\coloneqq\mathrm{Seg}_{\mu}(Y_1\times\cdots\times Y_s)\subset\PP(S^\mu V),
\end{equation}
where $Y_j=Q$ if $j\in J$ and $Y_j=\PP(V)$ otherwise. For all $1\le j\le s$ and all $J\subset[s]$ with $j$ elements, define $f_{\mu,J}$ to be the equation of the dual variety of $X_{\mu,J}$, when it is a hypersurface, otherwise $f_{\mu,J}\coloneqq 1$. When $\mu=1^d$, we use the notation $f_{d,J}\coloneqq f_{1^d,J}$. For example, assume $\mu=1^2$. Then $f_2\coloneqq f_{1^2}$ is the determinant of a $2\times 2$ matrix, $f_{2,\{1\}}=f_{2,\{2\}}=1$ and finally $f_{2,\{1,2\}}$ is the equation of the dual variety of $\mathrm{Seg}(Q\times Q)$, of degree four. Generally, $f_\mu\coloneqq f_{\mu,\emptyset}$ is the equation of the dual variety of $X_\mu=X_{\mu,\emptyset}$ when it is a hypersurface, usually called {\em $\mu$-discriminant} of a $\mu$-symmetric tensor. For $\mu=1^d$ the $\mu$-discriminant is known as the {\em hyperdeterminant} of a tensor, whereas for $\mu=(d)$ the $\mu$-discriminant is addressed simply as the {\em discriminant} of a symmetric tensor.

Applying the results on the lowest coefficient of ED polynomials in \cite{OS} and the inspiring work by Oeding \cite{O} on symmetrizations of the $\mu$-discriminant of a $\mu$-symmetric tensor, in this paper we determine an explicit expression, involving powers of the polynomials $f_{\mu,J}$, for the highest coefficient of the ED polynomial of the dual variety of $X_\mu$. This leads to the following closed formula for the product of the singular values of a general $\mu$-symmetric binary tensor $t\in S^\mu V$, which generalizes \cite[Main Theorem]{Sod} in the context of binary forms.

\begin{thm*}\label{thm: Main Theorem}
Consider an integer $d\ge 1$ and a partition $\mu=(\mu_1,\ldots,\mu_s)\vdash d$. If the $\mu$-symmetric binary tensor $t\in S^\mu V$ admits the maximum number $N=s!\mu_1\cdots\mu_s$ of singular values, counted with multiplicity (hypothesis verified for a general $t$), their squared product is
\begin{equation}\label{eq: the product formula}
(\sigma_1\cdots\sigma_N)^2=\prod_{J\subset[s]}f_{\mu,J}(t)^{2-\sum_{k\in J}\mu_k}.
\end{equation} 
\end{thm*}

The right-hand side of (\ref{eq: the product formula}) should be interpreted as the ratio between the lowest and the highest coefficient of the ED polynomial of the dual variety of $X_\mu$ at $t\in S^\mu V$. Depending on the sign of their exponent, the polynomials $f_{\mu,J}$ appear in the numerator or in the denominator of this ratio, otherwise they do not appear at all if their exponent is zero. If $\mu=1^d$, we get the following corollary.

\begin{corollary}\label{cor: Main Theorem, Segre case}
Consider $d\ge 1$ and $t\in V^{\otimes d}$. If $t$ admits the maximum number $d!$ of singular values, counted with multiplicity (hypothesis verified for a general $t$), their squared product is
\begin{equation}\label{eq: the product formula, Segre case}
(\sigma_1\cdots\sigma_{d!})^2=\prod_{j=0}^dg_j(t)^{2-j},\quad g_j\coloneqq\prod_{\substack{J\subset[d]\\ |J|=j}}f_{d,J}\ \mathrm{for\ all}\ 0\le j\le d.
\end{equation}
\end{corollary}

In Proposition \ref{SOS} we show that, for every non-empty subset $J\subset[s]$, $f_{\mu,J}$ is a sum of squares (SOS) polynomial. In particular, $f_{\mu,J}(t)>0$ for any non-zero real $\mu$-symmetric binary tensor $t$. This fact confirms the following known result.

\begin{proposition}{\cite[Proposition 2]{Lim}}
If the dual variety of $X_\mu$ is a hypersurface, then $0$ is a singular value of $t\in S^\mu V_\R$ if and only if $f_\mu(t)=0$.
\end{proposition}

Our paper is organized as follows. In Section \ref{sec: computation} we recall the technique of computation of the ED polynomial of an algebraic variety, and we adapt it in the case of the dual variety of the Segre-Veronese variety $X_\mu$. In Section \ref{sec: extreme coefficients} we apply the general theory of \cite{OS} for studying the polynomial $\mathrm{EDpoly}_{X_\mu^\vee,t}(0)$. Besides that, we determine the set of $\mu$-symmetric binary tensors which fail to have the maximum number of singular values, that is the vanishing locus of the highest coefficient of $\mathrm{EDpoly}_{X_\mu^\vee,t}(\epsilon^2)$. Inspired by the results of \cite{O,HO}, in Section \ref{sec: main} we develop a technique of partial symmetrizations in order to determine all the multiplicities of the factors appearing in the highest coefficient of $\mathrm{EDpoly}_{X_\mu^\vee,t}(\epsilon^2)$, thus proving the Main Theorem. Furthermore, we show some nice identities involving the degrees of $\mu$-discriminants which descend immediately from the Main Theorem. In Section \ref{sec: SOS}, we verify that $f_{\mu,J}$ is a SOS polynomial for all non-empty subsets $J$ and we provide some examples. In Section \ref{sec: 222}, we consider a general $2\times 2\times 2$ tensor $t$ and we compute simbolically all the coefficients of $\mathrm{EDpoly}_{X_d^\vee,t}(\epsilon^2)$ in terms of $\mathrm{SO}(V)^3$-invariants. This is useful for studying more in detail the $6=3!$ singular values of $t$, even when $t$ is partially symmetric. Note that in this case the formula (\ref{eq: the product formula, Segre case}) simplifies as{}
\[
(\sigma_1\cdots\sigma_6)^2=\frac{g_0^2\cdot g_1}{g_3},\quad g_0=f_3,\quad g_1=f_{3,\{1\}}\cdot f_{3,\{2\}}\cdot f_{3,\{3\}},\quad g_2=1,\quad g_3=f_{3,\{1,2,3\}}.
\]

\section{Computation of the ED polynomial of a Segre-Veronese variety}\label{sec: computation}

In this Section we show two equivalent approaches for computing the ED polynomial of the dual variety of a Segre-Veronese variety $X_\mu$ at a given binary tensor $t\in S^\mu V$. The first way follows the original setting about $\epsilon$-offsets of a variety and the ideas explained in \cite{DHOST}, and uses the Pythagorean Theorem. The second one applies directly Theorem \ref{thm: Lim Qi}.

First we set our notation. Having fixed coordinates $\{t_{i_1\cdots i_d}\}$ for the tensor space $V^{\otimes d}$, we can write equations for the subspace $S^\mu V\subset V^{\otimes d}$ and then define a new set of coordinates for $S^\mu V$. Since the ED degree and, in turn, the ED polynomial of $X_\mu$ depend strongly on the metric $\widetilde{q}$, we pay a particular attention to the polynomial defining the quadratic form $\widetilde{q}$, written in the new set of coordinates.

The partition $\mu=(\mu_1,\ldots,\mu_s)\vdash d$ and the inclusion $S^\mu V\subset V^{\otimes d}$ induce a partition of the set $[d]\coloneqq\{1,\ldots,d\}$ into mutually disjoint subsets $I_{\mu,1},\ldots, I_{\mu,s}$, meaning that we identify all copies of $V$ indexed by elements of $I_{\mu,k}$, for all $1\le k\le s$. We denote by $\Sigma_{\mu,k}$ the group of permutations of $[d]$ that involve only elements of $I_{\mu,k}$. If a binary tensor $t=(t_{i_1\ldots i_d})\in V^{\otimes d}$ is $\mu$-symmetric, then
\begin{equation}\label{eq: equations S mu V}
t_{i_1\cdots i_d}=t_{i_{\sigma_k(1)}\cdots i_{\sigma_k(d)}}\quad\mbox{for all}\ \sigma_k\in\Sigma_{\mu,k}\ \mbox{and all}\ (i_1,\ldots,i_d)\in\{0,1\}^{\times d}.
\end{equation}
The relations in (\ref{eq: equations S mu V}) are precisely the equations for the subspace $S^\mu V\subset V^{\otimes d}$. We use the smaller set of coordinates $\{c_{\omega_1\cdots\omega_s}\}$ for $S^\mu V$, where $(\omega_1,\ldots,\omega_s)$ varies in $\mathcal{P}_\mu\coloneqq\prod_{l=1}^s\{0,\ldots,\mu_l\}$. The old variable $t_{i_1\cdots i_d}$ corresponds to the new variable $c_{\omega_1\cdots\omega_s}$ if and only if $\sum_{l\in I_{\mu,k}}i_l=\omega_k$ for all $1\le k\le s$. Therefore, the number of old coordinates $t_{i_1\cdots i_d}$ which coincide on $S^\mu V$ with the new coordinate $c_{\omega_1\cdots\omega_s}$ is equal to $\prod_{k=1}^s\binom{\mu_k}{\omega_k}$. Moreover, the restriction to $S^\mu V$ of the polynomial in (\ref{eq: definition quadratic form on tensors}) defining $\widetilde{q}$ is written as
\[
\widetilde{q}(t)=\sum_{(\omega_1,\ldots,\omega_s)\in\mathcal{P}_\mu}\left[\prod_{k=1}^s\binom{\mu_k}{\omega_k}\right]c^2_{\omega_1\cdots\omega_s}\ \mathrm{for\ all}\ t=(c_{\omega_1\cdots\omega_s})\in S^\mu V.
\]

Note that for $\mu=(d)$ the above expression becomes
\begin{equation}\label{eq: metric symmetric tensors}
\widetilde{q}(t)=\sum_{j=0}^d\binom{d}{j}c_j^2\quad\mathrm{for\ all}\ t=(c_0,\ldots,c_d)\in S^d V.
\end{equation}

Since the next results hold true over an algebraically closed field, we need to consider the complex vector space $V\coloneqq V_\R\otimes\C\cong\C^2$ and the complex variety $X_\mu\coloneqq X_{\mu,\C}$. The squared distance between two complex vectors $x,y\in V$ or between two complex tensors $t,u\in S^\mu V$ is the complex value of the function $q(x-y)$ or $\widetilde{q}(t-u)$, respectively. We stress that the extensions of $q$ and $\widetilde{q}$ are not Hermitian forms. They are metrics only when restricted to $V_\R$ and $S^\mu V_\R$, respectively.

Let $\mathcal{I}_{X_\mu}\subset\C[\{c_{\omega_1\cdots\omega_s}\}]$ be the radical ideal defining the Segre-Veronese variety $X_\mu\subset S^\mu V$. The variety $X_\mu$ has codimension $c=c(\mu)\coloneqq\prod_{k=1}^s(\mu_k+1)-(s+1)$ in $S^\mu V$. We indicate by $\mathrm{Jac}_{X_{\mu}}$ the Jacobian matrix of the partial derivatives of the minimal generators of $\mathcal{I}_{X_\mu}$. For any fixed $t\in S^\mu V$, the {\em $\epsilon$-hyperball} centered at $t$ is the hypersurface
\begin{equation}\label{eq: def hyperball}
B_\mu(t;\epsilon)\coloneqq\mathcal{V}(\widetilde{q}(t-z)-\epsilon^2)\subset S^\mu V.
\end{equation}
Given $\epsilon\in\C$, a binary tensor $t=(c_{\omega_1\cdots\omega_s})\in S^\mu V$ belongs to the $\epsilon$-offset of $X_\mu$ if and only if the $\epsilon$-hyperball centered at $t$ does not intersect $X_\mu$ transversally. More precisely, there exists $z=(z_{\omega_1\cdots\omega_s})\in X_\mu\cap B_\mu(t;\epsilon)$ such that $T_zX_\mu\subset T_zB_\mu(t;\epsilon)$. This is interpreted algebraically by the condition
\[
\mathrm{rank}
\left[
\begin{pmatrix}
t-z\\
\mathrm{Jac}_{X_\mu}(z)
\end{pmatrix}
\right]\le c_\mu\ .
\]
The {\em critical ideal} $\mathcal{I}_t\subset\C[\{z_{\omega_1\cdots\omega_s}\}\cup\{c_{\omega_1\cdots\omega_s}\}]$ is defined as (see \cite[equation (2.1)]{DHOST})
\begin{equation}\label{eq: critical ideal}
\mathcal{I}_t\coloneqq\left(\mathcal{I}_{X_\mu}+\left\langle (c+1)\times(c+1)-\mbox{minors of}
\begin{pmatrix}
t-z\\
\mathrm{Jac}_{X_\mu}(z)
\end{pmatrix}
\right\rangle\right)\colon\mathcal{I}_0^\infty,
\end{equation}
where $\mathcal{I}_0$ is the ideal generated by the coordinates $\{z_{\omega_1\cdots\omega_s}\}$, in particular $\mathcal{V}(\mathcal{I}_0)$ is the origin of $S^\mu V$. Note that $X_\mu\subset S^\mu V$ is smooth away from the origin. The following is a restatement of \cite[Lemma 2.1]{DHOST} in the context of $\mu$-symmetric binary tensors.

\begin{lemma}\label{lem:critical ideal}
For a general $t\in S^\mu V$, the variety of the critical ideal $\mathcal{I}_t$ is finite. It consists precisely of the critical binary tensors of rank one for $t$.
\end{lemma}

The following definition is an instance of \cite[Definition 2.2]{OS}.

\begin{definition}\label{def:psi}
The ideal $\mathcal{I}_t+(\widetilde{q}(t-z)-\epsilon^2)$ in the ring $\C[\{z_{\omega_1\cdots\omega_s}\}\cup\{t_{\omega_1\cdots\omega_s}\}\cup\{\epsilon\}]$ defines a variety of dimension $\dim(S^\mu V)=\prod_{k=1}^s(\mu_k+1)$ (see \cite[Theorem 4.1]{DHOST}). Since the polynomial ring is a UFD, its projection (eliminating $\{z_{\omega_1\cdots\omega_s}\}$) in $\C[\{t_{\omega_1\cdots\omega_s}\}\cup\{\epsilon\}]$ is generated by a single polynomial in $\epsilon^2$ . We denote this generator (defined up to a scalar factor) by $\mathrm{EDpoly}_{X_\mu,t}(\epsilon^2)$ and we call it the {\em Euclidean Distance polynomial (ED polynomial)} of $X_\mu$ at $t$.
\end{definition}

For any fixed $\epsilon\in \R$, the variety defined in $S^\mu V_\R$ by the vanishing of $\mathrm{EDpoly}_{X_\mu,t}(\epsilon^2)$ coincides with the clssical $\epsilon$-offset hypersurface of $X_{\mu,\R}$. In this paper, we consider $t\in S^\mu V$ as fixed and $\epsilon\in\C$ as a variable, hence we view $\mathrm{EDpoly}_{X_\mu,t}(\epsilon^2)$ as a univariate polynomial, in the notation of (\ref{eq: algebraic relation}). The first property of $\mathrm{EDpoly}_{X_\mu,t}(\epsilon^2)$ that we mention deals with its $\epsilon^2$-degree.

\begin{proposition}\cite[Theorem 2.7]{HW}\label{prop: HW} The $\epsilon^2$-degree of the ED polynomial of $X_\mu$ at $t\in S^\mu V$ is
\[
\deg \mathrm{EDpoly}_{X_\mu,t}(\epsilon^2) = \mathrm{EDdegree}(X_\mu) = s!\mu_1\cdots\mu_s,
\]
where the second equality follows by equation (\ref{eq: EDdegree Segre-Veronese}).
\end{proposition}

The following proposition is a particular case of \cite[Proposition 2.3]{OS} and describes which are the roots of the ED polynomial already defined.

\begin{proposition}\label{pro:roots charpoly}
For a general $t\in S^\mu V$, the roots of $\mathrm{EDpoly}_{X_\mu,t}(\epsilon^2)$ are precisely of the form $\epsilon^2=\widetilde{q}(t-z)$, where $z$ is a critical binary tensor of rank one for $t$ on $X_{\mu}$. In particular the distance $\epsilon$ from $X_{\mu,\R}$ to $t\in S^\mu V_\R$ is a root of $\mathrm{EDpoly}_{X_\mu,t}(\epsilon^2)$. Moreover $t\in S^\mu V_\R$ satisfies $\mathrm{EDpoly}_{X_\mu,t}(0)=0$ (namely it is a root of the lowest term of the ED polynomial) if and only if $t\in X_{\mu,\R}$.
\end{proposition}

Anyway, this is not the end of our construction. Indeed, as anticipated in the introduction, we consider the ED polynomial of the dual variety of $X_\mu$, rather than the ED polynomial of the variety $X_\mu$ itself. First we recall a definition.

\begin{definition}\label{def: dual variety} Let $Z\subset S^\mu V$ be an affine cone with the origin as vertex, meaning that if $z\in Z$ then $\lambda z\in Z$ for all $\lambda\in\C$. The {\em dual variety} of $Z$ is
\[
Z^\vee\coloneqq\overline{\bigcup_{z\in Z_{\mathrm{sm}}}N_zZ}\subset S^\mu V,
\]
where $N_zZ\coloneqq (T_zZ)^\perp=\{h\mid \widetilde{q}(h-z,y)=0\ \forall y\in T_zZ\}$ is the normal space at the smooth point $z\in Z$. We stress that we are identifying the vector space $S^\mu V$ and its dual with respect to the quadratic form $\widetilde{q}$. Moreover, note that $Z^\vee$ is the empty set if $Z$ is a non-reduced variety. 
\end{definition}

The Pythagorean Theorem and \cite[Theorem 5.2]{DHOST} tell us that the passage between the ED polynomial of $X_\mu$ and the ED polynomial of $X_\mu^\vee$ is just a variable reflection (see \cite[Theorem 3.2]{OS}). 
\begin{proposition}\cite[Theorem 3.2]{OS}\label{pro: pythagorean duality}
For any $t\in S^\mu V$,
\begin{equation}
\mathrm{EDpoly}_{X_\mu,t}(\epsilon^2)=\mathrm{EDpoly}_{X_\mu^\vee,t}(\widetilde{q}(t)-\epsilon^2).
\end{equation}
\end{proposition}

The next result clarifies the reason why we concentrate on ED polynomials of dual varieties of Segre-Veronese varieties.

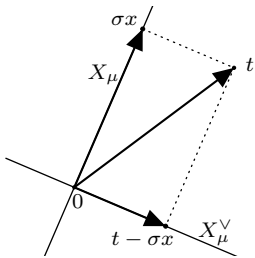
\begin{figure}[ht]
\begin{tikzpicture}[line cap=round,line join=round,>=triangle 45,x=1.0cm,y=1.0cm,scale=0.3]
\clip(-3.,-3.) rectangle (8.,8.);
\draw [line width=0.5pt,domain=-3.:8.] plot(\x,{(-0.--7.*\x)/3.});
\draw [line width=0.5pt,domain=-3.:8.] plot(\x,{(-0.-3.*\x)/7.});
\draw [->,line width=0.8pt,color=black] (0.,0.) -- (7.004482758620688,5.283793103448277);
\draw [->,line width=0.8pt,color=black] (0.,0.) -- (3.,7.);
\draw [->,line width=0.8pt,color=black] (0.,0.) -- (4.004482758620689,-1.7162068965517239);
\draw [line width=0.5pt,dotted] (3.,7.)-- (7.004482758620688,5.283793103448277);
\draw [line width=0.5pt,dotted] (7.004482758620688,5.283793103448277)-- (4.004482758620689,-1.7162068965517239);
\begin{scriptsize}
\draw [fill=black] (0.,0.) circle (2.0pt);
\draw[color=black] (0.15,-0.6) node {$0$};
\draw [fill=black] (3.,7.) circle (2.5pt);
\draw[color=black] (2.2,7.36) node {$\sigma x$};
\draw [fill=black] (4.004482758620689,-1.7162068965517239) circle (2.5pt);
\draw[color=black] (3,-2.3) node {$t-\sigma x$};
\draw [fill=black] (7.004482758620688,5.283793103448277) circle (2.0pt);
\draw[color=black] (7.7,5.44) node {$t$};
\draw[color=black] (6.2,-2) node {$X_{\mu}^\vee$};
\draw[color=black] (1.3,5) node {$X_{\mu}$};
\end{scriptsize}
\end{tikzpicture}
\caption{Singular tensors $\sigma x\in X_\mu$ and critical points $t-\sigma x\in X_\mu^\vee$ for the distance function $d_t$ on $X_\mu^\vee$ are in correspondence via the Pythagorean Theorem.}\label{fig: Pythagorean Theorem}
\end{figure}

\begin{proposition}\label{pro: roots and singular values}
For any $\mu$-symmetric binary tensor $t\in S^\mu V$ and any singular tensor $\sigma x\in S^\mu V$ for $t$, we have $\mathrm{EDpoly}_{X_\mu^\vee,t}(\sigma^2)=0$.
\end{proposition}
\proof
By Proposition \ref{pro:roots charpoly}, the roots of $\mathrm{EDpoly}_{X_\mu,t}(\epsilon^2)$ are of the form $\epsilon^2=\widetilde{q}(t-z)$, where $z$ is a critical binary tensor of rank one for $t$ on $X_{\mu}$. Moreover, by Theorem \ref{thm: Lim Qi} the non-isotropic critical tensors of rank one for $t$ correspond to the singular tensors for $t$. Then, consider a singular tensor $\sigma x$ for $t$. The root $\widetilde{q}(t-\sigma x)$ of $\mathrm{EDpoly}_{X_\mu,t}(\epsilon^2)$ corresponds, via Proposition \ref{pro: pythagorean duality}, with the root
\[
\widetilde{q}(t)-\widetilde{q}(t-\sigma x)=2\widetilde{q}(t,\sigma x)-\widetilde{q}(\sigma x)=2\sigma\widetilde{q}(t,x)-\sigma^2=2\sigma^2-\sigma^2=\sigma^2,
\]
of $\mathrm{EDpoly}_{X_\mu^\vee,t}(\epsilon^2)$ (see Figure \ref{fig: Pythagorean Theorem}), where we used the fact that $\sigma=\widetilde{q}(t,x)$ for any singular tensor $\sigma x$ for $t$, which is a direct consequence of equation (\ref{eq: singular vector tuple}).\qedhere
\endproof

\begin{remark}
The converse of Proposition \ref{pro: roots and singular values} is true only for general tensors. Indeed, there exist tensors $t\in S^\mu V$ such that some of the roots of $\mathrm{EDpoly}_{X_\mu^\vee,t}(\epsilon^2)$ do not correspond to singular values of $t$. In the symmetric case $\mu=(d)$, this phenomenon is studied in detail for example in \cite[Theorem 4]{Q2} and in \cite{LQZ}.
\end{remark}

\begin{remark}
On one hand, one may verify from Proposition \ref{pro: pythagorean duality} that
\begin{equation}\label{eq: identity leadcoefs}
\mathrm{EDpoly}_{X_{\mu},t}(\epsilon^2)=\mathrm{EDpoly}_{X_{\mu}^\vee,t}(\widetilde{q}(t)-\epsilon^2)=\sum_{k=0}^{N}(-1)^k\left[\sum_{j=k}^N\binom{j}{k}\widetilde{q}(t)^{j-k}a_j(t)\right]\epsilon^{2k}.
\end{equation}
Hence the highest terms of $\mathrm{EDpoly}_{X_{\mu},t}(\epsilon^2)$ and $\mathrm{EDpoly}_{X_{\mu}^\vee,t}(\widetilde{q}(t)-\epsilon^2)$ are equal to $a_N(t)$ up to sign. On the other hand, the corresponding lowest terms are not proportional.
\end{remark}

Summing up, a first way to compute the ED polynomial of $X_\mu^\vee$ at $t\in S^\mu V$ is by applying the original definition of ED polynomial of an algebraic variety together with Proposition \ref{pro: pythagorean duality}. The following is a Macaulay2 code \cite{GS} for computing the ED polynomial of $X_\mu^\vee$ in the symmetric case $\mu=(d)$, namely when $X_\mu$ is the rational normal curve of degree $d\ge 2$.

\begin{verbatim}
R = QQ[z_0..z_d, c_0..c_d, e];
RationalNormalCurve = minors(2, matrix{toList(z_0..z_(d-1)),toList(z_1..z_d)});
Jac = compress transpose jacobian RationalNormalCurve;
M = matrix{apply(d+1, j-> binomial(d, j)*(z_j-c_j))};
It = saturate( RationalNormalCurve + minors(d, M||Jac), ideal(toList(z_0..z_d)));
Hyperball = ideal(sum(d+1, j-> binomial(d,j)*z_j^2)-e^2);
EDpoly = (eliminate(toList(z_0..z_d), It+Hyperball))_0;
\end{verbatim}

The output \verb+EDpoly+ is the ED polynomial of the dual variety of \verb+RationalNormalCurve+ because we are intrinsically applying Proposition \ref{pro: pythagorean duality} in the definition of \verb+Hyperball+. Indeed, the usual relation $\widetilde{q}(t-z)-\epsilon^2$ is replaced by $\widetilde{q}(z)-\epsilon^2$ (see the definition in (\ref{eq: def hyperball})). Moreover, we stress that the metric $\widetilde{q}$ used in $S^d V$ is the one defined in equation (\ref{eq: metric symmetric tensors}). With this choice, the ED polynomial \verb+EDpoly+ of $X_{(d)}^\vee$ has degree $d=\mathrm{EDdegree}(X_{(d)})$ in \verb+e^2+.

Unfortunately, with this approach the symbolic computation of the ED polynomial is very hard even in the symmetric case $\mu=(d)$ for small values of $d$. The main reason lies in the computation of the critical ideal $\mathcal{I}_t$. Actually, Theorem \ref{thm: Lim Qi} and Proposition \ref{pro: roots and singular values} provide a more effective way for computing $\mathrm{EDpoly}_{X_{\mu}^\vee,t}(\epsilon^2)$, described in the following corollary.

\begin{corollary}\label{corol: second approach for computing EDpoly}
Let $\mu=(\mu_1,\ldots,\mu_s)\vdash d$. Fix coordinates $x_j=(x_{j,0},x_{j,1})$ for the $j$-th factor $V$ appearing in $\mathrm{Seg}_\mu(\PP(V)^{\times s})=X_\mu\subset\PP(S^\mu V$), for all $1\le j\le s$. Given a general $\mu$-symmetric binary tensor $t=(c_{\omega_1\cdots\omega_s})\in S^\mu V$, define $\mathcal{J}_t\subset\C[\{x_{j,0},x_{j,1}\},\{c_{\omega_1\cdots\omega_s}\},\epsilon]$ to be the ideal generated by all the relations in equation (\ref{eq: singular vector tuple}), when restricted to $S^\mu V\subset V^{\otimes d}$. Then
\[
\left(\mathcal{J}_t + \left\langle q(x_j)-1\mid 1\le j\le s\right\rangle\right)\cap \C[\{c_{\omega_1\cdots\omega_s}\},\epsilon]=\left\langle\mathrm{EDpoly}_{X_{\mu}^\vee,t}(\epsilon^2)\right\rangle.
\]
\end{corollary}

Below we give a second and more efficient Macaulay2 code for computing the ED polynomial of $X_\mu^\vee$ in the symmetric case $\mu=(d)$, with the alternative approach stated in Corollary \ref{corol: second approach for computing EDpoly}: 
\begin{verbatim}
R = QQ[x_0, x_1, u_0..u_d, e];
t = sum(d+1, j-> binomial(d,j)*u_j*x_0^(d-j)*x_1^j);
I = ideal(first entries((1/d)*diff(matrix{{x_0, x_1}}, t)-e*matrix{{x_0, x_1}}));
EDpoly = (eliminate({x_0, x_1},I + ideal(x_0^2+x_1^2-1)))_0;
\end{verbatim}

Note that in this case equations (\ref{eq: singular vector tuple}) simplify as in (\ref{eq: eigenvector eigenvalue}), which in turn correspond to the following system, where we interpret $t$ as a binary form of degree $d$:
\[
\frac{1}{d}\nabla t(x_0,x_1)=\lambda (x_0,x_1).
\]
For more details about the output of the above code we refer to Section \ref{sec: 222}. Moreover, the case of the ED polynomial of the dual of the rational normal curve of degree $d$ is studied in \cite{Sod} in the more general context of Veronese varieties.

\section{On the vanishing loci of the extreme coefficients}\label{sec: extreme coefficients}

In this and in the following sections, we focus on the lowest and highest coefficients of the ED polynomial of $X_\mu^\vee$ at a given $\mu$-symmetric tensor $t\in S^\mu V$. We consider $\mathrm{EDpoly}_{X_\mu^\vee,t}(\epsilon^2)$ written as in (\ref{eq: algebraic relation}). First of all, we investigate the vanishing loci of $a_0(t)$ and $a_N(t)$. The computation of the exponents of the factors in $a_0(t)$ and $a_N(t)$ is postponed to Section \ref{sec: main}.

A crucial role in the proof of the Main Theorem is played by the family of varieties $X_{\mu,J}$ with $J\subset[s]$ introduced in (\ref{eq: X mu j}). First we determine for which $J\subset[s]$ the variety $X_{\mu,J}^\vee$ is a hypersurface. 

\begin{lemma}\label{lem: when are hypersurfaces}
Consider the integers $0\le l\le s$ and a set of indices $J=\{k_1,\ldots,k_l\}\subset[s]$. The variety $X_{\mu,J}^\vee$ is a hypersurface unless $l=s-1$ and $\mu_j=1$ for $j\notin J$.
\end{lemma}
\proof
Without loss of generality, we may assume that $J=[l]\subset[s]$. The proof is a slight modification of \cite[V, Corollary 5.11]{GKZ}.
We need to verify that
\begin{align*}
\dim(Q)+\mathrm{codim}(Q^\vee)-1 &\le l\dim(Q)+(s-l)\dim(\PP(V))\quad\forall\ 1\le j\le l\\
\dim(\PP(V))+\mathrm{codim}(\PP(V)^\vee)-1 &\le l\dim(Q)+(s-l)\dim(\PP(V))\quad\forall\ l+1\le j\le s,
\end{align*}
where the first $l$ inequalities are related to the isotropic quadric $Q=\{[1,\sqrt{-1}],[1,-\sqrt{-1}]\}\subset\PP(V)$ in the Veronese embedding into $\PP(S^{\mu_j}V)$, whereas the remaining $s-l$ inequalities correspond to the projective space $\PP(V)$ in the Veronese embedding into $\PP(S^{\mu_j}V)$.

Pick any $1\le j\le l$. For any $\mu_j\ge 1$, $[\mathrm{Seg}_{\mu_j}(Q)]^\vee\subset\PP(S^{\mu_j}V)$ is the union of two conjugate hyperplanes. Hence the first $l$ inequalities are all equal to the condition $s-l\ge 0$, which is trivially satisfied. Now suppose that $l<s$ and consider the remaining $s-l$ inequalities. Pick any $l+1\le j\le s$. On one hand, if $\mu_j>1$, then $[\mathrm{Seg}_{\mu_j}(\PP(V))]^\vee\subset\PP(S^{\mu_j}V)$ is a hypersurface. Hence the corresponding inequality simplifies to $s-l\ge 1$, which is trivially satisfied. On the other hand, if $\mu_j=1$, then $\PP(V)^\vee=\emptyset$ and therefore $\mathrm{codim}(\PP(V)^\perp)-1=1$. Thus the corresponding inequality simplifies to $s-l\ge 2$. Therefore, in this case $X_{\mu,J}^\vee$ is not a hypersurface if and only if $l=s-1$.
\endproof

The vanishing locus of the lowest coefficient of the ED polynomial of an algebraic variety is completely described in the following result.

\begin{proposition}\label{lem: constterm}
The set of tensors $t\in S^\mu V$ which admit a $\mu$-symmetric critical tensor of rank one $z$ such that $\widetilde{q}(t-z)=0$ is
\begin{equation}\label{eq: components constterm}
\mathcal{V}(a_0)=X_{\mu}^\vee\cup(X_{\mu}\cap \widetilde{Q})^\vee.
\end{equation}
Moreover, if $X_{\mu}^\vee$ is a hypersurface, its equation $f_\mu$ appears with multiplicity two in $a_0$.
\end{proposition}

\proof
The identity (\ref{eq: components constterm}) follows immediately by \cite[Corollary 5.5]{OS}, hence we only need to prove that $f_\mu$ appears with multiplicity two in $a_0$. Note that we cannot apply directly \cite[Corollary 7.1]{OS}, since $X_\mu$ is not transversal to the isotropic quadric $\widetilde{Q}\subset\PP(S^\mu V)$. Nevertheless,
quadric hypersurfaces of $\PP(S^\mu V)$ that are transversal to $X_\mu$ and $X_\mu^\vee$ form a dense open subset $U\subset\PP(S^2S^\mu V)$. In particular, $\widetilde{Q}$ is the limit of a sequence $\{\widetilde{Q}_j\}\subset U$. Let $\mathrm{EDpoly}^{(j)}_{X_\mu^\vee,t}(\epsilon^2)$ be the ED polynomial of $X_\mu^\vee$ at $t\in S^\mu V$ with respect to te quadric $\widetilde{Q}_j$, for all $j$. By \cite[Corollary 7.1]{OS}, for all $j$ we have
\[
\mathrm{EDpoly}^{(j)}_{Y,y}(0)=f_\mu^2\cdot g_j,
\]
where $g_j$ is the equation of $(X_\mu\cap\widetilde{Q}_j)^\vee$. Moreover, by (\ref{eq: components constterm}) we know that $\mathrm{EDpoly}_{X_\mu^\vee,t}(0)=f_\mu^\alpha\cdot g^\beta$ for some non-negative integers $\alpha$ and $\beta$, where $g$ is the equation of $(X_\mu\cap\widetilde{Q})^\vee$. In particular,
\[
f_\mu^\alpha\cdot g^\beta\cdot h=\mathrm{EDpoly}_{X_\mu^\vee,t}(0)\cdot h=\lim_{j\to\infty}\mathrm{EDpoly}^{(j)}_{X_\mu^\vee,t}(0)=\lim_{j\to\infty}f_\mu^2\cdot g_j=f_\mu^2\cdot\lim_{j\to\infty}g_j,
\]
for some homogeneous polynomial $h$, possibly a scalar. In particular, $\alpha\ge 2$.

We show that actually $\alpha=2$. If $\alpha\ge 3$, then $f_\mu$ divides $\lim_{j\to\infty}g_j$, that is, $f_\mu$ divides $g$ or $f_\mu$ divides $h$. It remains to show that $f_\mu$ cannot divide $g$. In particular, our claim is that $\mathrm{codim}_\R[(X_\mu\cap\widetilde{Q})^\vee]\ge 2$. Consider a smooth point $z\in X_\mu\cap\widetilde{Q}$ and the normal space $S_{\mu,z}\coloneqq N_z(X_\mu\cap Q)$. Assume that $l_1,\ldots,l_r$ are the linear polynomials defining $S_{\mu,z}$. We denote by $\overline{S}_{\mu,z}$ the variety defined by $\bar{l}_1,\ldots,\bar{l}_n$, where the bar means complex conjugation. If $z\in\overline{S}_{\mu,z}$, then $q(z-\bar{z},y)=0$ for all $y\in T_z(X_\mu\cap\widetilde{Q})$. In particular, $q(\bar{z},z)=q(\bar{z},z)-q(z,z)=q(\bar{z}-z,z)=0$, contradiction. This implies that $S_{\mu,z}\neq\overline{S}_{\mu,z}$ and, in turn, that $\mathrm{codim}_\R(S_{\mu,z})\ge 2$. The claim follows by Definition \ref{def: dual variety}.
\endproof

Observe that an immediate consequence of Lemma \ref{lem: when are hypersurfaces} is that $X_\mu^\vee=X_{\mu,\emptyset}^\vee$ is always a hypersurface except for the trivial case $s=1$, $d=1$. Now we take a closer look at the other component $(X_\mu\cap\widetilde{Q})^\vee$ in (\ref{eq: components constterm}). By definition, the cone over the variety $X_\mu\cap\widetilde{Q}$, which we keep calling $X_\mu\cap\widetilde{Q}$, is isomorphic to $\{(x_1,\ldots,x_s)\in V^{\times s}\mid \widetilde{q}(x_1^{\mu_1}\otimes\cdots\otimes x_s^{\mu_s})=0\}$.
Define
\[
Y_{\mu,j}\coloneqq\{(x_1,\ldots,x_s)\in V^{\times s}\mid q(x_j)=0\},\quad 1\le j\le s.
\]
Then clearly $X_\mu\cap\widetilde{Q}\cong Y_{\mu,1}\cup\cdots\cup Y_{\mu,s}$.
\begin{lemma}\label{lem: Y_j}
Let $\mu=(\mu_1,\ldots,\mu_s)\vdash d$. For all $1\le j\le s$, $(Y_{\mu,j})_{red}\cong X_{\mu,\{j\}}$, where $(Y_{\mu,j})_{red}$ denotes the reduced locus of $Y_{\mu,j}$. Moreover if $\mu_j>1$, then $Y_{\mu,j}^\vee=\emptyset$.
\end{lemma}
\proof
It follows immediately by the definition of $Y_{\mu,j}$ that its reduced locus is isomorphic to $X_{\mu,\{j\}}$. Consider any $1\le j\le s$ and $x=x_1^{\mu_1}\otimes\cdots\otimes x_s^{\mu_s}\in Y_{\mu,j}$. On one hand, $x\in\widetilde{Q}$ and the tangent space $T_t\widetilde{Q}$ is the hyperplane filled by all tensors $u$ such that $\widetilde{q}(x,u)=0$. On the other hand, $x\in X_\mu$ and the tangent space of $X_\mu$ at $x$ is
\[
T_x X_\mu =\left\langle x, v_1x_1^{\mu_1-1}\otimes\cdots\otimes x_s^{\mu_s}, \ldots, x_1^{\mu_1}\otimes\cdots\otimes v_sx_s^{\mu_s-1}\mid v_k\in V\ \mbox{for all}\ 1\le k\le s\right\rangle.
\]
For any $1\le k\le s$ pick a non-zero vector $v_k\in V$ and consider $x_1^{\mu_1}\otimes\cdots\otimes v_kx_k^{\mu_k-1}\otimes\cdots\otimes x_s^{\mu_s}\in T_x X_\mu$. Then we get
\[
\widetilde{q}(x,x_1^{\mu_1}\otimes\cdots\otimes v_kx_k^{\mu_k-1}\otimes\cdots\otimes x_s^{\mu_s})=q(x_1)^{\mu_1}\cdots q(v_k,x_k)\cdot q(x_k)^{\mu_k-1}\cdots q(x_s)^{\mu_s}\ \mbox{for all}\ 1\le k\le s.
\]
In particular, $\widetilde{q}(x,x_1^{\mu_1}\otimes\cdots\otimes v_kx_k^{\mu_k-1}\otimes\cdots\otimes x_s^{\mu_s})=0$ for all $k\neq j$. Now assume $k=j$. If $\mu_j=1$, then for a general $v_j$ we get $q(x_j,v_j)\neq 0$ and in turn $\widetilde{q}(x,x_1^{\mu_1}\otimes\cdots\otimes v_jx_j^{\mu_j-1}\otimes\cdots\otimes x_s^{\mu_s})\neq 0$. This implies that the general point $x\in Y_{\mu,j}$ is smooth if $\mu_j=1$. Otherwise if $\mu_1>1$, then $\widetilde{q}(x,x_1^{\mu_1}\otimes\cdots\otimes v_jx_j^{\mu_j-1}\otimes\cdots x_s^{\mu_s})=0$. This means that $T_x X_\mu\subset T_x \widetilde{Q}$, and every point $x\in Y_{\mu,j}$ is not smooth. Therefore, by the definition of dual variety we have that $Y_{\mu,j}^\vee=\emptyset$ if $\mu_j>1$.
\endproof
An immediate consequence of Lemma \ref{lem: Y_j} is the identity
\[
(X_\mu\cap\widetilde{Q})^\vee=\bigcup_{j\in [s]\colon \mu_j=1}X_{\mu,\{j\}}^\vee.
\]
In the second part of this section, we are interested in giving a complete description of the vanishing locus of the highest coefficient $a_N(t)$. In the following, $H_\infty\cong\PP(S^\mu V)$ denotes the hyperplane at infinity, while $Y_\infty\coloneqq Y\cap H_\infty$ for any variety $Y\subset S^\mu V$.

\begin{figure}[ht]
\begin{tikzpicture}[line cap=round,line join=round,>=triangle 45,x=1.0cm,y=1.0cm,scale=0.3]
\clip(-3.,-3.) rectangle (8.,8.);
\draw [line width=0.5pt,domain=-3.:8.] plot(\x,{(-0.--7.*\x)/3.});
\draw [line width=0.5pt,domain=-3.:8.] plot(\x,{(-0.-3.*\x)/7.});
\draw [->,line width=0.8pt,color=black] (0.,0.) -- (7.004482758620688,5.283793103448277);
\draw [->,line width=0.8pt,color=black] (0.,0.) -- (3.,7.);
\draw [->,line width=0.8pt,color=black] (0.,0.) -- (4.004482758620689,-1.7162068965517239);
\draw [line width=0.5pt,dotted] (3.,7.)-- (7.004482758620688,5.283793103448277);
\draw [line width=0.5pt,dotted] (7.004482758620688,5.283793103448277)-- (4.004482758620689,-1.7162068965517239);
\begin{scriptsize}
\draw [fill=black] (0.,0.) circle (2.0pt);
\draw[color=black] (0.2,-0.7) node {$0$};
\draw [fill=black] (3.,7.) circle (2.5pt);
\draw[color=black] (1.2,7.36) node {$t_k-f_k$};
\draw[color=black] (-7.64,3.05) node {$X'$};
\draw [fill=black] (4.004482758620689,-1.7162068965517239) circle (2.5pt);
\draw[color=black] (4,-2.5) node {$f_k$};
\draw [fill=black] (7.004482758620688,5.283793103448277) circle (2.0pt);
\draw[color=black] (7.7,5.44) node {$t_k$};
\draw[color=black] (6.2,-2) node {$X_{\mu}^\vee$};
\draw[color=black] (1.2,5) node {$X_{\mu}$};
\end{scriptsize}
\end{tikzpicture}
\caption{The sequences $\{f_k\}\subset X_{\mu}^\vee$ and $\{t_k-f_k\}\subset X_{\mu}$.}\label{fig: sequences}
\end{figure}
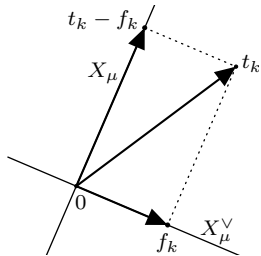

\begin{proposition}\label{prop: leadcoef of EDpoly}
The following inclusion holds true:
\[
\mathcal{V}(a_N)\subset\bigcup_{j\colon \mu_j>1}X_{\mu,\{j\}}^\vee\cup\bigcup_{|J|>1} X_{\mu, J}^\vee.
\]
\end{proposition} 
\proof
For any non-zero $\mu$-symmetric binary tensor $t\in S^\mu V$, we indicate by $[t]\in H_\infty$ the line spanned by $t$. Now assume that  $a_N(t)=0$. From this, from \cite[Theorem 3.2]{OS} and $(\ref{eq: identity leadcoefs})$, there exists a sequence $\{t_k\}\subset S^\mu V$ such that $t_k\to t$ and two corresponding sequences $\{f_k\}\subset X_{\mu}^\vee$ and $\{t_k-f_k\}\subset X_{\mu}$ of critical points for $d_{t_k}$ on $X_{\mu}^\vee$ and $X_\mu$, respectively, such that $\mathrm{EDpoly}_{X_\mu^\vee,t_k}(\epsilon_k^2)=0$ and $\mathrm{EDpoly}_{X_{\mu},t_k}(\eta_k^2)=0$ when $\epsilon_k^2=\widetilde{q}(f_k-t_k)$ and $\eta_k^2=\widetilde{q}(f_k)$ diverge simultaneously (see Figure \ref{fig: sequences}). In particular, we have that $[t_k-f_k]\in [(T_{f_k}X_{\mu}^\vee)_\infty]^\perp=(N_{f_k}X_{\mu}^\vee)_\infty$ and $[f_k]\in[(T_{t_k-f_k}X_{\mu})_\infty]^\vee=(N_{t_k-f_k}X_{\mu})_\infty$ for all $k$, where the external duals are taken in the projective subspace $H_{\infty}$.

Up to subsequences, we may assume that
\begin{equation}\label{eq:proposition integrality 1}
\lim_{k\to\infty}[f_k]\eqqcolon[f]\in(X_{\mu})_{\infty}^\vee, \textrm{\ for some\ }f\in S^\mu V.
\end{equation}
In the topology of the compact space $\overline{X_{\mu}}=X_{\mu}\cup (X_{\mu})_\infty$ we still have $f_k\to[f]\in (X_{\mu})_{\infty}^\vee$, more precisely in $\PP(\C\oplus V)$ we have $\left[(1,f_k)\right]\to\left[(0, f)\right]$. Consequently, we have that
\[
\lim_{k\to\infty}[t_k-f_k]=[t-f]\in (X_{\mu})_{\infty}.
\]
Repeating the argument of \cite[Proposition 4.4]{OS}, one verifies that
\begin{enumerate}
	\item[$i)$] $\widetilde{q}([t-f])=0$, namely $[t-f]\in\widetilde{Q}_{\infty}$, where in this case $\widetilde{q}$ is the quadratic form defined in $H_\infty$.
	\item[$ii)$] $T_{[t-f]}(X_{\mu})_\infty\subset T_{[t-f]}\widetilde{Q}_\infty$.
\end{enumerate}
Remembering that $t-f$ is a decomposable tensor, then $[t-f]=[x_1^{\mu_1}\otimes\cdots\otimes x_t^{\mu_t}]$ for some vectors $x_i\in V$. By $i)$, necessarily $q(x_j)=0$ for some $j$, say $j=1$. Now there are two possible cases to study.

If $\mu_1=1$, then we may suppose that $q(x_1)=0$. The inclusion $ii)$ implies that for any $[v]\in T_{[t-f]}(X_{\mu})_\infty$, $\widetilde{q}([v],[t-f])=0$, where in this case $\widetilde{q}(\cdot,\cdot)$ is the scalar product on $V$ restricted to $H_\infty$.
More explicitly, we may write
\[
v=\sum_{i=1}^t x_1\otimes\cdots\otimes x_{i-1}^{\mu_{i-1}}\otimes\xi_i\cdot x_i^{\mu_i-1}\otimes x_{i+1}^{\mu_{i+1}}\otimes\cdots\otimes x_t^{\mu_t}
\]
for some $\xi_1,\ldots,\xi_t\in V$. Then
\begin{align*}
0 &=\widetilde{q}([v],[t-f])\\
&=\sum_{i=1}^t \widetilde{q}([x_1\otimes\cdots\otimes x_{i-1}^{\mu_{i-1}}\otimes\xi_i\cdot x_i^{\mu_i-1}\otimes x_{i+1}^{\mu_{i+1}}\otimes\cdots\otimes x_t^{\mu_t}],[x_1^{\mu_1}\otimes\cdots\otimes x_t^{\mu_t}])\\
&=\sum_{i=1}^t q(x_1)\cdots q(x_{i-1})^{\mu_{i-1}}\cdot q(\xi_i,x_i)\cdot q(x_i)^{\mu_i-1}\cdot q(x_{i+1})^{\mu_{i+1}}\cdots q(x_t)^{\mu_t}.
\end{align*}
By our assumption $q(x_1)=0$. Then necessarily we get the identity $q(\xi_1,x_1)\cdot q(x_2)^{\mu_2}\cdots q(x_t)^{\mu_t}=0$.
Taking $v$ sufficiently general, we may suppose that $q(\xi_1,x_1)\neq 0$. Therefore there exists at least one more index $i\neq 1$ such that $q(x_i)=0$. In particular $t-f\in X_{\mu,J}$, where $J=\{j\in[s]\mid q(x_j)=0\}$.

If $\mu_1>1$, then $t-f\in X_{\mu,\{1\}}$, otherwise $t-f\in X_{\mu,J}$ for some $J\subset [s]$ such that $|J|>1$.

Now assume that $\mu_1=1$ (the proof in the case $\mu_1>1$ is the same). We show that necessarily $t\in X_{\mu,J}^\vee$. By definition, $X_{\mu,J}^\vee=\overline{S_{\mu,J}}$, where
\[
S_{\mu,J} = \bigcup_{x\in (X_{\mu,J})_{\mathrm{sm}}}N_xX_{\mu,J}.
\]
Now consider $y\in T_{t-f}X_{\mu,J}$. By the previous claims, $\widetilde{q}(y,t-f)=0$. On the other hand, we have $\widetilde{q}(y,f)=0$, since $X_{\mu}$ and $X_{\mu}^\vee$ are cones. Then $\widetilde{q}(y,t)=\widetilde{q}(y,t-f)+\widetilde{q}(y,f)=0+0=0$. This means that $t\in S_{\mu,J}$, hence $t\in X_{\mu,J}^\vee$.
\endproof

Proposition \ref{prop: leadcoef of EDpoly} leaves this information: if a $\mu$-symmetric binary tensor $t\in S^\mu V$ does not admit the expected number of critical points, then it must have a critical point with a precise isotropic structure.

\section{Proof of the Main Theorem}\label{sec: main}

For all $J\subset[s]$, we recall that $f_{\mu,J}$ denotes the equation of $X_{\mu,J}^\vee$, when it is a hypersurface. Otherwise, we set $f_{\mu,J}\coloneqq 1$. Moreover, we use the notation $f_{\mu}\coloneqq f_{\mu,\emptyset}$ for the $\mu$-discriminant. Proposition \ref{prop: leadcoef of EDpoly} carries the following fact: if $a_N=g_1^{\beta_1}\cdots g_r^{\beta_r}$ is the irreducible factorization of the highest coefficient $a_N$ of $\mathrm{EDpoly}_{X_\mu^\vee,t}(\epsilon^2)$, then the $g_k$'s are either proportional to $f_{\mu,\{j\}}$ for some $1\le j\le s$ such that $\mu_j>1$ or to $f_{\mu,J}$ for some $J\subset[s]$ such that $|J|>1$. In particular, we may write, up to scalars,
\begin{equation}\label{eq: writing the leadcoef}
a_N(t)=\prod_{j\colon \mu_j>1}f_{\mu,\{j\}}(t)^{\alpha_{\mu,\{j\}}}\cdot\prod_{|J|>1}f_{\mu,J}(t)^{\alpha_{\mu,J}},\quad \alpha_{\mu,J}\ge 0.
\end{equation}
Moreover, from Lemma \ref{lem: constterm} we have that, up to scalars,
\begin{equation}\label{eq: writing the constterm}
a_0(t)=f_{\mu}(t)^2\cdot\prod_{j\colon \mu_j=1}f_{\mu,\{j\}}(t)^{-\alpha_{\mu,\{j\}}},\quad \alpha_{\mu,\{j\}}\le 0\quad\mbox{when}\quad\mu_j=1.
\end{equation}
When $\mu=1^d$, we observe that $\mathrm{EDpoly}_{X_d^\vee,t}(\epsilon^2)$ is invariant under the action of the symmetric group $\Sigma_d$ on the entries of $t=(t_{i_1\cdots i_d})$. More precisely, we use the notation $\alpha_{d,j}\coloneqq\alpha_{1^d,J}$ for all subsets $J\subset[d]$. The reason for the negative sign in the notation for $\alpha_{\mu,J}$ in (\ref{eq: writing the constterm}) is explained in Proposition \ref{pro: linear system}.

A non trivial task is showing what are the exponents $\alpha_{\mu,J}$ appearing in the expressions of $a_N$ and $a_0$. On one hand, Corollary \ref{corol: preservation of exponents} simplifies a lot our problem, stating that $\alpha_{\mu,J}\in\{\alpha_{d,1},\ldots,\alpha_{d,d}\}$ for all $\mu\vdash d$ and all $J\subset[s]$. On the other hand, Corollary \ref{corol: case a_d} and Lemmas \ref{lem: case s=2} and \ref{lem: case s=3} produce linearly independent equations on the remaining unknowns $\alpha_{d,1},\ldots,\alpha_{d,d}$. Finally, the last step discussed in Proposition \ref{pro: linear system} concludes the proof of the Main Theorem.

In Lemma \ref{lem: when are hypersurfaces} we showed that the $\mu$-discriminant $f_\mu$ is always non trivial except for the case $d=1$. Its degree is recalled in the following result.

\begin{theorem}{\cite[XIII, Theorem 2.4]{GKZ}}\label{thm: degree mu discriminant}
Let $d\ge 1$ and consider $\mu=(\mu_1,\ldots,\mu_s)\vdash d$. Define $e_j(\mu_1,\ldots,\mu_s)=\sum_{1\le k_1<\cdots< k_j\le s}\mu_{k_1}\cdots \mu_{k_j}$ to be the elementary symmetric polynomial of degree $j$. Let $\delta_\mu\coloneqq\deg(f_\mu)$ (it is usually denoted by $N(1^d;\mu)$). Then
\begin{equation}\label{eq: degree mu discriminant}
\delta_\mu=\sum_{j=0}^s(-2)^{s-j}(j+1)!e_j(\mu).
\end{equation}
\end{theorem}

Note that in the non-symmetric case $\mu=1^d$, we have $e_j(1^d)=\binom{d}{j}$ for all $0\le j\le d$. Hence we recover the degree $\delta_d\coloneqq \delta_{1^d}$ of the hyperdeterminant of a $d$-dimensional binary tensor:
\begin{equation}\label{eq: degree hyperdeterminant}
\delta_d=\sum_{j=0}^d(-2)^{d-j}\binom{d}{j}(j+1)!.
\end{equation}
An almost immediate consequence of Theorem \ref{thm: degree mu discriminant} is a formula for the degree of $X_{\mu,J}^\vee$.

\begin{corollary}\label{corol: degree varieties f mu J}
For any subset $J\subset [s]$, let $\mu(J)\vdash d-\sum_{j\in J}\mu_j$ be the partition whose summands are all the $\mu_k$ such that $k\in[s]\setminus J$. In particular, $\mu(\emptyset)=\mu$. Then
\begin{equation}\label{eq: degree f mu j}
\deg(f_{\mu,J})=2^{|J|}\delta_{\mu(J)}=2^{|J|}\sum_{j=0}^{s-|J|}(-2)^{s-|J|-j}(j+1)!e_j(\mu(J)).
\end{equation}
\end{corollary}

The degrees $\deg(f_{\mu,J})$ appear in the following identities, which descend from relations (\ref{eq: writing the leadcoef}), (\ref{eq: writing the constterm}) and the identity $\deg(a_N)+2\mathrm{EDdegree}(X_\mu)=\deg(a_0)$, where we already know that $\alpha_{\mu,\emptyset}=\alpha_{\mu}=-2$:

\begin{equation}\label{eq: identity between degrees}
\sum_{J\subset[s]}\alpha_{\mu,J}\deg(f_{\mu,J})+2\mathrm{EDdegree}(X_\mu)=0\quad\mbox{for all}\quad\mu\vdash d.
\end{equation}

The main idea of the proof of the Main Theorem is related to partial symmetrizations of the ED polynomial of $X_d^\vee$. To this aim, we recall some definitions and preliminary results.

\begin{definition}\label{def: refinement}
Let $\mu=(\mu_1,\ldots,\mu_s)$ be a partition of $d$. A {\em symmetrization} of $\mu$ is any partition $\lambda=(\lambda_1,\ldots,\lambda_r)$ of $d$ with
\begin{equation}\label{eq: refinement}
\lambda_j=\mu_{i_{j,1}}+\cdots+\mu_{i_{j,l_j}}\quad\mbox{for all}\ 1\le j\le r,
\end{equation}
where $\mu=(\mu_{i_{1,1}},\ldots,\mu_{i_{1,l_1}},\ldots,\mu_{i_{r,1}},\ldots,\mu_{i_{r,l_r}})$ after a possible permutation. We write $\lambda\prec\mu$ to indicate that $\lambda$ is a symmetrization of $\mu$. We stress that different choices of $\mu_i$ appearing in different sums (\ref{eq: refinement}) yield different symmetrizations of $\mu$, even if some of the $\mu_i$ are equal.
\end{definition}

Given two partitions $\lambda=(\lambda_1,\ldots,\lambda_r)$ and $\mu=(\mu_1,\ldots,\mu_s)$ such that $\lambda\prec\mu$, we may consider the inclusion $S^{\lambda}V\subset S^{\mu}V$. Since the group $\mathrm{GL}(V)$ is reductive, there exists a unique $\mathrm{GL}(V)$-invariant complement to $S^{\lambda}V$ in $S^{\mu}V$, denoted by $W^{\lambda,\mu}$. We have a natural projection $\pi_{\lambda,\mu}:\PP(S^{\mu}V)\dasharrow\PP(S^{\lambda}V)$ from $W^{\lambda,\mu}$, whose definition on decomposable elements is
\[
\pi_{\lambda,\mu}\left([a_{i_{1,1}}^{\mu_{i_{1,1}}}\otimes\cdots\otimes a_{i_{1,t_1}}^{\mu_{i_{1,t_1}}}\otimes\cdots\otimes a_{i_{s,1}}^{\mu_{i_{s,1}}}\otimes\cdots\otimes a_{i_{t,t_s}}^{\mu_{i_{t,t_s}}}]\right)\coloneqq[a_{i_{1,1}}^{\mu_{i_{1,1}}}\cdots a_{i_{1,t_1}}^{\mu_{i_{1,t_1}}}\otimes\cdots\otimes a_{i_{s,1}}^{\mu_{i_{s,1}}}\cdots a_{i_{t,t_s}}^{\mu_{i_{t,t_s}}}].
\]
The projection $\pi_{\lambda,\mu}$ induces another projection $\PP(S^e(S^{\mu}V))\dasharrow\PP(S^e(S^{\lambda}V))$, which we keep calling $\pi_{\lambda,\mu}$. The {\em $\lambda$-symmetrization} of a degree $e$ homogeneous polynomial $f\in S^e(S^{\mu}V)$ is the image $\pi_{\lambda,\mu}(f)$ under the map already defined. The following result is an instance of \cite[Theorem 2.2]{HO}, and an almost immediate consequence of \cite[Theorem 1.2]{O}. It relates the $\lambda$-symmetrization of the $\mu$-discriminant (the equation of $X_\mu^\vee$), where $\lambda\prec\mu$ are two partitions of $d$.

\begin{proposition}\label{prop: symmetrization X mu X lambda}
Let $\lambda\prec\mu$ be two partitions of $d$. Then
\[
X_\lambda^\vee\subset X_\mu^\vee\cap\PP(S^\lambda V).
\]
Moreover, $f_\lambda$ is a factor of multiplicity one in $\pi_{W^{\lambda,\mu}}(f_\mu)$.
\end{proposition}
\proof
By \cite[Theorem 1.2]{O}, we have the inclusions $X_{(d)}^\vee\subset X_\mu^\vee\cap\PP(S^d V)$ and $X_{(d)}^\vee\subset X_\lambda^\vee\cap\PP(S^d V)$. We stress that, taking into account Definition \ref{def: dual variety}, we are identifying $\PP((S^\mu V/S^d V)^\perp)$ with $S^dV$ and by abuse of notation write $X_{(d)}^\vee\subset X_\mu^\vee\cap \PP(S^d V)$. Again by \cite[Theorem 1.2]{O}, the discriminant $f_{(d)}$ is a factor of multiplicity one in both the polynomials $\pi_{(d),\mu}(f_\mu)$ and $\pi_{(d),\lambda}(f_\lambda)$. 
\endproof

\begin{remark}\label{rmk: isomorphims between Seg Ver varieties}
Fix a partition $\mu\vdash d$ and a vector space $S^\mu V$. As a completion of the remark in Definition \ref{def: refinement}, we stress that distinct symmetrizations $\lambda_1\neq\lambda_2$ of $\mu$ yield isomorphic but distinct subspaces $S^{\lambda_1}V\neq S^{\lambda_2}V$ and, in turn, distinct (but isomorphic) varieties $X_{\lambda_1}$ and $X_{\lambda_2}$, even if $\lambda_1=\lambda_2$ as partitions of $d$. For example, when we write $X_{(2,1)}\subset X_{(1,1,1)}$ we do take into account which components of $\mu=(1,1,1)$ we are summing to get $\lambda=(2,1)$. Nevertheless, we chose to omit this assumption in our notation. In terms of projections $\pi_{\lambda,\mu}$, different symmetrizations $\lambda_1\neq\lambda_2$ of $\mu$ yield different maps $\pi_{\lambda_1,\mu}\neq\pi_{\lambda_2,\mu}$. 
\end{remark}

The following key fact shows how the previous result shifts from dual Segre-Veronese varieties to their respective ED polynomials.

\begin{proposition}\label{pro: ed poly divides ed poly}
Let $\lambda\prec\mu$ be two partitions of $d$ and let $t\in S^\lambda V$. Then the ED polynomial of $X_\lambda^\vee$ at $t$ divides with multiplicity one the $\lambda$-symmetrization of the ED polynomial of $X_\mu^\vee$ at $t$. 
\end{proposition}
\proof
Let $t\in S^\lambda V$ and let $x\in X_\lambda\subset X_\mu$ be a $\lambda$-symmetric singular tensor for $t$. We show that $x$ is also a $\mu$-symmetric singular tensor for $t$. According to the decomposition $S^\mu V=S^\lambda V\oplus W^{\lambda,\mu}$, the tangent space of $X_\mu$ at $x$ decomposes in a good way as $T_xX_\mu=T_xX_\lambda\oplus W$ for some subspace $W\subset W^{\lambda,\mu}$. In particular, any tangent vector of $X_\mu$ at $x$ may be written in a unique way as $y=y_\lambda+w$ for some $y_\lambda\in T_xX_\lambda$ and $w\in W$. Then we have
\[
\widetilde{q}(t-x,y)=\widetilde{q}(t-x,y_\lambda)+\widetilde{q}(t-x,w)=0+0=0
\]
for all $y\in T_xX_\mu$. Thanks to Proposition \ref{pro: roots and singular values}, this fact means, at the level of ED polynomials, that there exists an integer $\beta\ge 0$ such that $\mathrm{EDpoly}_{X_\mu^\vee,t}(\epsilon^2)=[\mathrm{EDpoly}_{X_\lambda^\vee,t}(\epsilon^2)]^\beta\cdot h$, and $\mathrm{EDpoly}_{X_\lambda^\vee,t}(\epsilon^2)$ is not a factor of $h$. By Lemma \ref{lem: constterm}, the equations of $X_\mu^\vee$ and $X_\lambda^\vee$, namely $f_\mu$ and $f_\lambda$, appear with multiplicity $2$ in the lowest terms of $\mathrm{EDpoly}_{X_\mu^\vee,t}(\epsilon^2)$ and $\mathrm{EDpoly}_{X_\lambda^\vee,t}(\epsilon^2)$. Moreover, by Proposition \ref{prop: symmetrization X mu X lambda}, $f_\lambda$ is a factor of multiplicity one in $\pi_{\lambda,\mu}(f_\mu)$. This implies that $\beta=1$.
\endproof

\begin{definition}\label{def: compatible}
Let $\mu=(\mu_1,\ldots,\mu_s)\vdash d$. Consider a subset of indices $J\subset[s]$. We say that a partition $\lambda=(\lambda_1,\ldots,\lambda_r)\prec \mu$ is {\em compatible} with $J$ if for all $1\le j\le r$ we may write
\[
\lambda_j=\mu_{i_{j,1}}+\cdots+\mu_{i_{j,l_j}}
\]
and the subset of indices $I_{\lambda,j}\coloneqq\{i_{j,1},\ldots,i_{j,l_j}\}$ is either contained in $J$ or in $[s]\setminus J$. Moreover, we define $J_\lambda\coloneqq\{j\in[r]\mid I_{\lambda,j}\subset J\}$.
\end{definition}

\begin{example}
Let $\mu=(\mu_1,\ldots,\mu_s)\vdash d$. Consider a non-empty subset $J=\{p_1,\ldots,p_n\}\subset[s]$ and its complement $[s]\setminus J=\{q_1,\ldots,q_{s-n}\}$. Then any of the following partitions $\lambda\prec \mu$, which we use in Corollary \ref{corol: preservation of exponents} and Lemmas \ref{lem: case s=2} and \ref{lem: case s=3}, is compatible with $J$:
\begin{enumerate}
\item $\lambda=(d)=(\mu_1+\cdots+\mu_s)$ if $J=[s]$,
\item $\lambda=(\lambda_1,\lambda_2)=(\mu_{p_1}+\cdots+\mu_{p_n},\mu_{q_1}+\cdots+\mu_{q_{s-n}})$,
\item  $\lambda=(\lambda_1,\lambda_2,\lambda_3)=(\mu_{p_1}+\cdots+\mu_{p_{n-1}},\mu_{p_n},\mu_{q_1}+\cdots+\mu_{q_{s-n}})$,
\item  $\lambda=(\lambda_1,\lambda_2,\lambda_3)=(\mu_{p_1}+\cdots+\mu_{p_n},\mu_{q_1}+\cdots+\mu_{q_{s-n-1}},\mu_{q_{s-n}})$.
\end{enumerate}
\end{example}

Definition \ref{def: compatible} is useful for introducing the following variation of Proposition \ref{prop: symmetrization X mu X lambda}.

\begin{proposition}\label{pro: symmetrization into two parts}
Let $\mu=(\mu_1,\ldots,\mu_s)$ be a partition of $d$. Consider a subset of indices $J\subset[s]$ and a partition $\lambda=(\lambda_1,\cdots,\lambda_r)\prec\mu$ compatible with $J$. If $X_{\mu,J}^\vee$ is a hypersurface, then $X_{\lambda,J_\lambda}^\vee$ is a hypersurface too, and
\begin{equation}\label{eq: inclusion symmetrization parts}
X_{\lambda,J_\lambda}^\vee\subset X_{\mu,J}^\vee\cap\PP(S^\lambda V).
\end{equation}
Moreover, $f_{\lambda,J_\lambda}$ is a factor of multiplicity one in $\pi_{\lambda,\mu}(f_{\mu,J})$.
\end{proposition}

\proof
Since the partition $\lambda$ is compatible with $J$, we conclude immediately from Definition \ref{def: compatible} that $X_{\lambda,J_\lambda}=X_{\mu,J}\cap\PP(S^\lambda V)$. Looking at its definition in (\ref{eq: X mu j}), we observe that $X_{\mu,J}$ is isomorphic to the union of $2^{|J|}$ copies of $X_{\mu(J)}$, where the partition $\mu(J)\vdash d-\sum_{j\in J}\mu_j$ was introduced in Corollary \ref{corol: degree varieties f mu J}. Analogously, $X_{\lambda,J_{\lambda}}$ is isomorphic to the union of $2^{|J_\lambda|}$ copies of $X_{\lambda(J_\lambda)}$.

Now let $C^{(\lambda)}\cong X_{\lambda(J_\lambda)}$ be a component of $X_{\lambda,J_{\lambda}}$. We may write $C^{(\lambda)}=\mathrm{Seg}_\lambda(W_1\times\cdots\times W_r)\subset\PP(S^\lambda V)$, where $W_j$ is one of the two points of $Q$ if $j\in J_\lambda$, otherwise $W_j=\PP(V)$. The variety $C^{(\lambda)}$ corresponds precisely to the $\lambda$-symmetrization of the component $C^{(\mu)}\cong X_{\mu(J)}$ of $X_{\mu,J}$, written as $C^{(\mu)}=\mathrm{Seg}_\mu(Z_1\times\cdots\times Z_s)\subset\PP(S^\mu V)$ with $Z_k=W_j$ for all $k\in I_{\lambda,j}$ and all $j\in J_\lambda$ (see Definition \ref{def: compatible}), otherwise $Z_k=\PP(V)$. The thesis follows by Proposition \ref{prop: symmetrization X mu X lambda}.
\endproof

Propositions \ref{pro: ed poly divides ed poly} and \ref{pro: symmetrization into two parts} yield the following useful corollary for the proof of the Main Theorem.

\begin{corollary}\label{corol: preservation of exponents}
Let $\mu=(\mu_1,\ldots,\mu_s)$ be a partition of $d$. Then $\alpha_{\mu,J}=\alpha_{d,\sum_{j\in J}\mu_j}$ for all $J\subset[s]$.
\end{corollary}
\proof
Let $K\subset[d]$ such that $|K|=\sum_{j\in J}\mu_j$. By Proposition \ref{pro: symmetrization into two parts}, $f_{\mu,J}$ is a factor of multiplicity one in $\pi_{\mu,1^d}(f_{d,K})$. Moreover, given any tensor $t\in S^\mu V$, by Proposition \ref{pro: ed poly divides ed poly} the ED polynomial of $X_\mu^\vee$ at $t$ divides the ED polynomial of $X_d^\vee$ at $t$ with multiplicity one. Therefore, the exponents of $f_{d,K}$ and $f_{\mu,J}$, which are respectively $\alpha_{d,|K|}=\alpha_{d,\sum_{j\in J}\mu_j}$ and $\alpha_{\mu,J}$, must coincide.
\endproof

Our method of partial symmetrizations starts from the easiest case, namely the symmetric case $\mu=(d)$, explained in the following result.

\begin{theorem}{\cite[Theorem 4.15]{Sod}}\label{thm: main theorem symmetric tensors}
Assume $s=1$, and let $d\ge 1$ be an integer. On one hand, the highest coefficient of $\mathrm{EDpoly}_{X_{(d)}^\vee,t}(\epsilon^2)$ is $a_N(t)=\Delta_Q(t)^{d-2}$, where $\Delta_Q(t)=f_{(d),\{1\}}$ is the equation of the dual of the Veronese embedding into $\PP(S^dV)$ of $Q\subset V$. On the other hand, the lowest coefficient is $c_0(t)=\Delta_d(t)^2$, where $\Delta_d(t)=f_{(d)}$ is the discriminant of the form $t$.
\end{theorem}

\begin{corollary}\label{corol: case a_d}
For any partition $\mu=(\mu_1,\ldots,\mu_s)\vdash d$, we have the relation
\begin{equation}\label{eq: equation partition 1 factor}
\alpha_{\mu,[s]}=\alpha_{(d),\{1\}}=\alpha_{d,d}=d-2.
\end{equation}
\end{corollary}

Last corollary solves the first problem of determining the exponent $\alpha_{d,d}$ of $f_{\mu,[s]}$, for any partition $\mu=(\mu_1,\ldots,\mu_s)\vdash d$. The next two technical lemmas furnish a bunch of linear conditions involving the remaining exponents $\alpha_{d,j}$, with $1\le j\le d-1$.

\begin{lemma}\label{lem: case s=2}
Let $d\ge 2$ and consider the partition $\mu=(k,d-k)\vdash d$ for all $1\le k\le \lfloor d/2 \rfloor$. The highest coefficient $a_N(t)$ and the lowest coefficient $a_0(t)$ of $\mathrm{EDpoly}_{X_\mu^\vee,t}(\epsilon^2)$ are respectively
\begin{align*}
\hspace{-1cm}&(1)\quad a_N(t)=1,\quad a_0(t)=f_\mu^2\quad\mbox{if}\ d=2,\\
\hspace{-1cm}&(2)\quad a_N(t)=f_{\mu,[2]}^{\alpha_{d,d}},\quad a_0(t)=f_\mu^2\cdot f_{\mu,\{1\}}^{-\alpha_{d,1}}\quad\mbox{if}\ d>2\ \mbox{and}\ k=1,\\
\hspace{-1cm}&(3)\quad a_N(t)=f_{\mu,[2]}^{\alpha_{d,d}}\cdot f_{\mu,\{1\}}^{\alpha_{d,k}}\cdot f_{\mu,\{2\}}^{\alpha_{d,d-k}},\quad a_0(t)=f_\mu^2\quad\mbox{if}\ d>2\ \mbox{and}\ 2\le k\le \lfloor d/2 \rfloor.
\end{align*}
In particular, for all\quad $1\le k\le \lfloor d/2 \rfloor$ we have the relation
\begin{equation}\label{eq: equations partitions 2 factors}
(d-k-1)\alpha_{d,k}+(k-1)\alpha_{d,d-k}+\alpha_{d,d}=2(k(d-k)-d+1).
\end{equation}
\end{lemma}

\proof
In the case $(1)$, $\mathrm{EDpoly}_{X_\mu^\vee,t}(\epsilon^2)=\det(\epsilon I-t)\det(\epsilon I+t)$. This is a monic polynomial, and $f_\mu=\det(t)$ is the determinant of the $2\times 2$ matrix representing $t$.

Consider case $(2)$. From Lemma \ref{lem: when are hypersurfaces} we have that $X_{\mu,\{2\}}^\vee$ is not a hypersurface, therefore $f_{\mu,\{2\}}=1$. By Corollary \ref{corol: preservation of exponents} we have that $\alpha_{\mu,[2]}=\alpha_{d,d}$, whereas $\alpha_{\mu,\{1\}}=\alpha_{d,1}$. Therefore equations (\ref{eq: writing the leadcoef}) and (\ref{eq: writing the constterm}) become respectively $a_N(t)=f_{\mu,[2]}^{\alpha_{d,d}}$ and $a_0(t)=f_{\mu}^2\cdot f_{\mu,\{1\}}^{-\alpha_{d,1}}$. Equation (\ref{eq: identity between degrees}) yields the identity
\[
\deg(f_{\mu,[2]})\alpha_{d,d}+\deg(f_{\mu,\{1\}})\alpha_{d,1}-2\deg(f_\mu)+4(d-1)=0.
\]
On one hand, $\deg(f_\mu)=2(d-2)$ by Theorem \ref{thm: degree mu discriminant}. On the other hand, $\deg(f_{\mu,[2]})=4$ and $\deg(f_{\mu,\{1\}})=4(d-2)$ by Corollary \ref{corol: degree varieties f mu J}. Hence we get relation (\ref{eq: equations partitions 2 factors}) for $k=1$.

Consider case $(3)$. Equations (\ref{eq: writing the leadcoef}) and (\ref{eq: writing the constterm}) become respectively $a_N(t)=f_{\mu,[2]}^{\alpha_{\mu,[2]}}\cdot f_{\mu,\{1\}}^{\alpha_{\mu,\{1\}}}\cdot f_{\mu,\{2\}}^{\alpha_{\mu,\{2\}}}$, $a_0(t)=f_{\mu}^2$. Again by Corollary \ref{corol: preservation of exponents} we have that $\alpha_{\mu,[2]}=\alpha_{d,d}$, while $\alpha_{\mu,\{1\}}=\alpha_{d,k}$ and $\alpha_{\mu,\{2\}}=\alpha_{d,d-k}$. Then relation (\ref{eq: equations partitions 2 factors}) follows by (\ref{eq: identity between degrees}) after applying Theorem \ref{thm: degree mu discriminant} and Corollary \ref{corol: degree varieties f mu J}.
\endproof

\begin{lemma}\label{lem: case s=3}
Let $d\ge 3$ and consider the partition $\mu=(k,d-k-1,1)\vdash d$ for all $1\le k\le \lfloor (d-1)/2 \rfloor$. The highest coefficient $a_N=a_N(t)$ and the lowest coefficient $a_0=a_0(t)$ of $\mathrm{EDpoly}_{X_\mu^\vee,t}(\epsilon^2)$ are respectively
\begin{align*}
& (1)\ a_N=f_{\mu,[3]}^{\alpha_{d,d}},\quad a_0=f_\mu^2\cdot f_{\mu,\{1\}}^{-\alpha_{d,1}}\cdot f_{\mu,\{2\}}^{-\alpha_{d,1}}\cdot f_{\mu,\{3\}}^{-\alpha_{d,1}}\quad\mbox{if}\ d=3,\\
& (2)\ a_N=f_{\mu,[3]}^{\alpha_{d,d}}\cdot f_{\mu,\{1,3\}}^{\alpha_{d,2}}\cdot f_{\mu,\{2\}}^{\alpha_{d,d-2}},\quad a_0=f_\mu^2\cdot f_{\mu,\{1\}}^{-\alpha_{d,1}}\cdot f_{\mu,\{3\}}^{-\alpha_{d,1}}\quad\mbox{if}\ d>3\ \mbox{and}\ k=1,\\
& (3)\ a_N=f_{\mu,[3]}^{\alpha_{d,d}}\cdot f_{\mu,\{1,3\}}^{\alpha_{d,k+1}}\cdot f_{\mu,\{2,3\}}^{\alpha_{d,d-k}}\cdot f_{\mu,\{1\}}^{\alpha_{d,k}}\cdot f_{\mu,\{2\}}^{\alpha_{d,d-k-1}},\quad a_0=f_\mu^2\cdot f_{\mu,\{3\}}^{-\alpha_{d,1}}\quad\mbox{if}\ d>3\ \mbox{and}\ 2\le k\le \lfloor \frac{d-1}{2} \rfloor.
\end{align*}
In particular, for all\quad $1\le k\le \lfloor (d-1)/2 \rfloor$ we have the relation
\begin{equation}\label{eq: equations partitions 3 factors}
(d-k-1)\alpha_{d,k}+2(d-k-2)\alpha_{d,k+1}+k\alpha_{d,d-k-1}+2(k-1)\alpha_{d,d-k}+2\alpha_{d,d}=2(3k(d-k-1)-2d+3).
\end{equation}
\end{lemma}

\proof
Consider case $(1)$. By Lemma \ref{lem: when are hypersurfaces}, the varieties $X_{\mu,\{1,2\}}^\vee$, $X_{\mu,\{1,3\}}^\vee$ and $X_{\mu,\{2,3\}}^\vee$ are not hypersurfaces, therefore $f_{\mu,\{1,2\}}=f_{\mu,\{1,3\}}=f_{\mu,\{2,3\}}=1$. Moreover, by Corollary \ref{corol: preservation of exponents} we have $\alpha_{\mu,[3]}=\alpha_{d,d}$ and $\alpha_{\mu,\{1\}}=\alpha_{\mu,\{2\}}=\alpha_{\mu,\{3\}}=\alpha_{d,1}$. We refer the reader to Section \ref{sec: 222} for a detailed treatise on this specific example.

Now consider case $(2)$. By Lemma \ref{lem: when are hypersurfaces}, the varieties $X_{\mu,\{1,2\}}^\vee$ and $X_{\mu,\{2,3\}}^\vee$ are not hypersurfaces, hence $f_{\mu,\{1,2\}}=f_{\mu,\{2,3\}}=1$. Moreover, by Corollary \ref{corol: preservation of exponents} we have that $\alpha_{\mu,[3]}=\alpha_{d,d}$, $\alpha_{\mu,\{1,3\}}=\alpha_{d,2}$, $\alpha_{\mu,\{2\}}=\alpha_{d,d-2}$ and $\alpha_{\mu,\{1\}}=\alpha_{\mu,\{3\}}=\alpha_{d,1}$. By Corollary \ref{corol: degree varieties f mu J} we have that $\deg(f_{\mu,[3]})=8$, $\deg(f_{\mu,\{1,3\}})=8(d-3)$, $\deg(f_{\mu,\{2\}})=4$, $\deg(f_{\mu,\{1\}})=4(d-2)$, $\deg(f_{\mu,\{3\}})=4(d-2)$ and finally $\deg(f_\mu)=4(2d-5)$. Putting all this information together, we obtain equation (\ref{eq: equations partitions 3 factors}) in the case $k=1$.

Finally consider case $(3)$. Analogously to case $(2)$, we have that the variety $X_{\mu,\{1,2\}}^\vee$ is not a hypersurface, hence $f_{\mu,\{1,2\}}=1$. In addition, $\alpha_{\mu,[3]}=\alpha_{d,d}$, $\alpha_{\mu,\{1,3\}}=\alpha_{d,k+1}$, $\alpha_{\mu,\{2,3\}}=\alpha_{d,d-k}$, $\alpha_{\mu,\{1\}}=\alpha_{d,k}$, $\alpha_{\mu,\{2\}}=\alpha_{d,d-k-1}$ and $\alpha_{\mu,\{3\}}=\alpha_{d,1}$. Moreover, $\deg(f_{\mu,[3]})=8$, $\deg(f_{\mu,\{1,3\}})=8(d-k-2)$, $\deg(f_{\mu,\{2,3\}})=8(k-1)$, $\deg(f_{\mu,\{1\}})=4(d-k-1)$, $\deg(f_{\mu,\{2\}})=4k$, $\deg(f_{\mu,\{3\}})=4(3k(d-k-1)-2d+4)$ and finally $\deg(f_\mu)=12k(d-k-1)-4d+4$. Thus we obtain equation (\ref{eq: equations partitions 3 factors}).
\endproof

\begin{proposition}\label{pro: linear system}
The linear system $\mathcal{S}_d$ defined by equations (\ref{eq: equation partition 1 factor}), (\ref{eq: equations partitions 2 factors}) and (\ref{eq: equations partitions 3 factors}) admits the only solution $(\alpha_{d,1},\ldots,\alpha_{d,d})=(-1,0,1,\ldots,d-2)$, provided that $\alpha_{d,d-1}\coloneqq d-3$.
\end{proposition}
\proof
The linear system $\mathcal{S}_d$ has $1+\lfloor d/2 \rfloor+\lfloor (d-1)/2 \rfloor=d$ equations in $d$ unknowns $\alpha_{d,1},\ldots,\alpha_{d,d}$. Observe that no equation involves the unknown $\alpha_{d-1}$, hence the rank of the matrix of coefficients of $\mathcal{S}_d$ is at most $d-1$. The geometrical reason is that, for any partition $\mu=(\mu_1,\ldots,\mu_s)\vdash d$, with $\mu_1\ge\cdots\ge\mu_s$, the only subset $J\subset[s]$ such that $\alpha_{\mu,J}=\alpha_{d,d-1}$ is $J=[s-1]$, by Corollary \ref{corol: preservation of exponents}. Indeed, by Lemma \ref{lem: when are hypersurfaces}, the corresponding dual variety $X_{\mu,J}^\vee$ is not a hypersurface, thus $f_{\mu,J}=1$ and the exponent $\alpha_{d,d-1}$ remains undetermined. To be consistent with the higher dimensional results stated in \cite{Sod2}, we define $\alpha_{d,d-1}\coloneqq d-3$. After substituting in (\ref{eq: equation partition 1 factor}), (\ref{eq: equations partitions 2 factors}) and (\ref{eq: equations partitions 3 factors}), we see that the vector $(\alpha_{d,1},\ldots,\alpha_{d,d})=(-1,0,1,\ldots,d-2)$ is a solution of $\mathcal{S}_d$.

It remains to show that the matrix of coefficients of $\mathcal{S}_d$ is of maximal rank $d-1$. First of all, one clearly checks that the set of equations coming from (\ref{eq: equations partitions 2 factors}) are pairwise linearly independent. The same holds for the set of equations coming from (\ref{eq: equations partitions 2 factors}). Moreover, any equation coming from either (\ref{eq: equations partitions 2 factors}) or (\ref{eq: equations partitions 3 factors}) is linearly independent with (\ref{eq: equation partition 1 factor}). It may happen that the $k$-th equation in (\ref{eq: equations partitions 3 factors}) is a linear combination of the $k$-th and $(k+1)$-th equations in (\ref{eq: equations partitions 2 factors}). This happens only if the maximal minors of the submatrix
\[
\mathcal{M}_{d,k}=
\begin{pmatrix}
d-k-1 & 0 & 0 & k-1\\
0 & d-k-2 & k & 0\\
d-k-1 & 2(d-k-2) & k & 2(k-1)
\end{pmatrix}
\]
obtained extracting the coefficients of $\alpha_k$, $\alpha_{k+1}$, $\alpha_{d-k-1}$ and $\alpha_{d-k}$ from the three mentioned equations, vanish simultaneously. One may easily verify that this is impossible for $1\le k\le \lfloor (d-1)/2 \rfloor$, because the maximal minors $m_{d,k}^{(j_1,j_2,j_3)}$ obtained picking the columns $j_1,j_2,j_3$ of $\mathcal{M}_{d,k}$ are respectively
\begin{align*}
m_{d,k}^{(1,2,3)} &= -k\left(d-k-2\right)\left(d-k-1\right),
& m_{d,k}^{(1,3,4)} &= \left(k-1\right)\left(d-k-2\right)\left(d-k-1\right),\\
m_{d,k}^{(1,2,4)} &= k\left(k-1\right)\left(d-k-1\right),
& m_{d,k}^{(2,3,4)} &= -k\left(k-1\right)\left(d-k-2\right).\qedhere
\end{align*}
\endproof

\proof[Proof of the Main Theorem]
Let $d\ge 1$ be an integer and let $\mu=(\mu_1,\ldots,\mu_s)$ be a partition of $d$. We only need to show that the highest coefficient $a_N=a_N(t)$ and the lowest coefficient $a_0=a_0(t)$ of the ED polynomial of $X_{\mu}^\vee$ at $t\in S^\mu V$ are respectively
\begin{equation}\label{eq: final expression extreme coefficients}
a_N=\prod_{j\colon \mu_j>1}f_{\mu,\{j\}}^{\mu_j-2}\cdot\prod_{|J|>1}f_{\mu,J}^{\sum_{k\in J}\mu_k-2},\quad a_0=f_{\mu}^2\cdot\prod_{j\colon \mu_j=1}f_{\mu,\{j\}}.
\end{equation}
By Corollary \ref{corol: preservation of exponents}, for every non-empty subset $J\subset[s]$, the exponent of $f_{\mu,J}$ is $\alpha_{\mu,J}=\alpha_{d,\sum_{k\in J}\mu_k}$. Moreover, by Proposition \ref{pro: linear system}, we have that $\alpha_{d,\sum_{k\in J}\mu_k}=\sum_{k\in J}\mu_k-2$, thus completing the proof.
\endproof

\begin{remark}
As one may foresee from the statement of the Main Theorem and the different steps of its proof, this result could be extended to tensors of format $n^d$ for $n\ge 2$ or even of format $n_1\times\cdots\times n_d$ for $n_1,\ldots,n_d\ge 2$. Indeed, in this paper we focus on tensors of binary format, and leave to the paper in preparation \cite{Sod2} other more general considerations on ED polynomials of Segre-Veronese varieties.
\end{remark}

An immediate consequence of the Main Theorem is that, for any partition $\mu\vdash d$, we may write an identity involving the ED degree of $X_\mu$ and the degrees of the varieties $X_{\mu,J}^\vee$, as pointed out below.
\begin{corollary} Consider a partition $\mu=(\mu_1,\ldots,\mu_s)\vdash d$. Let $\mu(J)$ be the partition defined in Corollary \ref{corol: degree varieties f mu J} for all $J\subset[s]$. Recalling that the degrees of the $\mu$-discriminant $\delta_{\mu}$ and of the hyperdeterminant $\delta_j$ are defined in (\ref{eq: degree mu discriminant}) and (\ref{eq: degree hyperdeterminant}), respectively, then
\begin{align*}\label{eq: identity degrees}
\mathrm{EDdegree}(X_\mu)&=s!\mu_1\cdots\mu_s=\sum_{J\subset[s]}2^{|J|-1}\left(2-\sum_{k\in J}\mu_k\right)\delta_{\mu(J)},\\
\mathrm{EDdegree}(X_d)&=d!=\sum_{j=0}^d\binom{d}{j}(2-j)2^{j-1}\delta_{d-j}.
\end{align*}
\end{corollary}

\section{On the non-negativity of the polynomials $f_{\mu,J}$}\label{sec: SOS}

This section investigates the non-negativity of the various factors $f_{\mu,J}$ appearing in the extreme coefficients of $\mathrm{EDpoly}_{X_\mu^\vee,t}(\epsilon^2)$. Actually, they are (products of) SOS polynomials. This fact is useful for the considerations about tensors of format $2\times 2\times 2$ made in Section \ref{sec: 222}.

\begin{proposition}\label{SOS}
Let $J\subset [s]$, $J\neq\emptyset$. If $J=[s]$, then $f_{\mu,J}$ is the product of $d$ SOS polynomials. If $J\neq[s]$, then $f_{\mu,J}$ is a SOS polynomial. In particular, $f_{\mu,J}$ is a non-negative polynomial for all $J\neq\emptyset$.
\end{proposition}
\proof
For $J=[s]$, the thesis follows since $X_{\mu,J}^\vee\subset\PP(V^{\otimes d})$ is the union of $d$ pairwise conjugate hyperplanes. Now let $\mu=1^d$ and $J=\{1\}$. More explicitly,
\begin{equation}\label{eq: two components Segre}
X_{d,\{1\}}^\vee=\left[\mathrm{Seg}([(1,\sqrt{-1})]\times\PP(V)^{\times(d-1)})\right]^\vee\cup\left[\mathrm{Seg}([(1,-\sqrt{-1})]\times\PP(V)^{\times(d-1)})\right]^\vee\subset\PP(V^{\otimes d}).
\end{equation}
Thus $X_{d,\{1\}}^\vee$ is isomorphic to two copies of $X_{d-1}^\vee\subset\PP(V^{\otimes(d-1)})$. By Lemma \ref{lem: when are hypersurfaces}, the varieties $X_{d,\{1\}}^\vee$ and $X_{d-1}^\vee$ are hypersurfaces when $d\ge 3$. If $\{a_{i_1\cdots i_d}\}$ and $\{b_{j_2\cdots j_d}\}$ are homogeneous coordinates for $\PP(V^{\otimes d})$ and $\PP(V^{\otimes(d-1)})$ respectively, the equations of the two components of $X_{d,\{1\}}^\vee$ in (\ref{eq: two components Segre}) are 
\begin{align}\label{eq: f + -}
\begin{split}
f_{d-1}^{+}(\{a_{j_1\cdots j_d}\})&\coloneqq f_{d-1}(\{b_{j_2\cdots j_d}\})_{\mkern 1mu \vrule height 2ex\mkern2mu \{b_{j_2\cdots j_d}\ =\ a_{0i_2\cdots i_d}+\sqrt{-1}a_{1i_2\cdots i_d}\}},\\
f_{d-1}^{-}(\{a_{j_1\cdots j_d}\})&\coloneqq f_{d-1}(\{b_{j_2\cdots j_d}\})_{\mkern 1mu \vrule height 2ex\mkern2mu \{b_{j_2\cdots j_d}\ =\ a_{0i_2\cdots i_d}-\sqrt{-1}a_{1i_2\cdots i_d}\}}.
\end{split}
\end{align}
In particular, $f_{d-1}^{+}$ and $f_{d-1}^{-}$ are conjugate polynomials and their product is the equation $f_{d,\{1\}}$ of $X_{d,\{1\}}^\vee$. Therefore $f_{d,\{1\}}$ is the sum of two squared polynomials. In the same fashion, we show that $f_{d,\{1,2\}}=f_{d,\{2\}}^{+}\cdot f_{d,\{2\}}^{-}$, where the factors $f_{d,\{2\}}^{+}$ and $f_{d,\{2\}}^{-}$ are defined as in (\ref{eq: f + -}). Therefore, $f_{d,\{1,2\}}$ is again a sum of two squared polynomials. More in general, the iteration of this argument shows that, possibly after a permutation of the indices, the polynomial $f_{d,J}$ is a sum of two squared polynomials for any subset $J\subset[d]$.

Now consider a partition $\mu=(\mu_1,\ldots,\mu_s)\vdash d$ and a non-empty subset $J\subset[s]$. On one hand, By Proposition \ref{pro: symmetrization into two parts} there exists a subset $\tilde{J}\subset[d]$ with $|\tilde{J}|=\sum_{j\in J}\mu_j$ such that $f_{\mu,J}$ divides $\pi_{\mu,1^d}(f_{d,\tilde{J}})$ with multiplicity one. On the other hand, the first part of the proof implies that $f_{d,\tilde{J}}=h_1^2+h_2^2$ for some homogeneous polynomials $h_1$ and $h_2$. Summing up, there exist two homogeneous polynomials $h_1'$ and $h_2'$ such that $h_1'\pm\sqrt{-1}h_2'$ divides $\pi_{\mu,1^d}(h_1\pm\sqrt{-1}h_2)$ with multiplicity one and $f_{\mu,J}=h_1'^2+h_2'^2=(h_1'+\sqrt{-1}h_2')(h_1'-\sqrt{-1}h_2')$, thus completing the proof.
\endproof

\begin{problem}
In particular, Proposition \ref{SOS} tells us that, if $(X_\mu\cap\widetilde{Q})^\vee$ is a hypersurface, then its equation is a SOS polynomial. More in general, it would be interesting to solve the following problem suggested by Bernd Sturmfels: assuming that the variety $(X^\vee\cap Q)^\vee$ appearing in \cite[Proposition 5.4]{OS} is a hypersurface, is its equation a non-negative polynomial? If so, is it a SOS polynomial?
\end{problem}

\begin{remark}\label{rmk: sos constterm pictures}
Looking closely at the polynomials defined in (\ref{eq: f + -}), one may see that for all $d\ge 3$ and for all $1\le j\le d$, the polynomial $f_{d,\{j\}}$ may be written as
\begin{equation}\label{eq: product hyperdeterminants}
f_{d,\{j\}}(t)=\mathrm{Det}(t_j^{(0)}+\sqrt{-1}t_j^{(1)})\cdot\mathrm{Det}(t_j^{(0)}-\sqrt{-1}t_j^{(1)}),
\end{equation}
where $t_j^{(0)}$ and $t_j^{(1)}$ are the tensors in $V^{\otimes(d-1)}$ obtained considering in $t$ the slices $\{t_{i_1\cdots i_d}\}$ with $i_j=0$ and $i_j=1$, respectively. The cases $d=3$ and $d=4$ are depicted in Figures \ref{fig: case d=3 factors constant term} and \ref{fig: case d=4 factors constant term}, respectively.

\begin{figure}[ht]
\begin{center}
\begin{minipage}[c]{.3\textwidth}
\centering
\begin{tikzpicture}[line cap=round,line join=round,>=triangle 45,x=1.0cm,y=1.0cm, scale=0.3]
\clip(-4.5,-4.1) rectangle (4.5,4.35);
\fill[line width=0.5pt,fill=black,fill opacity=0.1] (-0.541995,0.6821570000000001) -- (1.92516,1.1486005000000001) -- (1.92516,-1.5141485) -- (-0.541995,-1.980592) -- cycle;
\fill[line width=0.5pt,fill=black,fill opacity=0.1] (-1.92516,-1.1486) -- (0.5419949999999997,-0.6821565000000002) -- (0.5419950000000003,1.9805925) -- (-1.92516,1.514149) -- cycle;
\fill[line width=0.5pt,fill=black,fill opacity=0.1] (-1.92516,-1.1486) -- (-0.541995,-1.980592) -- (-0.541995,0.6821570000000001) -- (-1.92516,1.514149) -- cycle;
\fill[line width=0.5pt,fill=black,fill opacity=0.1] (0.5419950000000003,1.9805925) -- (1.92516,1.1486005000000001) -- (1.92516,-1.5141485) -- (0.5419949999999997,-0.6821565000000002) -- cycle;
\fill[line width=0.5pt,fill=black,fill opacity=0.1] (-0.541995,-1.980592) -- (1.92516,-1.5141485) -- (0.5419949999999997,-0.6821565000000002) -- (-1.92516,-1.1486) -- cycle;
\fill[line width=0.5pt,fill=black,fill opacity=0.1] (-1.92516,1.514149) -- (0.5419950000000003,1.9805925) -- (1.92516,1.1486005000000001) -- (-0.541995,0.6821570000000001) -- cycle;

\fill[line width=0.5pt,fill=black,fill opacity=0.1] (-3.85032,3.028298) -- (-1.08399,1.3643140000000002) -- (3.85032,2.2972010000000003) -- (1.08399,3.961185) -- (-3.85032,3.028298) -- cycle;
\fill[line width=0.5pt,fill=black,fill opacity=0.1] (-3.85032,-2.2972) -- (-1.08399,-3.961184) -- (3.85032,-3.028297) -- (1.08399,-1.3643130000000003) -- (-3.85032,-2.2972) -- cycle;

\fill[line width=0.5pt,fill=black,fill opacity=0.1] (-3.85032,3.028298) -- (-1.08399,1.3643140000000002) -- (-1.08399,-3.961184) -- (-3.85032,-2.2972) -- (-3.85032,3.028298) -- cycle;
\fill[line width=0.5pt,fill=black,fill opacity=0.1] (1.08399,3.961185) -- (3.85032,2.2972010000000003) -- (3.85032,-3.028297) -- (1.08399,-1.3643130000000003) -- (1.08399,3.961185) -- cycle;

\fill[line width=0.5pt,fill=black,fill opacity=0.1] (-1.08399,1.3643140000000002) -- (3.85032,2.2972010000000003) -- (3.85032,-3.028297) -- (-1.08399,-3.961184) -- (-1.08399,1.3643140000000002) -- cycle;
\fill[line width=0.5pt,fill=black,fill opacity=0.1] (-3.85032,3.028298) -- (1.08399,3.961185) -- (1.08399,-1.3643130000000003) -- (-3.85032,-2.2972) -- (-3.85032,3.028298) -- cycle;

\draw [line width=0.5pt] (3.85032,2.2972010000000003)-- (-1.08399,1.3643140000000002);
\draw [line width=0.5pt] (-1.08399,1.3643140000000002)-- (-1.08399,-3.961184);
\draw [line width=0.5pt] (-1.92516,1.514149)-- (-1.08399,1.008174144765809);
\draw [line width=0.5pt,dotted] (-1.08399,1.008174144765809)-- (-0.541995,0.6821570000000001);
\draw [line width=0.5pt] (0.5419950000000003,1.9805925)-- (0.9326826447970984,1.7455887260684537);
\draw [line width=0.5pt,dotted] (0.9326826447970984,1.7455887260684537)-- (1.92516,1.1486005000000001);
\draw [line width=0.5pt] (-1.92516,-1.1486)-- (-1.08399,-1.6545748552341912);
\draw [line width=0.5pt,dotted] (-1.08399,-1.6545748552341912)-- (-0.541995,-1.980592);
\draw [line width=0.5pt,dotted] (-1.92516,1.514149)-- (-1.4746776447970984,1.5993177739315465);
\draw [line width=0.5pt] (-1.4746776447970984,1.5993177739315465)-- (0.5419950000000003,1.9805925);
\draw [line width=0.5pt] (3.85032,2.2972010000000003)-- (3.85032,-3.028297);
\draw [line width=0.5pt] (3.85032,-3.028297)-- (-1.08399,-3.961184);
\draw [line width=0.5pt] (-0.541995,-1.980592)-- (-1.08399,-3.961184);
\draw [line width=0.5pt] (-1.08399,-3.961184)-- (-3.85032,-2.2972);
\draw [line width=0.5pt] (-3.85032,3.028298)-- (-3.85032,-2.2972);
\draw [line width=0.5pt] (-3.85032,3.028298)-- (-1.08399,1.3643140000000002);
\draw [line width=0.5pt] (-1.08399,1.3643140000000002)-- (-0.541995,0.6821570000000001);
\draw [line width=0.5pt] (-3.85032,3.028298)-- (-1.92516,1.514149);
\draw [line width=0.5pt] (1.08399,3.961185)-- (-3.85032,3.028298);
\draw [line width=0.5pt] (1.08399,3.961185)-- (3.85032,2.2972010000000003);
\draw [line width=0.5pt] (3.85032,2.2972010000000003)-- (1.92516,1.1486005000000001);
\draw [line width=0.5pt, dotted] (-3.85032,-2.2972)-- (1.08399,-1.3643130000000003);
\draw [line width=0.5pt] (-1.92516,-1.1486)-- (-3.85032,-2.2972);
\draw [line width=0.5pt] (1.92516,-1.5141485)-- (3.85032,-3.028297);
\draw [line width=0.5pt, dotted] (3.85032,-3.028297)-- (1.08399,-1.3643130000000003);
\draw [line width=0.5pt, dotted] (-1.92516,-1.1486)-- (0.5419949999999997,-0.6821565000000002);
\draw [line width=0.5pt] (-1.92516,1.514149)-- (-1.92516,-1.1486);
\draw [line width=0.5pt] (1.92516,1.1486005000000001)-- (1.92516,-1.5141485);
\draw [line width=0.5pt, dotted] (1.92516,-1.5141485)-- (0.5419949999999997,-0.6821565000000002);
\draw [line width=0.5pt, dotted] (0.5419949999999997,-0.6821565000000002)-- (0.5419950000000003,1.9805925);
\draw [line width=0.5pt] (0.5419950000000003,1.9805925)-- (1.08399,3.961185);
\draw [line width=0.5pt, dotted] (1.08399,3.961185)-- (1.08399,-1.3643130000000003);
\draw [line width=0.5pt, dotted] (1.08399,-1.3643130000000003)-- (0.5419949999999997,-0.6821565000000002);
\draw [line width=0.5pt, dotted] (0.5419949999999997,-0.6821565000000002)-- (0.5419950000000003,1.9805925);
\draw [line width=0.5pt] (-0.541995,0.6821570000000001)-- (1.92516,1.1486005000000001);
\draw [line width=0.5pt] (1.92516,1.1486005000000001)-- (1.92516,-1.5141485);
\draw [line width=0.5pt] (1.92516,-1.5141485)-- (-0.541995,-1.980592);
\draw [line width=0.5pt] (-0.541995,-1.980592)-- (-0.541995,0.6821570000000001);
\begin{tiny}
\draw [fill=black] (-1.92516,-1.1486) circle (2pt);
\draw [fill=black] (-0.541995,-1.980592) circle (2pt);
\draw [fill=black] (1.92516,-1.5141485) circle (2pt);
\draw [fill=black] (-1.92516,1.514149) circle (2pt);
\draw [fill=black] (-0.541995,0.6821570000000001) circle (2pt);
\draw [fill=black] (1.92516,1.1486005000000001) circle (2pt);
\draw [fill=black] (0.5419950000000003,1.9805925) circle (2pt);
\draw [fill=black] (0.5419949999999997,-0.6821565000000002) circle (2pt);
\draw [fill=black] (-3.85032,-2.2972) circle (2pt);
\draw [fill=black] (-1.08399,-3.961184) circle (2pt);
\draw [fill=black] (3.85032,-3.028297) circle (2pt);
\draw [fill=black] (-3.85032,3.028298) circle (2pt);
\draw [fill=black] (-1.08399,1.3643140000000002) circle (2pt);
\draw [fill=black] (3.85032,2.2972010000000003) circle (2pt);
\draw [fill=black] (1.08399,3.961185) circle (2pt);
\draw [fill=black] (1.08399,-1.3643130000000003) circle (2pt);
\draw[color=black] (0.2,0) node {$t_1^{(1)}$};
\draw[color=black] (-3.75,3.75) node {$t_1^{(0)}$};
\end{tiny}
\end{tikzpicture}
\end{minipage}
\hspace{-60pt}
\begin{minipage}[c]{.3\textwidth}
\centering
\begin{tikzpicture}[line cap=round,line join=round,>=triangle 45,x=1.0cm,y=1.0cm, scale=0.3]
\clip(-4.5,-4.1) rectangle (4.5,4.35);

\fill[line width=0.5pt,fill=black,fill opacity=0] (-1.08399,1.3643140000000002) -- (-0.541995,0.6821570000000001) -- (1.92516,1.1486005000000001) -- (3.85032,2.2972010000000003) -- cycle;
\fill[line width=0.5pt,fill=black,fill opacity=0] (-0.541995,0.6821570000000001) -- (1.92516,1.1486005000000001) -- (1.92516,-1.5141485) -- (-0.541995,-1.980592) -- cycle;
\fill[line width=0.5pt,fill=black,fill opacity=0] (-0.541995,-1.980592) -- (1.92516,-1.5141485) -- (3.85032,-3.028297) -- (-1.08399,-3.961184) -- cycle;

\fill[line width=0.5pt,fill=black,fill opacity=0] (-3.85032,3.028298) -- (-1.92516,1.514149) -- (0.5419950000000003,1.9805925) -- (1.08399,3.961185) -- cycle;
\fill[line width=0.5pt,fill=black,fill opacity=0] (-3.85032,-2.2972) -- (1.08399,-1.3643130000000003) -- (0.5419949999999997,-0.6821565000000002) -- (-1.92516,-1.1486) -- cycle;
\fill[line width=0.5pt,fill=black,fill opacity=0] (-1.92516,-1.1486) -- (0.5419949999999997,-0.6821565000000002) -- (0.5419950000000003,1.9805925) -- (-1.92516,1.514149) -- cycle;

\fill[line width=0.5pt,fill=black,fill opacity=0.1] (-1.92516,-1.1486) -- (-1.92516,1.514149) -- (-3.85032,3.028298) -- (-3.85032,-2.2972) -- cycle;
\fill[line width=0.5pt,fill=black,fill opacity=0.1] (-1.08399,1.3643140000000002) -- (-0.541995,0.6821570000000001) -- (-0.541995,-1.980592) -- (-1.08399,-3.961184) -- cycle;
\fill[line width=0.5pt,fill=black,fill opacity=0.1] (-1.92516,-1.1486) -- (-0.541995,-1.980592) -- (-1.08399,-3.961184) -- (-3.85032,-2.2972) -- cycle;
\fill[line width=0.5pt,fill=black,fill opacity=0.1] (-1.92516,-1.1486) -- (-0.541995,-1.980592) -- (-0.541995,0.6821570000000001) -- (-1.92516,1.514149) -- cycle;
\fill[line width=0.5pt,fill=black,fill opacity=0.1] (-3.85032,3.028298) -- (-1.08399,1.3643140000000002) -- (-0.541995,0.6821570000000001) -- (-1.92516,1.514149) -- cycle;

\fill[line width=0.5pt,fill=black,fill opacity=0.1] (1.92516,1.1486005000000001) -- (3.85032,2.2972010000000003) -- (3.85032,-3.028297) -- (1.92516,-1.5141485) -- cycle;
\fill[line width=0.5pt,fill=black,fill opacity=0.1] (0.5419950000000003,1.9805925) -- (1.08399,3.961185) -- (1.08399,-1.3643130000000003) -- (0.5419949999999997,-0.6821565000000002) -- cycle;
\fill[line width=0.5pt,fill=black,fill opacity=0.1] (0.5419950000000003,1.9805925) -- (1.08399,3.961185) -- (3.85032,2.2972010000000003) -- (1.92516,1.1486005000000001) -- cycle;
\fill[line width=0.5pt,fill=black,fill opacity=0.1] (0.5419949999999997,-0.6821565000000002) -- (1.92516,-1.5141485) -- (3.85032,-3.028297) -- (1.08399,-1.3643130000000003) -- cycle;
\fill[line width=0.5pt,fill=black,fill opacity=0.1] (0.5419950000000003,1.9805925) -- (1.92516,1.1486005000000001) -- (1.92516,-1.5141485) -- (0.5419949999999997,-0.6821565000000002) -- cycle;

\fill[line width=0.5pt,fill=black,fill opacity=0.1] (-3.85032,3.028298) -- (-1.08399,1.3643140000000002) -- (-1.08399,-3.961184) -- (-3.85032,-2.2972) -- (-3.85032,3.028298) -- cycle;
\fill[line width=0.5pt,fill=black,fill opacity=0.1] (1.08399,3.961185) -- (3.85032,2.2972010000000003) -- (3.85032,-3.028297) -- (1.08399,-1.3643130000000003) -- (1.08399,3.961185) -- cycle;

\fill[line width=0.5pt,fill=black,fill opacity=0] (-1.92516,1.514149) -- (0.5419950000000003,1.9805925) -- (1.92516,1.1486005000000001) -- (-0.541995,0.6821570000000001) -- cycle;
\fill[line width=0.5pt,fill=black,fill opacity=0] (-0.541995,-1.980592) -- (1.92516,-1.5141485) -- (0.5419949999999997,-0.6821565000000002) -- (-1.92516,-1.1486) -- cycle;
\fill[line width=0.5pt,fill=black,fill opacity=0] (-0.541995,0.6821570000000001) -- (1.92516,1.1486005000000001) -- (1.92516,-1.5141485) -- (-0.541995,-1.980592) -- cycle;

\draw [line width=0.5pt] (3.85032,2.2972010000000003)-- (-1.08399,1.3643140000000002);
\draw [line width=0.5pt] (-1.08399,1.3643140000000002)-- (-1.08399,-3.961184);
\draw [line width=0.5pt] (-1.92516,1.514149)-- (-1.08399,1.008174144765809);
\draw [line width=0.5pt,dotted] (-1.08399,1.008174144765809)-- (-0.541995,0.6821570000000001);
\draw [line width=0.5pt] (0.5419950000000003,1.9805925)-- (0.9326826447970984,1.7455887260684537);
\draw [line width=0.5pt,dotted] (0.9326826447970984,1.7455887260684537)-- (1.92516,1.1486005000000001);
\draw [line width=0.5pt] (-1.92516,-1.1486)-- (-1.08399,-1.6545748552341912);
\draw [line width=0.5pt,dotted] (-1.08399,-1.6545748552341912)-- (-0.541995,-1.980592);
\draw [line width=0.5pt,dotted] (-1.92516,1.514149)-- (-1.4746776447970984,1.5993177739315465);
\draw [line width=0.5pt] (-1.4746776447970984,1.5993177739315465)-- (0.5419950000000003,1.9805925);
\draw [line width=0.5pt] (3.85032,2.2972010000000003)-- (3.85032,-3.028297);
\draw [line width=0.5pt] (3.85032,-3.028297)-- (-1.08399,-3.961184);
\draw [line width=0.5pt] (-0.541995,-1.980592)-- (-1.08399,-3.961184);
\draw [line width=0.5pt] (-1.08399,-3.961184)-- (-3.85032,-2.2972);
\draw [line width=0.5pt] (-3.85032,3.028298)-- (-3.85032,-2.2972);
\draw [line width=0.5pt] (-3.85032,3.028298)-- (-1.08399,1.3643140000000002);
\draw [line width=0.5pt] (-1.08399,1.3643140000000002)-- (-0.541995,0.6821570000000001);
\draw [line width=0.5pt] (-3.85032,3.028298)-- (-1.92516,1.514149);
\draw [line width=0.5pt] (1.08399,3.961185)-- (-3.85032,3.028298);
\draw [line width=0.5pt] (1.08399,3.961185)-- (3.85032,2.2972010000000003);
\draw [line width=0.5pt] (3.85032,2.2972010000000003)-- (1.92516,1.1486005000000001);
\draw [line width=0.5pt, dotted] (-3.85032,-2.2972)-- (1.08399,-1.3643130000000003);
\draw [line width=0.5pt] (-1.92516,-1.1486)-- (-3.85032,-2.2972);
\draw [line width=0.5pt] (1.92516,-1.5141485)-- (3.85032,-3.028297);
\draw [line width=0.5pt, dotted] (3.85032,-3.028297)-- (1.08399,-1.3643130000000003);
\draw [line width=0.5pt, dotted] (-1.92516,-1.1486)-- (0.5419949999999997,-0.6821565000000002);
\draw [line width=0.5pt] (-1.92516,1.514149)-- (-1.92516,-1.1486);
\draw [line width=0.5pt] (1.92516,1.1486005000000001)-- (1.92516,-1.5141485);
\draw [line width=0.5pt, dotted] (1.92516,-1.5141485)-- (0.5419949999999997,-0.6821565000000002);
\draw [line width=0.5pt, dotted] (0.5419949999999997,-0.6821565000000002)-- (0.5419950000000003,1.9805925);
\draw [line width=0.5pt] (0.5419950000000003,1.9805925)-- (1.08399,3.961185);
\draw [line width=0.5pt, dotted] (1.08399,3.961185)-- (1.08399,-1.3643130000000003);
\draw [line width=0.5pt, dotted] (1.08399,-1.3643130000000003)-- (0.5419949999999997,-0.6821565000000002);
\draw [line width=0.5pt, dotted] (0.5419949999999997,-0.6821565000000002)-- (0.5419950000000003,1.9805925);
\draw [line width=0.5pt] (-0.541995,0.6821570000000001)-- (1.92516,1.1486005000000001);
\draw [line width=0.5pt] (1.92516,1.1486005000000001)-- (1.92516,-1.5141485);
\draw [line width=0.5pt] (1.92516,-1.5141485)-- (-0.541995,-1.980592);
\draw [line width=0.5pt] (-0.541995,-1.980592)-- (-0.541995,0.6821570000000001);
\begin{tiny}
\draw [fill=black] (-1.92516,-1.1486) circle (2.5pt);
\draw [fill=black] (-0.541995,-1.980592) circle (2.5pt);
\draw [fill=black] (1.92516,-1.5141485) circle (2.5pt);
\draw [fill=black] (-1.92516,1.514149) circle (2.0pt);
\draw [fill=black] (-0.541995,0.6821570000000001) circle (2.0pt);
\draw [fill=black] (1.92516,1.1486005000000001) circle (2.0pt);
\draw [fill=black] (0.5419950000000003,1.9805925) circle (2.0pt);
\draw [fill=black] (0.5419949999999997,-0.6821565000000002) circle (2.0pt);
\draw [fill=black] (-3.85032,-2.2972) circle (2.5pt);
\draw [fill=black] (-1.08399,-3.961184) circle (2.5pt);
\draw [fill=black] (3.85032,-3.028297) circle (2.5pt);
\draw [fill=black] (-3.85032,3.028298) circle (2.0pt);
\draw [fill=black] (-1.08399,1.3643140000000002) circle (2.0pt);
\draw [fill=black] (3.85032,2.2972010000000003) circle (2.0pt);
\draw [fill=black] (1.08399,3.961185) circle (2.0pt);
\draw [fill=black] (1.08399,-1.3643130000000003) circle (2.0pt);
\draw[color=black] (-2.9,-1) node {$t_2^{(0)}$};
\draw[color=black] (3,-1.3) node {$t_2^{(1)}$};
\end{tiny}
\end{tikzpicture}
\end{minipage}  
\hspace{-60pt}
\begin{minipage}[c]{.3\textwidth}
\centering
\begin{tikzpicture}[line cap=round,line join=round,>=triangle 45,x=1.0cm,y=1.0cm, scale=0.3]
\clip(-4.5,-4.1) rectangle (4.5,4.35);
\fill[line width=0.5pt,fill=black,fill opacity=0.1] (-3.85032,3.028298) -- (-1.08399,1.3643140000000002) -- (3.85032,2.2972010000000003) -- (1.08399,3.961185) -- (-3.85032,3.028298) -- cycle;
\fill[line width=0.5pt,fill=black,fill opacity=0.1] (-1.08399,1.3643140000000002) -- (-0.541995,0.6821570000000001) -- (1.92516,1.1486005000000001) -- (3.85032,2.2972010000000003) -- cycle;
\fill[line width=0.5pt,fill=black,fill opacity=0.1] (-3.85032,3.028298) -- (-1.92516,1.514149) -- (0.5419950000000003,1.9805925) -- (1.08399,3.961185) -- cycle;
\fill[line width=0.5pt,fill=black,fill opacity=0.1] (-3.85032,3.028298) -- (-1.08399,1.3643140000000002) -- (-0.541995,0.6821570000000001) -- (-1.92516,1.514149) -- cycle;
\fill[line width=0.5pt,fill=black,fill opacity=0.1] (0.5419950000000003,1.9805925) -- (1.08399,3.961185) -- (3.85032,2.2972010000000003) -- (1.92516,1.1486005000000001) -- cycle;
\fill[line width=0.5pt,fill=black,fill opacity=0.1] (-1.92516,1.514149) -- (0.5419950000000003,1.9805925) -- (1.92516,1.1486005000000001) -- (-0.541995,0.6821570000000001) -- cycle;

\fill[line width=0.5pt,fill=black,fill opacity=0.1] (-0.541995,-1.980592) -- (1.92516,-1.5141485) -- (0.5419949999999997,-0.6821565000000002) -- (-1.92516,-1.1486) -- cycle;
\fill[line width=0.5pt,fill=black,fill opacity=0.1] (-1.92516,-1.1486) -- (-0.541995,-1.980592) -- (-1.08399,-3.961184) -- (-3.85032,-2.2972) -- cycle;
\fill[line width=0.5pt,fill=black,fill opacity=0.1] (-0.541995,-1.980592) -- (1.92516,-1.5141485) -- (3.85032,-3.028297) -- (-1.08399,-3.961184) -- cycle;
\fill[line width=0.5pt,fill=black,fill opacity=0.1] (0.5419949999999997,-0.6821565000000002) -- (1.92516,-1.5141485) -- (3.85032,-3.028297) -- (1.08399,-1.3643130000000003) -- cycle;
\fill[line width=0.5pt,fill=black,fill opacity=0.1] (-3.85032,-2.2972) -- (1.08399,-1.3643130000000003) -- (0.5419949999999997,-0.6821565000000002) -- (-1.92516,-1.1486) -- cycle;
\fill[line width=0.5pt,fill=black,fill opacity=0.1] (-3.85032,-2.2972) -- (-1.08399,-3.961184) -- (3.85032,-3.028297) -- (1.08399,-1.3643130000000003) -- (-3.85032,-2.2972) -- cycle;

\fill[line width=0.5pt,fill=black,fill opacity=0] (-0.541995,0.6821570000000001) -- (1.92516,1.1486005000000001) -- (1.92516,-1.5141485) -- (-0.541995,-1.980592) -- cycle;
\fill[line width=0.5pt,fill=black,fill opacity=0] (-1.92516,-1.1486) -- (0.5419949999999997,-0.6821565000000002) -- (0.5419950000000003,1.9805925) -- (-1.92516,1.514149) -- cycle;
\fill[line width=0.5pt,fill=black,fill opacity=0] (-1.92516,-1.1486) -- (-1.92516,1.514149) -- (-3.85032,3.028298) -- (-3.85032,-2.2972) -- cycle;
\fill[line width=0.5pt,fill=black,fill opacity=0] (-1.08399,1.3643140000000002) -- (-0.541995,0.6821570000000001) -- (-0.541995,-1.980592) -- (-1.08399,-3.961184) -- cycle;
\fill[line width=0.5pt,fill=black,fill opacity=0] (-1.92516,-1.1486) -- (-0.541995,-1.980592) -- (-0.541995,0.6821570000000001) -- (-1.92516,1.514149) -- cycle;
\fill[line width=0.5pt,fill=black,fill opacity=0] (1.92516,1.1486005000000001) -- (3.85032,2.2972010000000003) -- (3.85032,-3.028297) -- (1.92516,-1.5141485) -- cycle;
\fill[line width=0.5pt,fill=black,fill opacity=0] (0.5419950000000003,1.9805925) -- (1.08399,3.961185) -- (1.08399,-1.3643130000000003) -- (0.5419949999999997,-0.6821565000000002) -- cycle;
\fill[line width=0.5pt,fill=black,fill opacity=0] (0.5419950000000003,1.9805925) -- (1.92516,1.1486005000000001) -- (1.92516,-1.5141485) -- (0.5419949999999997,-0.6821565000000002) -- cycle;
\fill[line width=0.5pt,fill=black,fill opacity=0] (-0.541995,0.6821570000000001) -- (1.92516,1.1486005000000001) -- (1.92516,-1.5141485) -- (-0.541995,-1.980592) -- cycle;

\draw [line width=0.5pt] (3.85032,2.2972010000000003)-- (-1.08399,1.3643140000000002);
\draw [line width=0.5pt] (-1.08399,1.3643140000000002)-- (-1.08399,-3.961184);
\draw [line width=0.5pt] (-1.92516,1.514149)-- (-1.08399,1.008174144765809);
\draw [line width=0.5pt,dotted] (-1.08399,1.008174144765809)-- (-0.541995,0.6821570000000001);
\draw [line width=0.5pt] (0.5419950000000003,1.9805925)-- (0.9326826447970984,1.7455887260684537);
\draw [line width=0.5pt,dotted] (0.9326826447970984,1.7455887260684537)-- (1.92516,1.1486005000000001);
\draw [line width=0.5pt] (-1.92516,-1.1486)-- (-1.08399,-1.6545748552341912);
\draw [line width=0.5pt,dotted] (-1.08399,-1.6545748552341912)-- (-0.541995,-1.980592);
\draw [line width=0.5pt,dotted] (-1.92516,1.514149)-- (-1.4746776447970984,1.5993177739315465);
\draw [line width=0.5pt] (-1.4746776447970984,1.5993177739315465)-- (0.5419950000000003,1.9805925);
\draw [line width=0.5pt] (3.85032,2.2972010000000003)-- (3.85032,-3.028297);
\draw [line width=0.5pt] (3.85032,-3.028297)-- (-1.08399,-3.961184);
\draw [line width=0.5pt] (-0.541995,-1.980592)-- (-1.08399,-3.961184);
\draw [line width=0.5pt] (-1.08399,-3.961184)-- (-3.85032,-2.2972);
\draw [line width=0.5pt] (-3.85032,3.028298)-- (-3.85032,-2.2972);
\draw [line width=0.5pt] (-3.85032,3.028298)-- (-1.08399,1.3643140000000002);
\draw [line width=0.5pt] (-1.08399,1.3643140000000002)-- (-0.541995,0.6821570000000001);
\draw [line width=0.5pt] (-3.85032,3.028298)-- (-1.92516,1.514149);
\draw [line width=0.5pt] (1.08399,3.961185)-- (-3.85032,3.028298);
\draw [line width=0.5pt] (1.08399,3.961185)-- (3.85032,2.2972010000000003);
\draw [line width=0.5pt] (3.85032,2.2972010000000003)-- (1.92516,1.1486005000000001);
\draw [line width=0.5pt, dotted] (-3.85032,-2.2972)-- (1.08399,-1.3643130000000003);
\draw [line width=0.5pt] (-1.92516,-1.1486)-- (-3.85032,-2.2972);
\draw [line width=0.5pt] (1.92516,-1.5141485)-- (3.85032,-3.028297);
\draw [line width=0.5pt, dotted] (3.85032,-3.028297)-- (1.08399,-1.3643130000000003);
\draw [line width=0.5pt, dotted] (-1.92516,-1.1486)-- (0.5419949999999997,-0.6821565000000002);
\draw [line width=0.5pt] (-1.92516,1.514149)-- (-1.92516,-1.1486);
\draw [line width=0.5pt] (1.92516,1.1486005000000001)-- (1.92516,-1.5141485);
\draw [line width=0.5pt, dotted] (1.92516,-1.5141485)-- (0.5419949999999997,-0.6821565000000002);
\draw [line width=0.5pt, dotted] (0.5419949999999997,-0.6821565000000002)-- (0.5419950000000003,1.9805925);
\draw [line width=0.5pt] (0.5419950000000003,1.9805925)-- (1.08399,3.961185);
\draw [line width=0.5pt, dotted] (1.08399,3.961185)-- (1.08399,-1.3643130000000003);
\draw [line width=0.5pt, dotted] (1.08399,-1.3643130000000003)-- (0.5419949999999997,-0.6821565000000002);
\draw [line width=0.5pt, dotted] (0.5419949999999997,-0.6821565000000002)-- (0.5419950000000003,1.9805925);
\draw [line width=0.5pt] (-0.541995,0.6821570000000001)-- (1.92516,1.1486005000000001);
\draw [line width=0.5pt] (1.92516,1.1486005000000001)-- (1.92516,-1.5141485);
\draw [line width=0.5pt] (1.92516,-1.5141485)-- (-0.541995,-1.980592);
\draw [line width=0.5pt] (-0.541995,-1.980592)-- (-0.541995,0.6821570000000001);
\begin{tiny}
\draw [fill=black] (-1.92516,-1.1486) circle (2.5pt);
\draw [fill=black] (-0.541995,-1.980592) circle (2.5pt);
\draw [fill=black] (1.92516,-1.5141485) circle (2.5pt);
\draw [fill=black] (-1.92516,1.514149) circle (2.0pt);
\draw [fill=black] (-0.541995,0.6821570000000001) circle (2.0pt);
\draw [fill=black] (1.92516,1.1486005000000001) circle (2.0pt);
\draw [fill=black] (0.5419950000000003,1.9805925) circle (2.0pt);
\draw [fill=black] (0.5419949999999997,-0.6821565000000002) circle (2.0pt);

\draw [fill=black] (-3.85032,-2.2972) circle (2.5pt);
\draw [fill=black] (-1.08399,-3.961184) circle (2.5pt);
\draw [fill=black] (3.85032,-3.028297) circle (2.5pt);
\draw [fill=black] (-3.85032,3.028298) circle (2.0pt);
\draw [fill=black] (-1.08399,1.3643140000000002) circle (2.0pt);
\draw [fill=black] (3.85032,2.2972010000000003) circle (2.0pt);
\draw [fill=black] (1.08399,3.961185) circle (2.0pt);
\draw [fill=black] (1.08399,-1.3643130000000003) circle (2.0pt);
\draw[color=black] (-1,2.75) node {$t_3^{(0)}$};
\draw[color=black] (1,-2.75) node {$t_3^{(1)}$};
\end{tiny}
\end{tikzpicture}
\end{minipage}
\hspace{-60pt}
\begin{minipage}[c]{.3\textwidth}
\centering
\begin{tikzpicture}[line cap=round,line join=round,>=triangle 45,x=1.0cm,y=1.0cm, scale=0.3]
\clip(-4.5,-4.1) rectangle (4.5,4.35);

\fill[line width=0.5pt,fill=black,fill opacity=0.1] (-1.08399,1.3643140000000002) -- (-0.541995,0.6821570000000001) -- (1.92516,1.1486005000000001) -- (3.85032,2.2972010000000003) -- cycle;
\fill[line width=0.5pt,fill=black,fill opacity=0.1] (-1.08399,1.3643140000000002) -- (-0.541995,0.6821570000000001) -- (-0.541995,-1.980592) -- (-1.08399,-3.961184) -- cycle;
\fill[line width=0.5pt,fill=black,fill opacity=0.1] (-0.541995,0.6821570000000001) -- (1.92516,1.1486005000000001) -- (1.92516,-1.5141485) -- (-0.541995,-1.980592) -- cycle;
\fill[line width=0.5pt,fill=black,fill opacity=0.1] (1.92516,1.1486005000000001) -- (3.85032,2.2972010000000003) -- (3.85032,-3.028297) -- (1.92516,-1.5141485) -- cycle;
\fill[line width=0.5pt,fill=black,fill opacity=0.1] (-0.541995,-1.980592) -- (1.92516,-1.5141485) -- (3.85032,-3.028297) -- (-1.08399,-3.961184) -- cycle;

\fill[line width=0.5pt,fill=black,fill opacity=0.1] (-3.85032,3.028298) -- (-1.92516,1.514149) -- (0.5419950000000003,1.9805925) -- (1.08399,3.961185) -- cycle;
\fill[line width=0.5pt,fill=black,fill opacity=0.1] (-3.85032,-2.2972) -- (1.08399,-1.3643130000000003) -- (0.5419949999999997,-0.6821565000000002) -- (-1.92516,-1.1486) -- cycle;
\fill[line width=0.5pt,fill=black,fill opacity=0.1] (-1.92516,-1.1486) -- (0.5419949999999997,-0.6821565000000002) -- (0.5419950000000003,1.9805925) -- (-1.92516,1.514149) -- cycle;
\fill[line width=0.5pt,fill=black,fill opacity=0.1] (0.5419950000000003,1.9805925) -- (1.08399,3.961185) -- (1.08399,-1.3643130000000003) -- (0.5419949999999997,-0.6821565000000002) -- cycle;
\fill[line width=0.5pt,fill=black,fill opacity=0.1] (-1.92516,-1.1486) -- (-1.92516,1.514149) -- (-3.85032,3.028298) -- (-3.85032,-2.2972) -- cycle;

\fill[line width=0.5pt,fill=black,fill opacity=0.1] (-1.08399,1.3643140000000002) -- (3.85032,2.2972010000000003) -- (3.85032,-3.028297) -- (-1.08399,-3.961184) -- (-1.08399,1.3643140000000002) -- cycle;
\fill[line width=0.5pt,fill=black,fill opacity=0.1] (-3.85032,3.028298) -- (1.08399,3.961185) -- (1.08399,-1.3643130000000003) -- (-3.85032,-2.2972) -- (-3.85032,3.028298) -- cycle;

\fill[line width=0.5pt,fill=black,fill opacity=0] (-1.92516,-1.1486) -- (-0.541995,-1.980592) -- (-1.08399,-3.961184) -- (-3.85032,-2.2972) -- cycle;
\fill[line width=0.5pt,fill=black,fill opacity=0] (-1.92516,-1.1486) -- (-0.541995,-1.980592) -- (-0.541995,0.6821570000000001) -- (-1.92516,1.514149) -- cycle;
\fill[line width=0.5pt,fill=black,fill opacity=0] (-3.85032,3.028298) -- (-1.08399,1.3643140000000002) -- (-0.541995,0.6821570000000001) -- (-1.92516,1.514149) -- cycle;

\fill[line width=0.5pt,fill=black,fill opacity=0] (-1.92516,1.514149) -- (0.5419950000000003,1.9805925) -- (1.92516,1.1486005000000001) -- (-0.541995,0.6821570000000001) -- cycle;

\fill[line width=0.5pt,fill=black,fill opacity=0] (0.5419950000000003,1.9805925) -- (1.08399,3.961185) -- (3.85032,2.2972010000000003) -- (1.92516,1.1486005000000001) -- cycle;
\fill[line width=0.5pt,fill=black,fill opacity=0] (0.5419949999999997,-0.6821565000000002) -- (1.92516,-1.5141485) -- (3.85032,-3.028297) -- (1.08399,-1.3643130000000003) -- cycle;
\fill[line width=0.5pt,fill=black,fill opacity=0] (-0.541995,-1.980592) -- (1.92516,-1.5141485) -- (0.5419949999999997,-0.6821565000000002) -- (-1.92516,-1.1486) -- cycle;
\fill[line width=0.5pt,fill=black,fill opacity=0] (0.5419950000000003,1.9805925) -- (1.92516,1.1486005000000001) -- (1.92516,-1.5141485) -- (0.5419949999999997,-0.6821565000000002) -- cycle;
\fill[line width=0.5pt,fill=black,fill opacity=0] (-0.541995,0.6821570000000001) -- (1.92516,1.1486005000000001) -- (1.92516,-1.5141485) -- (-0.541995,-1.980592) -- cycle;

\draw [line width=0.5pt] (3.85032,2.2972010000000003)-- (-1.08399,1.3643140000000002);
\draw [line width=0.5pt] (-1.08399,1.3643140000000002)-- (-1.08399,-3.961184);
\draw [line width=0.5pt] (-1.92516,1.514149)-- (-1.08399,1.008174144765809);
\draw [line width=0.5pt,dotted] (-1.08399,1.008174144765809)-- (-0.541995,0.6821570000000001);
\draw [line width=0.5pt] (0.5419950000000003,1.9805925)-- (0.9326826447970984,1.7455887260684537);
\draw [line width=0.5pt,dotted] (0.9326826447970984,1.7455887260684537)-- (1.92516,1.1486005000000001);
\draw [line width=0.5pt] (-1.92516,-1.1486)-- (-1.08399,-1.6545748552341912);
\draw [line width=0.5pt,dotted] (-1.08399,-1.6545748552341912)-- (-0.541995,-1.980592);
\draw [line width=0.5pt,dotted] (-1.92516,1.514149)-- (-1.4746776447970984,1.5993177739315465);
\draw [line width=0.5pt] (-1.4746776447970984,1.5993177739315465)-- (0.5419950000000003,1.9805925);
\draw [line width=0.5pt] (3.85032,2.2972010000000003)-- (3.85032,-3.028297);
\draw [line width=0.5pt] (3.85032,-3.028297)-- (-1.08399,-3.961184);
\draw [line width=0.5pt] (-0.541995,-1.980592)-- (-1.08399,-3.961184);
\draw [line width=0.5pt] (-1.08399,-3.961184)-- (-3.85032,-2.2972);
\draw [line width=0.5pt] (-3.85032,3.028298)-- (-3.85032,-2.2972);
\draw [line width=0.5pt] (-3.85032,3.028298)-- (-1.08399,1.3643140000000002);
\draw [line width=0.5pt] (-1.08399,1.3643140000000002)-- (-0.541995,0.6821570000000001);
\draw [line width=0.5pt] (-3.85032,3.028298)-- (-1.92516,1.514149);
\draw [line width=0.5pt] (1.08399,3.961185)-- (-3.85032,3.028298);
\draw [line width=0.5pt] (1.08399,3.961185)-- (3.85032,2.2972010000000003);
\draw [line width=0.5pt] (3.85032,2.2972010000000003)-- (1.92516,1.1486005000000001);
\draw [line width=0.5pt, dotted] (-3.85032,-2.2972)-- (1.08399,-1.3643130000000003);
\draw [line width=0.5pt] (-1.92516,-1.1486)-- (-3.85032,-2.2972);
\draw [line width=0.5pt] (1.92516,-1.5141485)-- (3.85032,-3.028297);
\draw [line width=0.5pt, dotted] (3.85032,-3.028297)-- (1.08399,-1.3643130000000003);
\draw [line width=0.5pt, dotted] (-1.92516,-1.1486)-- (0.5419949999999997,-0.6821565000000002);
\draw [line width=0.5pt] (-1.92516,1.514149)-- (-1.92516,-1.1486);
\draw [line width=0.5pt] (1.92516,1.1486005000000001)-- (1.92516,-1.5141485);
\draw [line width=0.5pt, dotted] (1.92516,-1.5141485)-- (0.5419949999999997,-0.6821565000000002);
\draw [line width=0.5pt, dotted] (0.5419949999999997,-0.6821565000000002)-- (0.5419950000000003,1.9805925);
\draw [line width=0.5pt] (0.5419950000000003,1.9805925)-- (1.08399,3.961185);
\draw [line width=0.5pt, dotted] (1.08399,3.961185)-- (1.08399,-1.3643130000000003);
\draw [line width=0.5pt, dotted] (1.08399,-1.3643130000000003)-- (0.5419949999999997,-0.6821565000000002);
\draw [line width=0.5pt, dotted] (0.5419949999999997,-0.6821565000000002)-- (0.5419950000000003,1.9805925);
\draw [line width=0.5pt] (-0.541995,0.6821570000000001)-- (1.92516,1.1486005000000001);
\draw [line width=0.5pt] (1.92516,1.1486005000000001)-- (1.92516,-1.5141485);
\draw [line width=0.5pt] (1.92516,-1.5141485)-- (-0.541995,-1.980592);
\draw [line width=0.5pt] (-0.541995,-1.980592)-- (-0.541995,0.6821570000000001);
\begin{tiny}
\draw [fill=black] (-1.92516,-1.1486) circle (2.5pt);
\draw [fill=black] (-0.541995,-1.980592) circle (2.5pt);
\draw [fill=black] (1.92516,-1.5141485) circle (2.5pt);
\draw [fill=black] (-1.92516,1.514149) circle (2.0pt);
\draw [fill=black] (-0.541995,0.6821570000000001) circle (2.0pt);
\draw [fill=black] (1.92516,1.1486005000000001) circle (2.0pt);
\draw [fill=black] (0.5419950000000003,1.9805925) circle (2.0pt);
\draw [fill=black] (0.5419949999999997,-0.6821565000000002) circle (2.0pt);
\draw [fill=black] (-3.85032,-2.2972) circle (2.5pt);
\draw [fill=black] (-1.08399,-3.961184) circle (2.5pt);
\draw [fill=black] (3.85032,-3.028297) circle (2.5pt);
\draw [fill=black] (-3.85032,3.028298) circle (2.0pt);
\draw [fill=black] (-1.08399,1.3643140000000002) circle (2.0pt);
\draw [fill=black] (3.85032,2.2972010000000003) circle (2.0pt);
\draw [fill=black] (1.08399,3.961185) circle (2.0pt);
\draw [fill=black] (1.08399,-1.3643130000000003) circle (2.0pt);
\draw[color=black] (-2.9,-1) node {$t_4^{(0)}$};
\draw[color=black] (0.1,-2.8) node {$t_4^{(1)}$};
\end{tiny}
\end{tikzpicture}
\end{minipage}
\end{center}
\caption{The slices $t_j^{(0)}$ and $t_j^{(1)}$ appearing in the computation of $f_{4,\{j\}}$.}\label{fig: case d=4 factors constant term}
\end{figure}
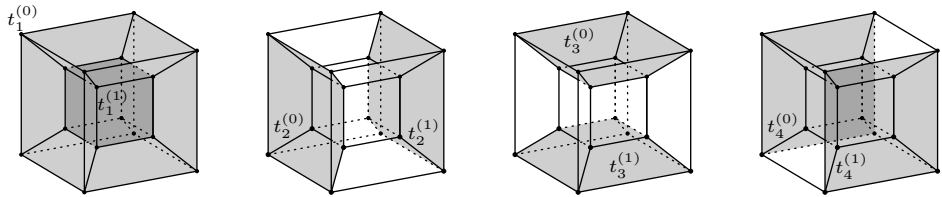

As pointed out in the proof of Proposition \ref{SOS}, formula (\ref{eq: product hyperdeterminants}) may be generalized to any polynomial $f_{d,J}\neq 1$. For example, below we interpret the equation $f_{4,\{1,2\}}$ of $X_{4,\{1,2\}}^\vee$ in terms of the tensors $t_{\{1,2\}}^{(rs)}\in V\otimes V$, with $r,s\in\{0,1\}$, obtained extracting from $t$ the slices $\{t_{rsij}\}$ highlighted in Figure \ref{fig: QQPP}.
{\scriptsize
\begin{align*}
f_{4,\{1,2\}}&=\det\left[\left(\left(t_{\{1,2\}}^{(00)}+\sqrt{-1}t_{\{1,2\}}^{(10)}\right)+\sqrt{-1}\left(t_{\{1,2\}}^{(01)}+\sqrt{-1}t_{\{1,2\}}^{(11)}\right)\right)\left(\left(t_{\{1,2\}}^{(00)}+\sqrt{-1}t_{\{1,2\}}^{(10)}\right)-\sqrt{-1}\left(t_{\{1,2\}}^{(01)}+\sqrt{-1}t_{\{1,2\}}^{(11)}\right)\right)\right]\\
&\hspace{10pt}\cdot\det\left[\left(\left(t_{\{1,2\}}^{(00)}-\sqrt{-1}t_{\{1,2\}}^{(10)}\right)+\sqrt{-1}\left(t_{\{1,2\}}^{(01)}-\sqrt{-1}t_{\{1,2\}}^{(11)}\right)\right)\left(\left(t_{\{1,2\}}^{(00)}-\sqrt{-1}t_{\{1,2\}}^{(10)}\right)-\sqrt{-1}\left(t_{\{1,2\}}^{(01)}-\sqrt{-1}t_{\{1,2\}}^{(11)}\right)\right)\right].
\end{align*}
}
\begin{figure}[ht]
\begin{tikzpicture}[line cap=round,line join=round,>=triangle 45,x=1.0cm,y=1.0cm, scale=0.7]
\clip(-4.9,-4.3) rectangle (4.9,4.3);
\fill[line width=0.5pt,fill=black,fill opacity=0.10000000149011612] (-1.08399,1.3643140000000002) -- (-0.541995,0.6821570000000001) -- (1.92516,1.1486005000000001) -- (3.85032,2.2972010000000003) -- cycle;
\fill[line width=0.5pt,fill=black,fill opacity=0.15000000596046448] (-0.541995,-1.980592) -- (1.92516,-1.5141485) -- (3.85032,-3.028297) -- (-1.08399,-3.961184) -- cycle;
\fill[line width=0.5pt,fill=black,fill opacity=0.10000000149011612] (-3.85032,3.028298) -- (-1.92516,1.514149) -- (0.5419950000000003,1.9805925) -- (1.08399,3.961185) -- cycle;
\fill[line width=0.5pt,fill=black,fill opacity=0.10000000149011612] (-3.85032,-2.2972) -- (1.08399,-1.3643130000000003) -- (0.5419949999999997,-0.6821565000000002) -- (-1.92516,-1.1486) -- cycle;

\draw [line width=0.5pt] (3.85032,2.2972010000000003)-- (-1.08399,1.3643140000000002);
\draw [line width=0.5pt] (-1.08399,1.3643140000000002)-- (-1.08399,-3.961184);
\draw [line width=0.5pt] (-1.92516,1.514149)-- (-1.08399,1.008174144765809);
\draw [line width=0.5pt,dotted] (-1.08399,1.008174144765809)-- (-0.541995,0.6821570000000001);
\draw [line width=0.5pt] (0.5419950000000003,1.9805925)-- (0.9326826447970984,1.7455887260684537);
\draw [line width=0.5pt,dotted] (0.9326826447970984,1.7455887260684537)-- (1.92516,1.1486005000000001);
\draw [line width=0.5pt] (-1.92516,-1.1486)-- (-1.08399,-1.6545748552341912);
\draw [line width=0.5pt,dotted] (-1.08399,-1.6545748552341912)-- (-0.541995,-1.980592);
\draw [line width=0.5pt,dotted] (-1.92516,1.514149)-- (-1.4746776447970984,1.5993177739315465);
\draw [line width=0.5pt] (-1.4746776447970984,1.5993177739315465)-- (0.5419950000000003,1.9805925);
\draw [line width=0.5pt] (3.85032,2.2972010000000003)-- (3.85032,-3.028297);
\draw [line width=0.5pt] (3.85032,-3.028297)-- (-1.08399,-3.961184);
\draw [line width=0.5pt] (-0.541995,-1.980592)-- (-1.08399,-3.961184);
\draw [line width=0.5pt] (-1.08399,-3.961184)-- (-3.85032,-2.2972);
\draw [line width=0.5pt] (-3.85032,3.028298)-- (-3.85032,-2.2972);
\draw [line width=0.5pt] (-3.85032,3.028298)-- (-1.08399,1.3643140000000002);
\draw [line width=0.5pt] (-1.08399,1.3643140000000002)-- (-0.541995,0.6821570000000001);
\draw [line width=0.5pt] (-3.85032,3.028298)-- (-1.92516,1.514149);
\draw [line width=0.5pt] (1.08399,3.961185)-- (-3.85032,3.028298);
\draw [line width=0.5pt] (1.08399,3.961185)-- (3.85032,2.2972010000000003);
\draw [line width=0.5pt] (3.85032,2.2972010000000003)-- (1.92516,1.1486005000000001);
\draw [line width=0.5pt, dotted] (-3.85032,-2.2972)-- (1.08399,-1.3643130000000003);
\draw [line width=0.5pt] (-1.92516,-1.1486)-- (-3.85032,-2.2972);
\draw [line width=0.5pt] (1.92516,-1.5141485)-- (3.85032,-3.028297);
\draw [line width=0.5pt, dotted] (3.85032,-3.028297)-- (1.08399,-1.3643130000000003);
\draw [line width=0.5pt, dotted] (-1.92516,-1.1486)-- (0.5419949999999997,-0.6821565000000002);
\draw [line width=0.5pt] (-1.92516,1.514149)-- (-1.92516,-1.1486);
\draw [line width=0.5pt] (1.92516,1.1486005000000001)-- (1.92516,-1.5141485);
\draw [line width=0.5pt, dotted] (1.92516,-1.5141485)-- (0.5419949999999997,-0.6821565000000002);
\draw [line width=0.5pt, dotted] (0.5419949999999997,-0.6821565000000002)-- (0.5419950000000003,1.9805925);
\draw [line width=0.5pt] (0.5419950000000003,1.9805925)-- (1.08399,3.961185);
\draw [line width=0.5pt, dotted] (1.08399,3.961185)-- (1.08399,-1.3643130000000003);
\draw [line width=0.5pt, dotted] (1.08399,-1.3643130000000003)-- (0.5419949999999997,-0.6821565000000002);
\draw [line width=0.5pt, dotted] (0.5419949999999997,-0.6821565000000002)-- (0.5419950000000003,1.9805925);
\draw [line width=0.5pt] (-0.541995,0.6821570000000001)-- (1.92516,1.1486005000000001);
\draw [line width=0.5pt] (1.92516,1.1486005000000001)-- (1.92516,-1.5141485);
\draw [line width=0.5pt] (1.92516,-1.5141485)-- (-0.541995,-1.980592);
\draw [line width=0.5pt] (-0.541995,-1.980592)-- (-0.541995,0.6821570000000001);
\begin{tiny}
\draw [fill=black] (-1.92516,-1.1486) circle (1pt);
\draw[color=black] (-2.45,-0.9) node {$1001$};
\draw [fill=black] (-1.92516,1.514149) circle (1pt);
\draw[color=black] (-2.45,1.3) node {$0001$};
\draw [fill=black] (-0.541995,-1.980592) circle (1pt);
\draw[color=black] (0,-2.2) node {$1101$};
\draw [fill=black] (1.92516,-1.5141485) circle (1pt);
\draw[color=black] (2.45,-1.3) node {$1111$};
\draw [fill=black] (-0.541995,0.682157) circle (1pt);
\draw[color=black] (0,0.5) node {$0101$};
\draw [fill=black] (1.92516,1.1486005) circle (1pt);
\draw[color=black] (2.45,0.9) node {$0111$};
\draw [fill=black] (0.541995,1.9805925) circle (1pt);
\draw[color=black] (0,2.2) node {$0011$};
\draw [fill=black] (0.5419949999999,-0.6821565) circle (1pt);
\draw[color=black] (0,-0.5) node {$1011$};
\draw [fill=black] (-3.85032,3.028298) circle (1pt);
\draw[color=black] (-4.4,2.8) node {$0000$};
\draw [fill=black] (3.85032,2.297201) circle (1pt);
\draw[color=black] (4.4,2.1) node {$0110$};
\draw [fill=black] (-1.08399,1.364314) circle (1pt);
\draw[color=black] (-0.45,1.241) node {$0100$};
\draw [fill=black] (1.08399,3.961185) circle (1pt);
\draw[color=black] (1.6,4.1) node {$0010$};
\draw [fill=black] (-1.08399,-3.961184) circle (1pt);
\draw[color=black] (-1.6,-4.1) node {$1100$};
\draw [fill=black] (-3.85032,-2.2972) circle (1pt);
\draw[color=black] (-4.4,-2.1) node {$1000$};
\draw [fill=black] (3.85032,-3.028297) circle (1pt);
\draw[color=black] (4.4,-2.8) node {$1110$};
\draw [fill=black] (1.08399,-1.3643130000000003) circle (1pt);
\draw[color=black] (0.45,-1.241) node {$1010$};
\end{tiny}
\begin{scriptsize}
\draw[color=black] (-1,2.8) node {$t_{\{1,2\}}^{(00)}$};
\draw[color=black] (1,-2.8) node {$t_{\{1,2\}}^{(11)}$};
\draw[color=black] (-1.9,-1.8) node {$t_{\{1,2\}}^{(10)}$};
\draw[color=black] (1.9,1.6) node {$t_{\{1,2\}}^{(01)}$};
\end{scriptsize}
\end{tikzpicture}
\caption{The slices $f_{\{1,2\}}^{(rs)}$ appearing in the expression of $f_{4,\{1,2\}}(t)$.}\label{fig: QQPP}
\end{figure}
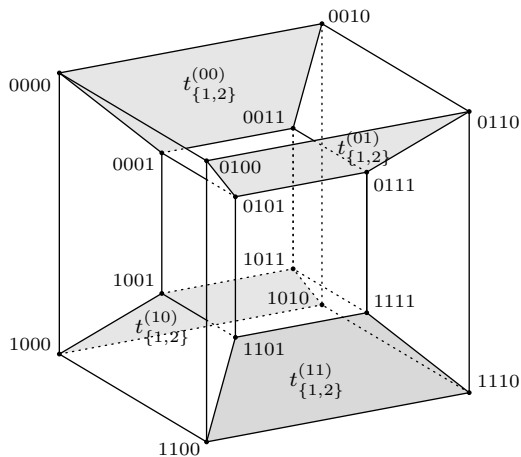
\end{remark}

\section{Explicit computations in the case $d=3$}\label{sec: 222}

We recall that $V$ is the complexification of the real vector space $V_\R$. In particular, $q$ is not a Hermitian form. For example, it is symmetric rather than conjugate-symmetric. The subgroup of all linear operators $A\in\mathrm{GL}(V)$ which preserve this form (i.e., such that $q(Ax,Ay)=q(x,y)$ for all $x,y\in V$) is the {\em complex orthogonal group} $\mathrm{O}(V)$. The group $\mathrm{SO}(V)$ is defined to be the set of all $A$ in $\mathrm{O}(V)$ with $\det(A)=1$ and is called {\em complex special orthogonal group} over $V$.

The action of $\mathrm{SO}(V)$ on $V$ induces another action of $\mathrm{SO}(V)^d$ on the tensor space $V^{\otimes d}$. Thanks to \cite[Proposition 2.11]{OS}, the coefficients of $\mathrm{EDpoly}_{X^\vee,t}(\epsilon^2)$ are $\mathrm{SO}(V)^d$-invariants. Indeed we are interested in computing a minimal generating set for the invariant ring $S(V^{\otimes d})^{\mathrm{SO}(V)^d}$.

As in the previous sections, $q$ is the standard Euclidean scalar product. We fix $x_{j,0},x_{j,1}$ as coordinates for the $j$-th copy of $V$ in $V^{\otimes d}$. Then, the associated quadratic form $q$ is in coordinates $x_{j,0}^2+x_{j,1}^2$ for all $1\le j\le d$. Now consider the change of coordinates $z_{j,0}=x_{j,0}+\sqrt{-1}x_{j,1}$, $z_{j,1}=x_{j,0}-\sqrt{-1}x_{j,1}$. In these new coordinates, the expression for the quadratic form $q$ on the $j$-th copy of $V$ in $V^{\otimes d}$ becomes $z_{j,0}z_{j,1}$. Moreover, each binary tensor $t=(t_{i_1\cdots i_d})\in V^{\otimes d}$ may be written as
\[
t=\sum_{(i_1,\ldots,i_d)\in \{0,1\}^d}t_{i_1\cdots i_d}x_{1,i_1}\cdots x_{d,i_d}=\sum_{(i_1,\ldots,i_d)\in \{0,1\}^d}u_{i_1\cdots i_d}z_{1,i_1}\cdots z_{d,i_d},
\]
for some coefficients $u_{i_1\cdots i_d}$ depending on the old set of coordinates $\{t_{i_1\cdots i_d}\}$ via the following relations:
\[
u_{i_1\cdots i_d}=\sum_{(j_1,\ldots,j_d)\in \{0,1\}^d}\left[\sqrt{-1}^{\sum_{l=1}^dj_l}(-1)^{\sum_{l=1}^d(1-i_l)j_l}\right]t_{j_1\cdots j_d}.
\]
One may verify by direct computation that, for all $(i_1,\ldots,i_d)\in\{0,1\}^d$, the complex conjugate of $u_{i_1\cdots i_d}$ is $u_{k_1\cdots k_d}$, where $k_l=1-i_l$.

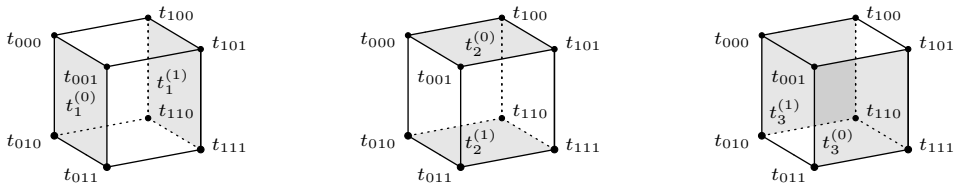
\begin{figure}[ht]
\begin{center}
    \begin{minipage}[c]{.3\textwidth}
    \centering
    \begin{tikzpicture}[line cap=round,line join=round,>=triangle 45,x=1.0cm,y=1.0cm, scale=0.5]
\clip(-3.3,-3) rectangle (3.3,3);

\fill[line width=0.5pt,fill opacity=0] (-0.541995,0.6821570000000001) -- (1.92516,1.1486005000000001) -- (1.92516,-1.5141485) -- (-0.541995,-1.980592) -- cycle;
\fill[line width=0.5pt,fill=black,fill opacity=0] (-1.92516,-1.1486) -- (0.5419949999999997,-0.6821565000000002) -- (0.5419950000000003,1.9805925) -- (-1.92516,1.514149) -- cycle;

\fill[line width=0.5pt,fill=black,fill opacity=0.1] (-1.92516,-1.1486) -- (-0.541995,-1.980592) -- (-0.541995,0.6821570000000001) -- (-1.92516,1.514149) -- cycle;
\fill[line width=0.5pt,fill=black,fill opacity=0.1] (0.5419950000000003,1.9805925) -- (1.92516,1.1486005000000001) -- (1.92516,-1.5141485) -- (0.5419949999999997,-0.6821565000000002) -- cycle;

\fill[line width=0.5pt,fill=black,fill opacity=0] (-1.92516,1.514149) -- (0.5419950000000003,1.9805925) -- (1.92516,1.1486005000000001) -- (-0.541995,0.6821570000000001) -- cycle;
\fill[line width=0.5pt,fill=black,fill opacity=0] (-0.541995,-1.980592) -- (1.92516,-1.5141485) -- (0.5419949999999997,-0.6821565000000002) -- (-1.92516,-1.1486) -- cycle;

\draw [line width=0.5pt] (-1.92516,1.514149)-- (-1.08399,1.008174144765809);
\draw [line width=0.5pt] (-1.08399,1.008174144765809)-- (-0.541995,0.6821570000000001);
\draw [line width=0.5pt] (0.5419950000000003,1.9805925)-- (0.9326826447970984,1.7455887260684537);
\draw [line width=0.5pt] (0.9326826447970984,1.7455887260684537)-- (1.92516,1.1486005000000001);
\draw [line width=0.5pt] (-1.92516,-1.1486)-- (-1.08399,-1.6545748552341912);
\draw [line width=0.5pt] (-1.08399,-1.6545748552341912)-- (-0.541995,-1.980592);
\draw [line width=0.5pt] (-1.92516,1.514149)-- (-1.4746776447970984,1.5993177739315465);
\draw [line width=0.5pt] (-1.4746776447970984,1.5993177739315465)-- (0.5419950000000003,1.9805925);
\draw [line width=0.5pt, dotted] (-1.92516,-1.1486)-- (0.5419949999999997,-0.6821565000000002);
\draw [line width=0.5pt] (-1.92516,1.514149)-- (-1.92516,-1.1486);
\draw [line width=0.5pt] (1.92516,1.1486005000000001)-- (1.92516,-1.5141485);
\draw [line width=0.5pt, dotted] (1.92516,-1.5141485)-- (0.5419949999999997,-0.6821565000000002);
\draw [line width=0.5pt, dotted] (0.5419949999999997,-0.6821565000000002)-- (0.5419950000000003,1.9805925);
\draw [line width=0.5pt, dotted] (0.5419949999999997,-0.6821565000000002)-- (0.5419950000000003,1.9805925);
\draw [line width=0.5pt] (-0.541995,0.6821570000000001)-- (1.92516,1.1486005000000001);
\draw [line width=0.5pt] (1.92516,1.1486005000000001)-- (1.92516,-1.5141485);
\draw [line width=0.5pt] (1.92516,-1.5141485)-- (-0.541995,-1.980592);
\draw [line width=0.5pt] (-0.541995,-1.980592)-- (-0.541995,0.6821570000000001);
\begin{tiny}
\draw [fill=black] (-1.92516,-1.1486) circle (2.5pt);
\draw[color=black] (-2.7,-1.3) node {$t_{010}$};
\draw [fill=black] (-0.541995,-1.980592) circle (2.5pt);
\draw[color=black] (-1.2,-2.2) node {$t_{011}$};
\draw [fill=black] (1.92516,-1.5141485) circle (2.5pt);
\draw[color=black] (2.7,-1.4) node {$t_{111}$};
\draw [fill=black] (-1.92516,1.514149) circle (2.0pt);
\draw[color=black] (-2.7,1.4) node {$t_{000}$};
\draw [fill=black] (-0.541995,0.6821570000000001) circle (2.0pt);
\draw[color=black] (-1.2,0.4) node {$t_{001}$};
\draw [fill=black] (1.92516,1.1486005000000001) circle (2.0pt);
\draw[color=black] (2.7,1.3) node {$t_{101}$};
\draw [fill=black] (0.5419950000000003,1.9805925) circle (2.0pt);
\draw[color=black] (1.3,2.1) node {$t_{100}$};
\draw [fill=black] (0.5419949999999997,-0.6821565000000002) circle (2.0pt);
\draw[color=black] (1.3,-0.5) node {$t_{110}$};
\draw[color=black] (1.2,0.3) node {$t_1^{(1)}$};
\draw[color=black] (-1.2,-0.3) node {$t_1^{(0)}$};
\end{tiny}
\end{tikzpicture}
    \end{minipage}
    \hspace{-20pt}
    \begin{minipage}[c]{.3\textwidth}
        \centering
        \begin{tikzpicture}[line cap=round,line join=round,>=triangle 45,x=1.0cm,y=1.0cm, scale=0.5]
\clip(-3.3,-3) rectangle (3.3,3);

\fill[line width=0.5pt,fill=black,fill opacity=0] (-0.541995,0.6821570000000001) -- (1.92516,1.1486005000000001) -- (1.92516,-1.5141485) -- (-0.541995,-1.980592) -- cycle;
\fill[line width=0.5pt,fill=black,fill opacity=0] (-1.92516,-1.1486) -- (0.5419949999999997,-0.6821565000000002) -- (0.5419950000000003,1.9805925) -- (-1.92516,1.514149) -- cycle;

\fill[line width=0.5pt,fill=black,fill opacity=0] (-1.92516,-1.1486) -- (-0.541995,-1.980592) -- (-0.541995,0.6821570000000001) -- (-1.92516,1.514149) -- cycle;
\fill[line width=0.5pt,fill=black,fill opacity=0] (0.5419950000000003,1.9805925) -- (1.92516,1.1486005000000001) -- (1.92516,-1.5141485) -- (0.5419949999999997,-0.6821565000000002) -- cycle;

\fill[line width=0.5pt,fill=black,fill opacity=0.1] (-1.92516,1.514149) -- (0.5419950000000003,1.9805925) -- (1.92516,1.1486005000000001) -- (-0.541995,0.6821570000000001) -- cycle;
\fill[line width=0.5pt,fill=black,fill opacity=0.1] (-0.541995,-1.980592) -- (1.92516,-1.5141485) -- (0.5419949999999997,-0.6821565000000002) -- (-1.92516,-1.1486) -- cycle;

\draw [line width=0.5pt] (-1.92516,1.514149)-- (-1.08399,1.008174144765809);
\draw [line width=0.5pt] (-1.08399,1.008174144765809)-- (-0.541995,0.6821570000000001);
\draw [line width=0.5pt] (0.5419950000000003,1.9805925)-- (0.9326826447970984,1.7455887260684537);
\draw [line width=0.5pt] (0.9326826447970984,1.7455887260684537)-- (1.92516,1.1486005000000001);
\draw [line width=0.5pt] (-1.92516,-1.1486)-- (-1.08399,-1.6545748552341912);
\draw [line width=0.5pt] (-1.08399,-1.6545748552341912)-- (-0.541995,-1.980592);
\draw [line width=0.5pt] (-1.92516,1.514149)-- (-1.4746776447970984,1.5993177739315465);
\draw [line width=0.5pt] (-1.4746776447970984,1.5993177739315465)-- (0.5419950000000003,1.9805925);
\draw [line width=0.5pt, dotted] (-1.92516,-1.1486)-- (0.5419949999999997,-0.6821565000000002);
\draw [line width=0.5pt] (-1.92516,1.514149)-- (-1.92516,-1.1486);
\draw [line width=0.5pt] (1.92516,1.1486005000000001)-- (1.92516,-1.5141485);
\draw [line width=0.5pt, dotted] (1.92516,-1.5141485)-- (0.5419949999999997,-0.6821565000000002);
\draw [line width=0.5pt, dotted] (0.5419949999999997,-0.6821565000000002)-- (0.5419950000000003,1.9805925);
\draw [line width=0.5pt, dotted] (0.5419949999999997,-0.6821565000000002)-- (0.5419950000000003,1.9805925);
\draw [line width=0.5pt] (-0.541995,0.6821570000000001)-- (1.92516,1.1486005000000001);
\draw [line width=0.5pt] (1.92516,1.1486005000000001)-- (1.92516,-1.5141485);
\draw [line width=0.5pt] (1.92516,-1.5141485)-- (-0.541995,-1.980592);
\draw [line width=0.5pt] (-0.541995,-1.980592)-- (-0.541995,0.6821570000000001);
\begin{tiny}
\draw [fill=black] (-1.92516,-1.1486) circle (2.5pt);
\draw[color=black] (-2.7,-1.3) node {$t_{010}$};
\draw [fill=black] (-0.541995,-1.980592) circle (2.5pt);
\draw[color=black] (-1.2,-2.2) node {$t_{011}$};
\draw [fill=black] (1.92516,-1.5141485) circle (2.5pt);
\draw[color=black] (2.7,-1.4) node {$t_{111}$};
\draw [fill=black] (-1.92516,1.514149) circle (2.0pt);
\draw[color=black] (-2.7,1.4) node {$t_{000}$};
\draw [fill=black] (-0.541995,0.6821570000000001) circle (2.0pt);
\draw[color=black] (-1.2,0.4) node {$t_{001}$};
\draw [fill=black] (1.92516,1.1486005000000001) circle (2.0pt);
\draw[color=black] (2.7,1.3) node {$t_{101}$};
\draw [fill=black] (0.5419950000000003,1.9805925) circle (2.0pt);
\draw[color=black] (1.3,2.1) node {$t_{100}$};
\draw [fill=black] (0.5419949999999997,-0.6821565000000002) circle (2.0pt);
\draw[color=black] (1.3,-0.5) node {$t_{110}$};

\draw[color=black] (0,1.3) node {$t_2^{(0)}$};
\draw[color=black] (0,-1.3) node {$t_2^{(1)}$};
\end{tiny}
\end{tikzpicture}
\end{minipage}
\hspace{-20pt}
\begin{minipage}[c]{.3\textwidth}
\centering
\begin{tikzpicture}[line cap=round,line join=round,>=triangle 45,x=1.0cm,y=1.0cm, scale=0.5]
\clip(-3.3,-3) rectangle (3.3,3);

\fill[line width=0.5pt,fill=black,fill opacity=0.1] (-0.541995,0.6821570000000001) -- (1.92516,1.1486005000000001) -- (1.92516,-1.5141485) -- (-0.541995,-1.980592) -- cycle;
\fill[line width=0.5pt,fill=black,fill opacity=0.1] (-1.92516,-1.1486) -- (0.5419949999999997,-0.6821565000000002) -- (0.5419950000000003,1.9805925) -- (-1.92516,1.514149) -- cycle;

\fill[line width=0.5pt,fill=black,fill opacity=0] (-1.92516,-1.1486) -- (-0.541995,-1.980592) -- (-0.541995,0.6821570000000001) -- (-1.92516,1.514149) -- cycle;
\fill[line width=0.5pt,fill=black,fill opacity=0] (0.5419950000000003,1.9805925) -- (1.92516,1.1486005000000001) -- (1.92516,-1.5141485) -- (0.5419949999999997,-0.6821565000000002) -- cycle;

\fill[line width=0.5pt,fill=black,fill opacity=0] (-1.92516,1.514149) -- (0.5419950000000003,1.9805925) -- (1.92516,1.1486005000000001) -- (-0.541995,0.6821570000000001) -- cycle;
\fill[line width=0.5pt,fill=black,fill opacity=0] (-0.541995,-1.980592) -- (1.92516,-1.5141485) -- (0.5419949999999997,-0.6821565000000002) -- (-1.92516,-1.1486) -- cycle;

\draw [line width=0.5pt] (-1.92516,1.514149)-- (-1.08399,1.008174144765809);
\draw [line width=0.5pt] (-1.08399,1.008174144765809)-- (-0.541995,0.6821570000000001);
\draw [line width=0.5pt] (0.5419950000000003,1.9805925)-- (0.9326826447970984,1.7455887260684537);
\draw [line width=0.5pt] (0.9326826447970984,1.7455887260684537)-- (1.92516,1.1486005000000001);
\draw [line width=0.5pt] (-1.92516,-1.1486)-- (-1.08399,-1.6545748552341912);
\draw [line width=0.5pt] (-1.08399,-1.6545748552341912)-- (-0.541995,-1.980592);
\draw [line width=0.5pt] (-1.92516,1.514149)-- (-1.4746776447970984,1.5993177739315465);
\draw [line width=0.5pt] (-1.4746776447970984,1.5993177739315465)-- (0.5419950000000003,1.9805925);
\draw [line width=0.5pt, dotted] (-1.92516,-1.1486)-- (0.5419949999999997,-0.6821565000000002);
\draw [line width=0.5pt] (-1.92516,1.514149)-- (-1.92516,-1.1486);
\draw [line width=0.5pt] (1.92516,1.1486005000000001)-- (1.92516,-1.5141485);
\draw [line width=0.5pt, dotted] (1.92516,-1.5141485)-- (0.5419949999999997,-0.6821565000000002);
\draw [line width=0.5pt, dotted] (0.5419949999999997,-0.6821565000000002)-- (0.5419950000000003,1.9805925);
\draw [line width=0.5pt, dotted] (0.5419949999999997,-0.6821565000000002)-- (0.5419950000000003,1.9805925);
\draw [line width=0.5pt] (-0.541995,0.6821570000000001)-- (1.92516,1.1486005000000001);
\draw [line width=0.5pt] (1.92516,1.1486005000000001)-- (1.92516,-1.5141485);
\draw [line width=0.5pt] (1.92516,-1.5141485)-- (-0.541995,-1.980592);
\draw [line width=0.5pt] (-0.541995,-1.980592)-- (-0.541995,0.6821570000000001);
\begin{tiny}
\draw [fill=black] (-1.92516,-1.1486) circle (2.5pt);
\draw[color=black] (-2.7,-1.3) node {$t_{010}$};
\draw [fill=black] (-0.541995,-1.980592) circle (2.5pt);
\draw[color=black] (-1.2,-2.2) node {$t_{011}$};
\draw [fill=black] (1.92516,-1.5141485) circle (2.5pt);
\draw[color=black] (2.7,-1.4) node {$t_{111}$};
\draw [fill=black] (-1.92516,1.514149) circle (2.0pt);
\draw[color=black] (-2.7,1.4) node {$t_{000}$};
\draw [fill=black] (-0.541995,0.6821570000000001) circle (2.0pt);
\draw[color=black] (-1.2,0.4) node {$t_{001}$};
\draw [fill=black] (1.92516,1.1486005000000001) circle (2.0pt);
\draw[color=black] (2.7,1.3) node {$t_{101}$};
\draw [fill=black] (0.5419950000000003,1.9805925) circle (2.0pt);
\draw[color=black] (1.3,2.1) node {$t_{100}$};
\draw [fill=black] (0.5419949999999997,-0.6821565000000002) circle (2.0pt);
\draw[color=black] (1.3,-0.5) node {$t_{110}$};

\draw[color=black] (-1.3,-0.5) node {$t_3^{(1)}$};
\draw[color=black] (0.1,-1.3) node {$t_3^{(0)}$};
\end{tiny}
\end{tikzpicture}
\end{minipage}
\caption{The slices $t_j^{(0)}$ and $t_j^{(1)}$ appearing in the computation of $f_{3,\{j\}}$.}\label{fig: case d=3 factors constant term}
\end{center}
\end{figure}

The new system of coordinates is more effective for computing $\mathrm{SO}(V)^d$-invariants. Indeed, the torus $\mathrm{SO}(V)^d\cong(\C^*)^d=(\C\setminus\{0\})^d$ acts on $V^{\otimes d}$ by rescaling each coordinate $u_{i_1\cdots i_d}$ to $\prod_{j=1}^d\xi_j^{(-1)^{i_j}}u_{i_1\cdots i_d}$ for some $(\xi_1,\ldots,\xi_d)\in(\C^*)^d$. Using \cite[Algorithm 1.4.5]{Stu}, we computed a minimal generating set of invariants of $S(V^{\otimes d})^{\mathrm{SO}(V)^d}$, a least for small values of $d$. Focusing on the case $d=3$, we get that
\[
S(V^{\otimes d})^{\mathrm{SO}(V)^d}\cong\C[\{u_{i_1\cdots i_d}\}]^{\mathrm{SO}(V)^d}\cong\C[\theta_1,\theta_2,\theta_3,\theta_4,\varphi_1,\varphi_2],
\]
where the $\theta_j$'s are four real invariants of degree two, whereas $\varphi_1$ and $\varphi_2$ are two non-real mutually conjugate invariants of degree four:
{\smaller
\begin{align}\label{eq: generators invariant ring}
\begin{split}
\theta_1=u_{0,0,0}u_{1,1,1},\quad \theta_2=u_{0,0,1}u_{1,1,0},&\quad \theta_3=u_{0,1,0}u_{1,0,1},\quad \theta_4=u_{0,1,1}u_{1,0,0},\\
\varphi_1=u_{0,0,1}u_{0,1,0}u_{1,0,0}u_{1,1,1},&\quad \varphi_2=u_{0,0,0}u_{0,1,1}u_{1,0,1}u_{1,1,0}.
\end{split}
\end{align}
}
In addition, the only relation among them is $\theta_1\theta_2\theta_3\theta_4-\varphi_1\varphi_2=0$. Since we are dealing with real binary tensors, the coefficients of the ED polynomial of $X_\mu^\vee$ at $t$ are all real polynomials in the entries $\{t_{i_1\cdots i_d}\}$ of $t$. Indeed, they are elements of $\R[\theta_1,\theta_2,\theta_3,\theta_4,\varphi]$, where $\varphi\coloneqq(\varphi_1+\varphi_2)/2$. In the old set of coordinates $\{t_{i_1\cdots i_d}\}$, these invariants become respectively
{\smaller
\begin{align}\label{eq: invariants p_j q_j 222 case}
\begin{split}
&\theta_1 = (t_{000}-t_{011}-t_{101}-t_{110})^2+(t_{111}-t_{100}-t_{010}-t_{001})^2,\\
&\theta_2 = (t_{000}+t_{011}+t_{101}-t_{110})^2+(t_{111}+t_{100}+t_{010}-t_{001})^2,\\
&\theta_3 = (t_{000}+t_{011}-t_{101}+t_{110})^2+(t_{111}+t_{100}-t_{010}+t_{001})^2,\\
&\theta_4 = (t_{000}-t_{011}+t_{101}+t_{110})^2+(t_{111}-t_{100}+t_{010}+t_{001})^2,\\
&\varphi = t_{000}^{4}+2t_{000}^{2}t_{001}^{2}+t_{001}^{4}+2t_{000}^{2}t_{010}^{2}-2t_{001}^{2}t_{010}^{2}+t_{010}^{4}+8t_{000}t_{001}t_{010}t_{011}-2t_{000}^{2}t_{011}^{2}+2t_{001}^{2}t_{011}^{2}+2t_{010}^{2}t_{011}^{2}\\
&\ +t_{011}^{4}+2t_{000}^{2}t_{100}^{2}-2t_{001}^{2}t_{100}^{2}-2t_{010}^{2}t_{100}^{2}-6t_{011}^{2}t_{100}^{2}+t_{100}^{4}+8t_{000}t_{001}t_{100}t_{101}+8t_{010}t_{011}t_{100}t_{101}-2t_{000}^{2}t_{101}^{2}\\
&\ +2t_{001}^{2}t_{101}^{2}-6t_{010}^{2}t_{101}^{2}-2t_{011}^{2}t_{101}^{2}+2t_{100}^{2}t_{101}^{2}+t_{101}^{4}+8t_{000}t_{010}t_{100}t_{110}+8t_{001}t_{011}t_{100}t_{110}+8t_{001}t_{010}t_{101}t_{110}\\
&\ -8t_{000}t_{011}t_{101}t_{110}-2t_{000}^{2}t_{110}^{2}-6t_{001}^{2}t_{110}^{2}+2t_{010}^{2}t_{110}^{2}-2t_{011}^{2}t_{110}^{2}+2t_{100}^{2}t_{110}^{2}-2t_{101}^{2}t_{110}^{2}+t_{110}^{4}-8t_{001}t_{010}t_{100}t_{111}\\
&\ +8t_{000}t_{011}t_{100}t_{111}+8t_{000}t_{010}t_{101}t_{111}+8t_{001}t_{011}t_{101}t_{111}+8t_{000}t_{001}t_{110}t_{111}+8t_{010}t_{011}t_{110}t_{111}+8t_{100}t_{101}t_{110}t_{111}\\
&\ -6t_{000}^{2}t_{111}^{2}-2t_{001}^{2}t_{111}^{2}-2t_{010}^{2}t_{111}^{2}+2t_{011}^{2}t_{111}^{2}-2t_{100}^{2}t_{111}^{2}+2t_{101}^{2}t_{111}^{2}+2t_{110}^{2}t_{111}^{2}+t_{111}^{4}.
\end{split}
\end{align}
}
Now we describe the ED polynomial of $X_3^\vee\subset\PP(V^{\otimes 3})\cong\PP^7$ at $t\in V^{\otimes 3}$, which has degree $6=\mathrm{EDdegree}(X_3)$ in $\epsilon^2$ and may be written as
\[
\mathrm{EDpoly}_{X^\vee,t}(\epsilon^2)=a_6(t)\epsilon^{12}+a_5(t)\epsilon^{10}+\cdots+a_0(t).
\]
From Corollary \ref{cor: Main Theorem, Segre case}, we obtain that the extreme coefficients are respectively
\begin{equation}\label{eq: factors constterm 222 case}
a_0=g_0^2\cdot g_1=\mathrm{Det}^2\cdot f_{3,\{1\}}\cdot f_{3,\{2\}}\cdot f_{3,\{3\}},\quad a_6=g_3=f_{3,[3]},
\end{equation}
where $\mathrm{Det}(t)$ is the hyperdeterminant of $t$, written explicitly as
\begin{align*}
\mathrm{Det} &= \frac{1}{64}\left[2(\theta_1\theta_2+\theta_1\theta_3+\theta_2\theta_3+\theta_1\theta_4+\theta_2\theta_4+\theta_3\theta_4)-(\theta_1^2+\theta_2^2+\theta_3^2+\theta_4^{2})-8\varphi\right]\\
 &=\left[\det
 \begin{pmatrix}
 t_{000} & t_{011}\\
 t_{100} & t_{111}
 \end{pmatrix}
 +\det
 \begin{pmatrix}
 t_{010} & t_{001}\\
 t_{110} & t_{101}
 \end{pmatrix}
 \right]^2-
 4\det
 \begin{pmatrix}
 t_{000} & t_{001}\\
 t_{100} & t_{101}
 \end{pmatrix}
 \det
 \begin{pmatrix}
 t_{010} & t_{011}\\
 t_{110} & t_{111}
 \end{pmatrix}.
\end{align*}
The other factors of $a_0$ and $a_6$ are written explicitly below with respect to both the generators $\theta_1,\ldots,\theta_4,\varphi$ and the coordinates of $t$, in the same fashion of Remark \ref{rmk: sos constterm pictures}.
\begin{align}\label{eq: three additional factors constterm 222 case}
\begin{split}
f_{3,\{1\}} &= \frac{1}{16}\left[\theta_1\theta_2+\theta_3\theta_4-2\varphi\right]=\det\left[\left(t_1^{(0)}+\sqrt{-1}t_1^{(1)}\right)\left(t_1^{(0)}-\sqrt{-1}t_1^{(1)}\right)\right],\\
f_{3,\{2\}} &= \frac{1}{16}\left[\theta_2\theta_3+\theta_1\theta_4-2\varphi\right]=\det\left[\left(t_2^{(0)}+\sqrt{-1}t_2^{(1)}\right)\left(t_2^{(0)}-\sqrt{-1}t_2^{(1)}\right)\right],\\
f_{3,\{3\}} &= \frac{1}{16}\left[\theta_1\theta_3+\theta_2\theta_4-2\varphi\right]=\det\left[\left(t_3^{(0)}+\sqrt{-1}t_3^{(1)}\right)\left(t_3^{(0)}-\sqrt{-1}t_3^{(1)}\right)\right].
\end{split}
\end{align}
On one hand, the three factors $f_{3,\{1\}}$, $f_{3,\{2\}}$ and $f_{3,\{3\}}$ of $g_1$ represent three quartic hypersurfaces in $\PP^7$. Each of them is the union of two conjugate quadric hypersurfaces. In turn, the singular locus of each of these quadric hypersurfaces has dimension three and meets the Segre variety $X_3$ in a quadric surface. Finally, these six quadric surfaces may be interpreted as the two dimensional ``faces'' of the three dimensional ``cube'' $X_{3,[3]}=\mathrm{Seg}(Q^{\times 3})$ of totally isotropic rank one tensors, as Figure \ref{fig: case d=3 factors constant term} may suggest. On the other hand, $a_6=\theta_1\theta_2\theta_3\theta_4$. Geometrically speaking, each polynomial $\theta_j$ defines a pair of conjugate hyperplanes dual to a pair of conjugate vertices of $X_{3,[3]}$.

In addition, we determined symbolically all the intermediate coefficients of $\mathrm{EDpoly}_{X^\vee,t}(\epsilon^2)$ with respect to the generators $\theta_1,\ldots,\theta_4,\varphi$. In particular, $\deg(a_j)=2(10-j)$ for all $j\in\{0,\ldots,6\}$. For example, the coefficient $a_5(t)$ is relevant since the ratio $a_5(t)/a_6(t)$ corresponds to the sum of the squares of the singular values of $t$, thanks to Proposition \ref{pro: roots and singular values}:
\[
a_5(t) = \frac{1}{8}\left[(\theta_1\theta_2\theta_3+\theta_1\theta_2\theta_4+\theta_1\theta_3\theta_4+\theta_2\theta_3\theta_4)\varphi-3\theta_1\theta_2\theta_3\theta_4(\theta_1+\theta_2+\theta_3+\theta_4)\right].
\]
In the following, we assume that $t\in V^{\otimes 3}$ is $\mu$-symmetric for $\mu\in\{(2,1),(3)\}$. Among the six critical binary tensors for $t$ on $X_3$, $\mathrm{EDdegree}(X_\mu)$ of them are $\mu$-symmetric. Below we describe the critical binary tensors that belong to $X_3\setminus X_\mu$.
\begin{proposition}\label{pro: reality issue}
\begin{enumerate}
	\item Let $\mu=(2,1)$ and let $t\in S^\mu V$ be general. Then $t$ admits four critical binary tensors on $X_\mu$. The remaining two critical binary tensors are $x\otimes y\otimes z$ and $y\otimes x\otimes z$ for some $x,y,z\in V$. If $t$ is real, the common singular value of the two critical points on $X\setminus X_\mu$ is real.
	\item Let $\mu=(3)$ and let $t\in S^\mu V$ be general. Then $t$ admits three critical binary tensors on $X_\mu$. The remaining three critical binary tensors are $x\otimes x\otimes y$, $x\otimes y\otimes x$ and $y\otimes x\otimes x$ for some $x,y\in V$. If $t$ is real, the common singular value of the three critical points on $X\setminus X_\mu$ is real.
\end{enumerate}
\end{proposition}
\proof
By Proposition \ref{pro: ed poly divides ed poly}, $\mathrm{EDpoly}_{X_\mu^\vee,t}(\epsilon^2)$ divides $\mathrm{EDpoly}_{X_3^\vee,t}(\epsilon^2)$ with multiplicity one when $t\in S^\mu V$.
Consider part $(1)$. Then $t$ admits $\mathrm{EDdegree}(X_\mu)=4$ (see (\ref{eq: EDdegree Segre-Veronese})) critical binary tensors corresponding to four singular vector triples $(x_j,x_j,y_j)$ for some $x_j,y_j\in V$, $j\in[4]$. Moreover, for any singular vector triple $(x,y,z)$ for $t$ with singular value $\sigma$ and $x\neq y$, the permutation $(y,x,z)$ is again a singular vector triple for $t$, and shares the same singular value $\sigma$. Hence, there is a linear polynomial $h(\epsilon^2)$ such that
\[
\mathrm{EDpoly}_{X_3^\vee,t}(\epsilon^2)=\mathrm{EDpoly}_{X_\mu^\vee,t}(\epsilon^2)\cdot h(\epsilon^2)^2.
\]
In conclusion, apart from the $\mu$-symmetric singular vector tuples, there is room left only for one more non-symmetric singular vector triple $(x,y,z)$ and its permutation $(y,x,z)$. In addition, if $t\in S^\mu V_\R$, the root of the linear polynomial $h(\epsilon^2)$ must be real.

Now consider part $(2)$. Then $t$ admits $\mathrm{EDdegree}(X_{(3)})=3$ (see (\ref{eq: EDdegree Segre-Veronese})) critical binary tensors corresponding to three singular vector triples $(x_j,x_j,x_j)$ for some $x_j\in V$, $j\in[3]$. With a similar argument of part $(1)$, we observe that there is a linear polynomial $\tilde{h}(\epsilon^2)$ such that
\[
\mathrm{EDpoly}_{X_3^\vee,t}(\epsilon^2)=\mathrm{EDpoly}_{X_\mu^\vee,t}(\epsilon^2)\cdot h'(\epsilon^2)^3.
\]
In conclusion, apart from the $\mu$-symmetric singular vector tuples, there is room left only for one more non-symmetric singular vector triple of the form $(x,x,y)$ for some $x,y\in V$, together with its permutations $(x,y,x)$ and $(y,x,x)$. Moreover, if $t$ has real entries, the root of $h'(\epsilon^2)$ must be real.\qedhere
\endproof

\begin{remark}
Consider Proposition \ref{pro: reality issue}(1). In this case the invariants $\theta_2$ and $\theta_3$ introduced in (\ref{eq: invariants p_j q_j 222 case}) coincide. This implies that the highest coefficient $a_6=\theta_1\theta_2\theta_3\theta_4$ of $\mathrm{EDpoly}_{X_3^\vee,t}(\epsilon^2)$ splits into two factors $\theta_1\theta_4$ and $\theta_2\theta_3=\theta_2^2$, which correspond to the highest coefficients of $\mathrm{EDpoly}_{X_\mu^\vee,t}(\epsilon^2)$ and $h(\epsilon^2)^2$, respectively. About the lowest coefficient $a_0$, from (\ref{eq: three additional factors constterm 222 case}) we see that in this case the polynomials $f_{3,\{1\}}$ and $f_{3,\{3\}}$ coincide. Indeed, the lowest coefficients of $\mathrm{EDpoly}_{X_\mu^\vee,t}(\epsilon^2)$ and $h(\epsilon^2)^2$ are respectively $\mathrm{Det}^2\cdot f_{3,\{2\}}$ and $f_{3,\{1\}}\cdot f_{3,\{3\}}=f_{3,\{1\}}^2$. More precisely, $\mathrm{Det}=f_\mu$ and $f_{3,\{2\}}=f_{\mu,\{2\}}$. We computed symbolically the ED polynomial of $X_3^\vee$ at a $\mu$-symmetric tensor $t$. In particular,
\[
h(\epsilon^2)=16\theta_2\epsilon^2-\theta_1\theta_2-\theta_2\theta_4+2\varphi.
\]
In addition, a consequence of Proposition \ref{SOS} is that, up to sign multiplication, the highest and lowest coefficients of $h(\epsilon^2)$ are SOS polynomials. In particular, the root of $h(\epsilon^2)$ may be written as
\[
\hspace{1.3cm}16\epsilon^2=\frac{\theta_1\theta_2+\theta_2\theta_4-2\varphi}{\theta_2}=\frac{(c_{01}c_{10}-c_{11}c_{20}-c_{00}c_{11}+c_{10}c_{21})^2+(c_{00}c_{21}-c_{01}c_{20})^2}{(c_{00}+c_{20})^2+(c_{01}+c_{21})^2},
\]
where we are using the $\mu$-symmetric variables $\{c_{ij}\}$ introduced at the beginning of Section \ref{sec: computation}.

Now consider Proposition \ref{pro: reality issue}(2). Looking at their definition in (\ref{eq: invariants p_j q_j 222 case}), in this case the invariants $\theta_2$, $\theta_3$ and $\theta_4$ coincide. Indeed the highest coefficient $a_6=\theta_1\theta_2\theta_3\theta_4$ of $\mathrm{EDpoly}_{X_3^\vee,t}(\epsilon^2)$ splits into two factors $\theta_1$ and $\theta_2\theta_3\theta_4=\theta_2^3$, which correspond to the highest coefficients of $\mathrm{EDpoly}_{X_\mu^\vee,t}(\epsilon^2)$ and $h'(\epsilon^2)^3$, respectively. About the lowest coefficient $a_0$, from (\ref{eq: three additional factors constterm 222 case}) we see that in this case the polynomials $f_{3,\{1\}}$, $f_{3,\{2\}}$ and $f_{3,\{3\}}$ coincide. Indeed, the lowest coefficients of $\mathrm{EDpoly}_{X_\mu^\vee,t}(\epsilon^2)$ and $h'(\epsilon^2)^3$ are respectively $\mathrm{Det}^2$ and $f_{3,\{1\}}\cdot f_{3,\{2\}}\cdot f_{3,\{3\}}=f_{3,\{1\}}^3$. More precisely, $\mathrm{Det}=f_\mu$. Moreover, in this case
\[
h'(\epsilon^2)=16\theta_2\epsilon^2-\theta_1\theta_2-\theta_2^2+2\varphi
\]
and the root of $h'(\epsilon^2)$ may be expressed as (using the coordinates $\{c_j\}$ of the symmetric tensor $t$)
\[
\hspace{1.3cm}16\epsilon^2=\frac{\theta_1\theta_2+\theta_2^2-2\varphi}{\theta_2}=\frac{(c_1^2-c_2^2-c_0c_2+c_1c_3)^2+(c_0c_3-c_1c_2)^2}{(c_0+c_2)^2+(c_1+c_3)^2}.
\]
\end{remark}

\begin{remark}
More generally, one may verify that for any partition $\mu\vdash d$ and for a general symmetric binary tensor $t\in S^dV$, the polynomial $\mathrm{EDpoly}_{X_\mu^\vee,t}(\epsilon^2)$ is divided by $\mathrm{EDpoly}_{X_{(d)}^\vee,t}(\epsilon^2)$ and by other factors. We observed that there is a precise relation between the factors of $\mathrm{EDpoly}_{X_\mu^\vee,t}(\epsilon^2)$ and the dual multiple root loci $[\mathrm{Chow}_{\lambda}(\PP(V))]^\vee$ for all $\lambda\prec\mu$, that somehow shifts the work by Oeding \cite{O} from symmetrizations of $\mu$-discriminants to symmetrizations of their respective ED polynomials. This and other related aspects are studied in the paper in preparation \cite{Sod2}.
\end{remark}

\section*{Acknowledgements}
The author is very grateful to his advisor Giorgio Ottaviani for valuable guidance. Moreover, he warmly thanks Fr\'ed\'eric Holweck, Antonio Lerario and Luke Oeding for their remarks. Luca Sodomaco is member of INDAM-GNSAGA.


\begin{thebibliography}{10}

\bibitem[CS]{CS} D. Cartwright, B. Sturmfels, The number of eigenvalues of a tensor, Linear Algebra Appl. {\bf 438} (2013), no. 2, 942--952.

\bibitem[DHOST]{DHOST} J. Draisma, E. Horobe\c{t}, G. Ottaviani,  B. Sturmfels, R. Thomas, The Euclidean Distance Degree of an algebraic variety, Found. Comput. Math. {\bf 16} (2016), no. 1, 99--149.

\bibitem[FO]{FO} S. Friedland, G. Ottaviani, The number of singular vector tuples and uniqueness of best rank one approximation of tensors, Found. Comput. Math. {\bf 14} (2014), 1209--1242.

\bibitem[GKZ]{GKZ} I. M. Gelfand, M. M. Kapranov, A. V. Zelevinsky, {\em Discriminants, Resultants and Multidimensional Determinants}, Birkh\"auser, Boston, 1994.

\bibitem[GS]{GS} D. Grayson, M. Stillman, {\em Macaulay2, a software system for research in algebraic geometry}. Available at \url{http://www.math.uiuc.edu/Macaulay2/}.

\bibitem[HHLQ]{HHLQ} S. Hu, Z. H. Huang, C. Ling, L. Qi, On determinants and eigenvalue theory of tensors, J. Symbolic
Comput. {\bf 50} (2013), 508--531.

\bibitem[HO]{HO} F. Holweck, L. Oeding, Hyperdeterminants from the $E_8$ discriminant, \arxiv{1810.05857}.

\bibitem[HW]{HW} E. Horobe\c{t}, M. Weinstein, Offset hypersurfaces and persistent homology of algebraic varieties, \arxiv{1803.07281}.

\bibitem[Lan]{Lan} J. M. Landsberg, {\em Tensors: Geometry and Applications}, Graduate studies in mathematics, American Mathematical
Society, Providence, 2011.
 
\bibitem[Lim]{Lim} L. H. Lim, Singular values and eigenvalues of tensors: a variational approach, Proc. IEEE Internat. Workshop on Comput. Advances in Multi-Sensor Adaptive Processing (CAMSAP 2005), 129--132.

\bibitem[LQZ]{LQZ} A. M. Li, L. Qi, B. Zhang,
E-characteristic polynomials of tensors, Commun. Math. Sci. {\bf 11} (2013), no. 1, 33--53.

\bibitem[NQWW]{NQWW} G. Ni, L. Qi, F. Wang, Y. Wang, The degree of the E-characteristic polynomial of an even order tensor, J. Math. Anal. Appl. {\bf 329} (2007), no. 2, 1218--1229.

\bibitem[O]{O} L. Oeding, Hyperdeterminants of polynomials, Adv. Math. {\bf 231} (2012), no. 3-4, 1308--1326.

\bibitem[OS]{OS} G. Ottaviani, L. Sodomaco, The distance function from a real algebraic variety, \arxiv{1807.10390}.

\bibitem[Q]{Q} L. Qi, Eigenvalues of a real supersymmetric tensor, J. Symbolic Comput. {\bf 40} (2005), no. 6, 1302--1324.

\bibitem[Q2]{Q2} L. Qi, Eigenvalues and invariants of tensors, J. Math. Anal. Appl. {\bf 325} (2007), no. 2, pp. 1363--1377.

\bibitem[QL]{QL} L. Qi, Z. Luo, {\em Tensor analysis: Spectral theory and special tensors}, SIAM, Philadelphia, 2017.

\bibitem[Sod]{Sod} L. Sodomaco, The product of the eigenvalues of a symmetric tensor, Linear Algebra Appl. {\bf 554} (2018), 224--248.

\bibitem[Sod2]{Sod2} L. Sodomaco, On the ED polynomial of a Segre-Veronese variety, in preparation.

\bibitem[Stu]{Stu} B. Sturmfels, {\em Algorithms in invariant theory}, Springer Science \& Business Media, 2008.

\end{thebibliography}
\end{document}